\documentclass[preprint,10pt]{elsart}
\usepackage{graphicx}
\usepackage{latexsym}
\usepackage{amssymb}
\usepackage{subfigure}
\usepackage[usenames,dvipsnames]{pstricks}
\usepackage{epsfig}

\newcommand{\PM}{Petviashvili }

\usepackage{graphicx,epsfig,amssymb,amsmath,amsfonts}

\begin{document}
\begin{frontmatter}
\title{On Petviashvili type methods for traveling wave computations:  Acceleration techniques}
\author{J. \'Alvarez}
\address{Department of Applied Mathematics,
University of Valladolid, Paseo del Cauce 59, 47011,
Valladolid, Spain.}
\address{
IMUVA, Institute of Mathematics of University of Valladolid; Spain.
Email: joralv@eii.uva.es}
\author{A. Dur\'an \thanksref{au}}
\address{Department of Applied Mathematics, University of
Valladolid, Paseo de Belen 15, 47011-Valladolid, Spain.}
\address{
IMUVA, Institute of Mathematics of University of Valladolid; Spain.
Email:
angel@mac.uva.es }

\thanks[au]{Corresponding author}

\begin{abstract}
In this paper a family of fixed point algorithms, generalizing the \PM method, is considered. A previous work studied the convergence of the methods. Presented here is a second part of the analysis, concerning the introduction of some acceleration techniques into the iterative procedures. The purpose of the research is two-fold: one is improving the performance of the methods in case of convergence and the second one is widening their application when generating traveling waves in nonlinear dispersive wave equations, transforming some divergent into convergent cases. Two families of acceleration techniques are considered: the vector extrapolation methods and the Anderson acceleration methods. A comparative study through several numerical experiments is carried out.
\end{abstract}
\end{frontmatter}

\section{Introduction}
\label{sec1}
In a previous paper \cite{alvarezd}, a family of fixed-point algorithms for the numerical approximation of nonlinear systems of the form
\begin{eqnarray}
Lu=N(u),\quad u\in \mathbb{R}^{m}, \quad m>1,\label{mm1}
\end{eqnarray}
was introduced. In (\ref{mm1}),  $L$ is a nonsingular $m\times m$ real matrix and
$N:\mathbb{R}^{m}\rightarrow \mathbb{R}^{m}$ is an homogeneous
function with degree $p, |p|>1$ (this means that $N(\lambda u)=\lambda^{p}N(u)$). Among other applications, systems of this form are very typical in the numerical generation of traveling waves and ground states in water wave problems and nonlinear optics. For the numerical approximation to solutions of (\ref{mm1}), the use of the classical fixed-point algorithm is not suitable. This is due to the fact that if $u^{*}$ is a solution and $S=L^{-1}N^{\prime}(u^{*})$ stands for the iteration matrix at $u^{*}$, then, since $N$ is homogeneous of degree $p$ then $N^{\prime}(u^{*})u^{*}=pN(u^{*})$ and therefore
\begin{eqnarray*}
Su^{*}=L^{-1}N^{\prime}(u^{*})u^{*}=pL^{-1}N(u^{*})=pu^{*},
\end{eqnarray*}
that is, $p$ is an eigenvalue of $S$ with magnitude above one. This makes the iteration not convergent in general.

As an alternative and based on the \PM method, \cite{petviashvili,pelinovskys,lakobay,lakobay2},  the following fixed-point algorithms were considered in \cite{alvarezd}:
\begin{eqnarray}
 Lu_{n+1}=s(u_{n})N(u_{n}), \quad n=0,1,\ldots,\label{mm2}
\end{eqnarray}
 from $u_{0}\neq 0$ and where $s:\mathbb{R}^{m}\rightarrow \mathbb{R}$ is a $C^{1}$ function satisfying the following properties:
\begin{itemize}
\item[(P1)] A set of fixed points of the iteration operator
\begin{eqnarray}
F(x)=s(x)L^{-1}N(x),\label{iterop}
\end{eqnarray}
coincides with a set of fixed points of (\ref{mm1}). This means that: (a) if $u^{*}$ is a solution of (\ref{mm1}) then $s(u^{*})=1$; (b) inversely, if the sequence $\{u_{n}\}_{n}$, generated by (\ref{mm2}), converges to some $y$, then $s(y)=1$ (and, consequently, $y$ is a solution of (\ref{mm1})).
\item[(P2)] $s$ is homogeneous with degree $q$ such that $|p+q|<1$.
\end{itemize}
The function $s$ is called the stabilizing factor of the method, inheriting the nomenclature of the \PM method. Actually, formula (\ref{mm2}) generalizes the \PM scheme, which corresponds to the choice
\begin{eqnarray}
s(x)=\left(\frac{\langle Lx,x\rangle}{\langle N(x),x\rangle}\right)^{\gamma},\quad q=\gamma(1-p),\label{mm3c}
\end{eqnarray}
The first part of the work, carried out in \cite{alvarezd} (see also \cite{alvarezd2015}), analyzed the convergence of (\ref{mm2}). The main conclusion was that, compared to the classical fixed-point algorithm, the stabilizing factor acts like a filter of the spectrum of the matrix $S$ in the sense that:
\begin{itemize}
\item The eigenvalue $\lambda=p$ of $S$ is transformed to the eigenvalue $\lambda=p+q$ of the iteration matrix $F^{\prime}(u^{*})$ of (\ref{iterop}).
\item The rest of the spectrum of $F^{\prime}(u^{*})$ is included into the spectrum of $S$.
\end{itemize}
Thus, the convergence of (\ref{mm2}) depends of the spectrum of $S$ different from $p$. From these conclusions, several results of convergence can be derived, see \cite{alvarezd} for details.

The motivation of this paper is two-fold. First, several  numerical experiments in the literature show that the \PM type algorithms are in sometimes computationally slower than other alternatives. In order to continue to benefit from the easy implementation (one of the advantages of the methods) the algorithms should improve their performance with the inclusion of some acceleration technique. A further motivation comes from the known mechanism of some extrapolation methods, \cite{Sidi2003}, to transform divergent into convergent cases. The application of this property to these \PM type methods may extend their use to compute traveling waves under more demanding conditions, for example in two dimensions or/and in case of highly oscillatory waves.

The literature on acceleration techniques is very rich with many different families and strategies, \cite{Sidi2003,BrezinskiR1991}. This paper will be focused on two types of procedures: the vector extrapolation methods, \cite{BrezinskiR1991,smithfs,jbilous,brezinski2,sidifs} and the Anderson mixing, \cite{anderson,walkern,ni,fangs}. We think that the first one is the most widely studied group; in particular, known convergence results for some of these methods will serve us to justify several examples of transformation from divergence to convergence when generating traveling waves iteratively. The second family accelerates the convergence by introducing the strategy of minimization of the residual in some norms at each step. They have been revealed efficient in, for example, electronic structure computations, \cite{anderson,Pulay1982} (see also \cite{ni,walkern} and references therein) and, to our knowledge, this is the first time they are applied to the numerical generation of traveling waves.

The main purpose of this paper is then exploring by numerical means the application of these two acceleration methods to the generation of traveling waves through several problems of interest and from the \PM type methods (\ref{mm2}) and their extended versions derived in \cite{alvarezd2014b}. With the case studies presented here we have tried to cover different situations of hard computation of the waves as an attempt to give some guidelines of application. In this sense, the paper provides several conclusions to be emphasized:
\begin{itemize}
\item The use of acceleration techniques is highly recommended here since it improves the performance in general and allows to extend the application of the methods to computationally harder situations, with especial emphasis on two-dimensional simulations and highly oscillatory wave generation.
\item By comparing the two families of acceleration techniques considered in this study, the vector extrapolation methods are in general more competitive for these problems compared to the Anderson acceleration methods. Among the vector extrapolation methods, the polynomial methods provide a better performance in general (some exceptions can be seen in the experiments below).
\item The main drawback of the Anderson acceleration methods  concerns the numerical treatment of the associated minimization problem, since most of the difficulties come from ill-conditioning. This might be improved by including suitable preconditioning techniques (here the methods were implemented in a standard way, \cite{walkern,ni,fangs}). However, it is remarkable that when the Anderson acceleration methods work, their performance is in general comparable to that of some vector extrapolation methods.
\end{itemize}

The structure of the paper is as follows. Section \ref{se2} is devoted to a description of the two families of acceleration techniques considered in this study. This also includes some comments on the implementation and convergence results. The application of both techniques to the methods (\ref{mm2}) is studied in Section \ref{se3} through a plethora of numerical experiments involving the computation of different types of traveling waves: ground states, classical and generalized solitary waves as well as periodic traveling waves. The numerical study will be focused on the two main motivations of the paper: the improvement of the efficiency and the extension of application of the methods to computationally harder problems and where the iteration is initially not convergent. Finally, Section \ref{se4} completes the computational study with some illustrations of the application of the acceleration to the extended versions of the algorithms (\ref{mm2}), treated in \cite{alvarezd2014b} and suitable when the nonlinearity in (\ref{mm1}) contains several homogeneous terms with different degree. Some concluding remarks are in Section \ref{se5}.

\section{Acceleration techniques}
\label{se2}
Besides the local character of the convergence, in some cases and compared to other alternatives, fixed point algorithms has the additional disadvantage of a slow performance. In what follows, several techniques of acceleration will be considered and applied to the methods (\ref{mm2}), with the aim of improving their efficiency. Furthermore, as in the case of the classical algorithm, \cite{smithfs}, some cases of divergence will be transformed to convergent iterations.

This section introduces two families of acceleration techniques: the vector extrapolation methods (VEM from now on) and the Anderson acceleration methods (AAM). We will include a description of the schemes (including some convergence results) and some comments on implementation.

\subsection{Vector extrapolation methods}
The first group of acceleration techniques consists of vector extrapolation methods (VEM). For a more detailed analysis and implementation of the methods see \cite{CabayJ1976,Eddy1979,Mesina1977,smithfs,jbilous,brezinski2,BrezinskiR1991,Sidi2003} and references therein. Here we will describe the general features of the procedures and their application to (\ref{mm2}).

Two families of VEM are typically emphasized in the literature. The first one covers the so-called polynomial methods; they include, as the most widely cited,  the minimal polynomial extrapolation (MPE), the reduced rank extrapolation (RRE) and the modified minimal polynomial extrapolation (MMPE) methods, \cite{CabayJ1976,Eddy1979,Mesina1977,smithfs,jbilous,sidifs,sidi}. The second family consists of the so-called $\epsilon$-algorithms; typical examples are the scalar and vector $\epsilon$-algorithms and the topological $\epsilon$-algorithm, \cite{brezinski1,brezinski2,jbilous,tan}. All the methods share of course the idea of introducing the extrapolation as a procedure to transform the original sequence $\{u_{n}\}$ of the involved iterative process by some strategy. The polynomial methods are usually described in terms of the transformation ($k\leq m$)
\begin{eqnarray}
&&T_{k}:\mathbb{R}^{m}\longrightarrow\mathbb{R}^{m}\label{fsec211}\\
&&u_{n}\longmapsto t_{n,k}=u_{n}-\Delta U_{n,k}\left(V_{n,k}^{*}\Delta^{2}U_{n,k}\right)^{+}V_{n´k}\Delta u_{n},\nonumber
\end{eqnarray}
where 
\begin{itemize}
\item $\Delta u_{n}=u_{n+1}-u_{n}, \Delta^{2}u_{n}=\Delta u_{n+1}-\Delta u_{n}$.
\item $\Delta^{i}U_{n.k}$ ($i=1,2$) denotes the $m\times k$ matrix of columns $\Delta^{i}u_{n},\ldots,\Delta^{i}u_{n+k-1}$.
\item $V_{n,k}$ stands for the $m\times k$ matrix of some columns $v_{1}^{(n)},\ldots,v_{k}^{(n)}$ with $V_{n,k}^{*}$ as the adjoint matrix of $V_{n,k}$ (conjugate transpose).
\end{itemize}
In (\ref{fsec211}), $A^{+}$ stands for the Moore-Penrose generalized inverse of $A$, defined as $A^{+}=(A^{*}A)^{-1}A^{*}$, \cite{Demmel,Meyer,GolubV}. Different choices of the vectors $v_{j}^{(n)}, j1,\ldots,k$ lead to the most widely used polynomial methods:
\begin{enumerate}
\item[(i)] Minimal polynomial extrapolation (MPE): $v_{j}^{(n)}=\Delta u_{n+j-1}, j=1,\ldots,k$.
\item[(ii)] Reduced rank extrapolation (RRE): $v_{j}^{(n)}=\Delta^{2} u_{n+j-1}, j=1,\ldots,k$. 
\item[(iii)] Modified minimal polynomial extrapolation (MMPE):  $v_{j}^{(n)}=v_{j}, j=1,\ldots,k$, for arbitrary, fixed, linearly independent vectors $v_{1},\ldots v_{k}\in \mathbb{R}^{m}$.
\end{enumerate}

The formulation of the VEM may follow an alternative approach, \cite{sidifs,smithfs}. The transformation (\ref{fsec211}) can be computed in the form
\begin{eqnarray}
t_{n.k}=\sum_{j=0}^{k}\gamma_{j}u_{n+j},\label{fsec214}\quad \sum_{j=0}^{k}\gamma_{j}=1,
\end{eqnarray}
where the coefficients $\gamma_{j}$ are obtained from the resolution (in some sense) of overdetermined, inconsistent systems
\begin{eqnarray}
\sum_{i=0}^{k-1}d_{i}w_{n+i}=\widetilde{w}_{n},\label{fsec215}
\end{eqnarray}
for some vectors $w_{j}, \widetilde{w}_{j}\in\mathbb{R}^{m}$. Different methods emerge by combining different choices of the norm where the residual vector $\sum_{i=0}^{k-1}d_{i}w_{n+1}=\widetilde{w}_{n}$ is minimized with suitable vectors $w_{j}, \widetilde{w}_{j}\in\mathbb{R}^{m}$. Thus, for example, assuming $k<m$, we have:
\begin{itemize}
\item RRE is obtained by writing (\ref{fsec214}) in the form
\begin{eqnarray*}
t_{n.k}=u_{n}-\sum_{j=0}^{k-1}\beta_{j}\Delta u_{n+j},
\end{eqnarray*}
where the $\beta_{j}$ solve (\ref{fsec215}) with $w_{j}=\Delta^{2}u_{j}, \widetilde{w}_{j}=\Delta u_{j}$ and the Euclidean norm with equal weights is used.
\item MPE is obtained by using (\ref{fsec214}) with
\begin{eqnarray*}
\gamma_{j}=\frac{c_{j}}{\sum_{i=0}^{k}c_{i}},\leq 0\leq j\leq k,
\end{eqnarray*}
where $c_{k}=1$ and the $c_{j}, 0\leq j\leq k-1$ solve (\ref{fsec215}) with $w_{j}=\Delta u_{j}, \widetilde{w}_{j}=-\Delta u_{j+k}$ in the sense of minimization with the Euclidean norm with equal weights.
\item MMPE is obtained from (\ref{fsec214}) but where instead of (\ref{fsec215}) a system of the form
\begin{eqnarray}
\sum_{i=0}^{k-1}d_{i}Q_{j}(w_{n+i})=Q_{j}(\widetilde{w}_{n}),\quad j=1,\ldots,k,\label{fsec216}
\end{eqnarray}
is used. In (\ref{fsec216}) $Q_{j}(y)=\langle e_{j},y\rangle=y_{j}$, being $y=(y_{1},\ldots,y_{m})^{T}$.
\end{itemize}
These formulations can be unified by  a  representation with determinants, \cite{smithfs,sidifs,sidi,SidiB1988}. This writes the extrapolation steps  $t_{n,k}$ in the form
\begin{eqnarray}
t_{n,k}=\frac{D(u_{n},u_{n+1},\ldots,u_{n+k})}{D(1,1,\ldots,1)},\label{fsec217}
\end{eqnarray}
with
\begin{eqnarray}
D(\sigma_{0},\ldots,\sigma_{k})=\left|\begin{matrix}
\sigma_{0}&\sigma_{1}&\cdots&\cdots&\sigma_{k}\\
u_{0,0}&u_{0,1}&\cdots&\cdots&u_{0,k}\\
\vdots&\vdots&\cdots&\cdots&\vdots\\
\vdots&\vdots&\cdots&\cdots&\vdots\\
u_{k-1,0}&u_{k-1,1}&\cdots&\cdots&u_{k-1,k}
\end{matrix}\right|,\label{fsec218}
\end{eqnarray}
where the $u_{i,j}$ are scalars that depend on the extrapolation method and where the expansion of (\ref{fsec218}) is in the sense
\begin{eqnarray}
D(\sigma_{0},\ldots,\sigma_{k})=\sum_{i=0}^{k}\sigma_{i}N_{i},\label{fsec219}
\end{eqnarray}
with $N_{i}$ the cofactor of $\sigma_{i}$ in the first row. Thus in the case of the numerator in (\ref{fsec217}), formula (\ref{fsec219}) is a vector, while in the case of the denominator in (\ref{fsec217}), formula (\ref{fsec219}) is a scalar. (See \cite{BrezinskiR2003} for the interpretation in terms of the Schur complement of a matrix.) The three previously mentioned polynomial methods correspond to the following choices of $u_{i,j}$:
\begin{itemize}
\item MPE: $u_{i,j}=\langle \Delta u_{n+i},\Delta u_{n+j}\rangle$.
\item RRE: $u_{i,j}=\langle \Delta^{2} u_{n+i},\Delta u_{n+j}\rangle$.
\item MMPE: $u_{i,j}=\langle e_{i+1},\Delta u_{n+j}\rangle$, where $e_{1},\ldots e_{k}$ are linearly independent vectors in $\mathbb{R}^{m}$.
\end{itemize}

A second family of VEM is called the $\epsilon$-algorithms. A description of them may start from the scalar $\epsilon$-algorithm of Wynn, \cite{Wynn1962,Wynn1966}. This scalar extrapolation method can be derived from the representation (\ref{fsec217}), (\ref{fsec218}) (in the scalar case) with $u_{i,j}=\Delta u_{n+i+j}, i=0,\ldots,k-1; j=0,\ldots,k$ (which are scalars in the scalar case).  The corresponding ratio of determinants
\begin{eqnarray}
t_{n,k}=e_{k}(u_{n})=\frac{D(u_{n},u_{n+1},\ldots,u_{n+k})}{D(1,1,\ldots,1)},\label{fsec2111}
\end{eqnarray}
is called the classical $e$- (or Shanks Schmidt $SS$) transform, \cite{Schmidt1941,Shanks1955,Wynn1956}. This ratio can be evaluated recursively for increasing $k$ and $n$ without the computation of determinants or Schur complements. The corresponding formulation is
\begin{eqnarray}
&&\epsilon_{-1}^{(n)}=0,\quad \epsilon_{0}^{(n)}=u_{n},\quad n=0,1,2,\ldots,\label{fsec2112}\\
&&\epsilon_{k+1}^{(n)}=\epsilon_{k-1}^{(n+1)}+(\epsilon_{k}^{(n+1)}-
\epsilon_{k}^{(n)})^{-1},\quad k,n=0,1,2,\ldots,\label{fsec2113}
\end{eqnarray} 
where $\epsilon_{2k}^{(n)}:=e_{k}(u_{n}), \epsilon_{2k+1}^{(n)}:=(e_{k}(\Delta u_{n}))^{-1}$ and works along diagonals on $n+k$ constant. Thus, from (\ref{fsec2111}), formulas (\ref{fsec2112}), (\ref{fsec2113}) compute each entry of a triangular array in terms of the previous entries.

The extension of the scalar $\epsilon$-algorithm to the vectorial case was carried out by Brezinski, \cite{Brezinski1980}, and Wynn, \cite{Wynn1964,Gekeler1972}, by using different definitions of \lq inverse\rq\ of a vector, see \cite{smithfs}. Wynn suggests to consider the transpose of the Moore-Penrose generalized inverse of a vector, 
\begin{eqnarray}
w^{-1}=\frac{w}{||w||^{2}},\label{fsec2114}
\end{eqnarray} 
leading to the vector $\epsilon$-algorithm (VEA), whose formulas are of the form (\ref{fsec2112}), (\ref{fsec2113}) where the scalars are substituted by vectors and (\ref{fsec2113}) makes use of (\ref{fsec2114}). This implies that the $e$-transform (\ref{fsec2111}) is understood in the above described vectorial sense. This is called the generalized Shanks Schmidt (GSS) transform, \cite{brezinski1}. On the other hand, Brezinski defines the inverse of  pair of vectors $(v,w)$ such that $\langle v,w\rangle\neq 0$ as the pair of vectors $(w^{-1}, v^{-1})$ where
\begin{eqnarray*}
w^{-1}=\frac{v}{\langle w,v\rangle},\quad 
v^{-1}=\frac{w}{\langle w,v\rangle}
\end{eqnarray*}
Thus, $v^{-1}$ is called the inverse of $v$ with respect to $w$ and viceversa. This definition leads to the so-called Topological $\epsilon$-algorithm (TEA), when an arbitrary vector $y$ is fixed and the inverses of $\Delta\epsilon_{2k}^{(n)}$ and $\Delta\epsilon_{2k+1}^{(n)}$ are considered with respect to $y$, that is
\begin{eqnarray*}
\left(\Delta\epsilon_{2k}^{(n)}\right)^{-1}=\frac{y}{\langle y,\Delta\epsilon_{2k}^{(n)}\rangle},\quad
\left(\Delta\epsilon_{2k+1}^{(n)}\right)^{-1}=\frac{y}{\langle y,\Delta\epsilon_{2k+1}^{(n)}\rangle}
\end{eqnarray*}
The recursive formulas are
\begin{eqnarray}
&&\epsilon_{-1}^{(n)}=0,\quad \epsilon_{0}^{(n)}=u_{n},\quad n=0,1,2,\ldots,\nonumber\\
&&\epsilon_{2k+1}^{(n)}=\epsilon_{2k-1}^{(n+1)}+
\left(\Delta\epsilon_{2k}^{(n)}\right)^{-1},\quad k,n=0,1,2,\ldots,\label{fsec2115}\\
&&\epsilon_{2k+2}^{(n)}=\epsilon_{2k}^{(n+1)}+\frac{\Delta\epsilon_{2k}^{(n)}}{\langle \Delta\epsilon_{2k+1}^{(n)},\Delta\epsilon_{2k}^{(n)}\rangle},\quad k,n=0,1,2,\ldots\nonumber
\end{eqnarray} 
Brezinski proved, \cite{brezinski1}, the connection with the GSS-transform, showing that 
\begin{eqnarray*}
\epsilon_{2k}^{(n)}=e_{k}(u_{n}),\quad 
\epsilon_{2k+1}^{(n)}=(e_{k}(u_{n}))^{-1}=\frac{y}{\langle y,e_{k}(u_{n})\rangle}.
\end{eqnarray*}
For an efficient implementation of (\ref{fsec2115}) see \cite{tan}. Thus (TEA) corresponds to take $u_{i,j}=Q(u_{n+i+j})=\langle y,u_{n+i+j}\rangle$ in (\ref{fsec217}), (\ref{fsec218}). 

The mechanism of working of the VEM can be described as follows, see \cite{sidifs,smithfs,Sidi2003} for details. One starts from assuming an asymptotic expression for the sequence $u_{n}$ of the form
\begin{eqnarray}
u_{n}\equiv u+\sum_{j=1}^{\infty} w_{j}\lambda_{j}^{n},\quad n\rightarrow\infty,\label{fsec212}
\end{eqnarray}
where $u\in\mathbb{R}^{m}$ and $\lambda_{j}\in\mathbb{C}$, ordered such that $|\lambda_{j}|\geq |\lambda_{j+1}|, \lambda_{j}\neq 0,1, \lambda_{i}\neq \lambda_{j}$ if $i\neq j$ with only a finite number of $\lambda_{j}$ having the same modulus. The expansion (\ref{fsec212}) can be generalized by considering, instead of constant vectors $w_{j}$, polynomials  $P_{j}(n)$  in $n$ with vector coefficients of the form
\begin{eqnarray}
P_{j}(n)=\sum_{l=0}^{p_{j}}v_{jl}\begin{pmatrix}m\\l\end{pmatrix},
\label{fsec212b}
\end{eqnarray}
with $\{v_{j0},\ldots,v_{jp_{j}}\}$ linearly independent in $\mathbb{R}^{m}$ and where if $|\lambda_{j}|= |\lambda_{j+1}|$ then $p_{j}\geq p_{j+1}$, \cite{SidiB1988}. For simplicity, the description below will make use of (\ref{fsec212}).

The asymptotic expansion (\ref{fsec212}) is considered in a general vector space (finite or infinite dimensional) where the iteration is defined. It is understood in the sense that any truncation differs from $u_{n}$ in less than some power of the next $\lambda$. This means that  for any positive integer $N$ there are $K>0$ and a positive integer $n_{0}$ that only depend on $N$ such that for every $n\geq n_{0}$
\begin{eqnarray*}
||u_{n}-u-\sum_{j=1}^{N-1} w_{j}\lambda_{j}^{n}||\leq K\lambda_{N}^{n}.\label{fsec213}
\end{eqnarray*}
In particular, the case $N=1$ allows to identify $u$ as limit or anti-limit of the sequence $u_{n}$, according to the size of $\lambda_{1}$. (That is, if $|\lambda_{1}|<1$ then $\lim_{n\rightarrow\infty}u_{n}$ exists and equals $u$. If $|\lambda_{1}|>1$ then $\lim_{n\rightarrow\infty}u_{n}$ does not exist and $u$ is called the anti-limit of the sequence $u_{n}$.) Under these conditions, several results of convergence for MPE, RRE, MMPE and TEA are obtained in the literature, \cite{sidifs,smithfs} and references therein. For these methods, one can find an extrapolation step $\kappa$ such that
\begin{eqnarray}
||t_{n,\kappa}-u||=O(\lambda_{\kappa+1}^{n}),\quad n\rightarrow\infty.\label{fsec213b}
\end{eqnarray}
The estimate (\ref{fsec213b}) may explain the convergent behaviour of the extrapolation in some cases. If the $\lambda$'s are identified as the eigenvalues of the linearization operator of the iteration at the limit (or anti-limit) $u$, then the extrapolation has the effect of translating the behaviour of the iteration to an eigenvalue $\lambda_{\kappa+1}$ that may be into the unit disk, even if the previous ones are out of it. Hence, $u$ may be anti-limit for the original iteration and the extrapolation converges to it. These results are extended to the defective linear case with more general polynomials (\ref{fsec212b}) in \cite{SidiB1988}.

The integer $\kappa$ is related to the concept of minimal polynomial $P(\lambda)$ of a matrix $A$ with respect to a vector $v$, \cite{JbilouS1991,jbilous,smithfs}; this is the unique polynomial of least degree such that
$$P(A)v=0.$$ Thus in the case of linear iteration with matrix $A$, $\kappa$ is taken to be the degree of the minimal polynomial of $A$ with respect to the first iteration $u_{0}$. In the case of a nonlinear system written in fixed point form
\begin{eqnarray}
x=\mathcal{F}(x),\label{fsec2110}
\end{eqnarray}
then $\kappa$ is theoretically defined as the degree of the minimal polynomial of $A={\mathcal{F}}^{\prime}(u^{*})$ with respect to $u_{0}$, where $u^{*}$ is a solution of (\ref{fsec2110}). In contrast with the linear case, there is no way to determine $\kappa$ in advance for the nonlinear case. This forces to consider several strategies for the choice and the corresponding implementation, see the discussion in \cite{smithfs} and the comments here below.

We also mention that in the linear case, the extrapolation methods MPE and RRE are mathematically equivalent to the method of Arnoldi, \cite{Saad1981}, and the GMRES, \cite{SaadS1986}, respectively, see \cite{Sidi1988}, while the MMPE is mathematically equivalent to the Hessenberg method, \cite{sidifs} and TEA to the method of Lanczos, \cite{Lanczos1952}, see \cite{Sidi1988,jbilous}. 

Efficient and stable implementation of the RRE, MPE and MMPE methods by using QR and LU factorizations can be seen in \cite{Sidi1991,JbilousS1999}. For the case of TEA, see \cite{tan} and \cite{BrezinskiR1974,sidi} for VEA. The implementation is usually carried out in a cycling mode. A cycle of the iteration is performed by the following steps : consider a method (\ref{mm2}) with a stabilizing factor $s$ satisfying (P1), (P2). Given $u_{0}\neq 0$ and a width of extrapolation $mw\geq 1$, for
$l=0,1,\ldots$, the advance $l\mapsto l+1$ is:
\begin{itemize}
\item[(A)] Set $t_{0}=u_{l}$ and compute $mw$ steps of the fixed-point algorithm:
\begin{eqnarray*}
Lt_{n+1}=s(t_{n})N(t_{n}), n=0,\ldots mw-1.\label{fsec2116}
\end{eqnarray*}
\item[(B)] Compute the extrapolation steps (\ref{fsec217}) with any of the methods described above and $ n=0,\ldots,mw$.
\item[(C)] Set $u_{l+1}=t_{mw,l}, t_{0}=u_{l+1}$ and go to step (A).
\end{itemize}
The cycle (A)-(B)-(C) is repeated until the error (residual or between two consecutive iterations) is below a prefixed tolerance, a maximum number of iterations is attained or the discrepancy between the stabilizing factor at the iterations and one is below a prefixed tolerance.

The width of extrapolation $mw$ depends on the choice of the technique: $mw=\kappa+1$ for MPE, RRE or MMPE and $mw=2\kappa$ for VEA or TEA, \cite{smithfs}. Since $\kappa$ is generally unknown, some strategy for the implementation must be adopted. In practice, as discussed in \cite{smithfs},  the methods are implemented with some (small) values of $mw$ and take that with the best performance. It may be also different for each cycle, although quadratic convergence is not expected if $\kappa$ is too small. This choice of $\kappa$ will be computationally studied in some examples in Section \ref{se3}.

The hypotheses for the expansion (\ref{fsec212}) include the conditions $\lambda_{j}\neq 1 \;\; \forall j$. In many problems for traveling wave generation, $\lambda=1$ appears as eigenvalue of the iterative technique (of fixed point type) although under especial circumstances that allow to extend the convergence results in some sense. This especial situation is related to the presence of symmetries in the equations for traveling waves. In order to extend the results of convergence to this case, one has to consider the orbits by the symmetry group of the equations and interpret the convergence in the orbital sense, that is a convergence for the orbits. The description of this orbital convergence can be seen in \cite{alvarezd}.

 Finally, local quadratic convergence is proved in \cite{smithfs} (see also \cite{leferrand,jbilous,JbilouS1991,Van1994}) for the four methods and VEA for a general nonlinear system (\ref{fsec2110}) and under the following hypotheses on $\mathcal{F}$:
\begin{itemize}
\item The Jacobian matrix ${\mathcal{F}}^{\prime}(u^{*})$ does not have $\lambda=1$ as eigenvalue.
\item $k$ is the degree of the minimal polynomial of ${\mathcal{F}}^{\prime}(u^{*})$ with respect to $u_{0}-u^{*}$.
\item the algorithm is implemented in the cycling mode where $\kappa$ is chosen in the $i$-th cycle as the degree of the minimal polynomial of ${\mathcal{F}}^{\prime}(u^{*})$ with respect to $t_{i-1}-u^{*}$.
\end{itemize}

This result can be extended to the case where ${\mathcal{F}}$ admits a $\nu$-parameter ($\nu\geq 1$) group of symmetries and, consequently, $\lambda=1$ is an eigenvalue of ${\mathcal{F}}^{\prime}(u^{*})$, by using the reduced system for the orbits of the group and in the orbital sense, \cite{alvarezd}. This would lead to quadratic orbital convergence but linear local convergence.

\subsection{Anderson acceleration methods}
A second family of acceleration techniques considered here is the so-called Anderson family or Anderson mixing, \cite{anderson}. It is widely used in electronic structure computations and only recently it has been analyzed in a more general context, \cite{walkern,ni,yangmlw}. (To our knowledge, this is the first time that AAM are applied to accelerate traveling wave computations.) The main goal of the approach consists of combining the iteration with a minimization problem for the residual at each step. For linear problems, this technique is essentially equivalent to the GMRES method, \cite{SaadS1986,Demmel,GolubV}. The stages for an iteration step are as follows:
Given $u_{0}\neq 0, nw\geq 1$,
set $Lu_{1}=s(u_{0})N(u_{0})$.
For
$k=1,2,\ldots$
\begin{itemize}
\item Set $n_{k}=\min \{nw,k\}$
\item Set $F_{k}=(f_{k-n_{k}},\ldots ,f_{k})$ where $f_{i}=Lu_{i}-s(u_{i})N(u_{i})$
\item Determine $\alpha^{(k)}=(\alpha_{0}^{(k)},\ldots,\alpha_{n_{k}}^{(k)})$ that solves
\begin{eqnarray}
\min_{\alpha=(\alpha_{0},\ldots,\alpha_{n_{k}})}||F_{k}\alpha||,\quad \sum_{i=0}^{n_{k}}\alpha_{i}=1\label{minimize}
\end{eqnarray}
\item Set
$$
Lu_{k+1}=\sum_{i=0}^{n_{k}}\alpha_{i}^{(k)}s(u_{k-n_{k}+i})N(u_{k-n_{k}+i})
$$
\end{itemize}
The resolution of the optimization problem (\ref{minimize}) is the source of the additional computational work of the acceleration. One way to reduce this extra effort is the so-called multisecant updating \cite{fangs,Eyert1996}. (This also clarifies the connection with quasi-Newton methods.) This technique consists of writing the problem in an equivalent form
\begin{eqnarray}
\min_{\gamma=(\gamma_{0},\ldots,\gamma_{n_{k}})}||f-{\mathcal F}_{k}\gamma||, {\mathcal F}_{k}=(\Delta f_{k-n_{k}},\ldots ,\Delta f_{k-1}), \Delta f_{i}=f_{i+1}-f_{i},\label{minip}
\end{eqnarray}
but with a more direct resolution, and determining the acceleration from it. The general step becomes:
\begin{itemize}
\item Set $n_{k}=\min \{nw,k\}$
\item Determine $\gamma^{(k)}=(\gamma_{0}^{(k)},\ldots,\gamma_{n_{k}}^{(k)})$ by solving (\ref{minip}).
\item Set $$\alpha_{0}^{(k)}=\gamma_{0}^{(k)}, \alpha_{i}^{(k)}=\gamma_{i}^{(k)}-
\gamma_{i}^{(k)}, 1\leq i\leq n_{k}-1, \alpha_{n_{k}}^{(k)}=1-\gamma_{n_{k}-1}^{(k)}$$
\item Set
$$
Lu_{k+1}=\sum_{i=0}^{n_{k}}\alpha_{i}^{(k)}s(u_{k-n_{k}+i})N(u_{k-n_{k}+i}).
$$
\end{itemize}
As mentioned in \cite{walkern}, if ${\mathcal F}_{k}$ is full-rank, the solution of the minimization problem can be written as $\gamma^{(k)}=({\mathcal F}_{k}^{T}{\mathcal F}_{k})^{-1}{\mathcal F}_{k}^{T}f_{k}$ and the Anderson acceleration has the alternative form
\begin{eqnarray}
Lu_{k+1}&=&Lu_{k}-G_{k}f_{k},\nonumber\\
G_{k}&=&-I+({\mathcal H}_{k}+{\mathcal F}_{k})({\mathcal F}_{k}^{T}{\mathcal F}_{k})^{-1}{\mathcal F}_{k}^{T},\quad
{\mathcal H}_{k}=(\Delta u_{k-m_{k}},\ldots,\Delta u_{k-1}),\nonumber\\
&& \Delta u_{i}=u_{i+1}-u_{i}.\label{aaalternative}
\end{eqnarray}
(Note that $G_{k}$ can be viewed as an approximate inverse of the Jacobian of $Lx-N(x)$). The formulation (\ref{aaalternative}) motivates the generalization of the Anderson mixing, \cite{fangs}. This is performed replacing  $({\mathcal F}_{k}^{T}{\mathcal F}_{k})^{-1}{\mathcal F}_{k}$ in (\ref{aaalternative}) by some ${\mathcal V}_{k}\in\mathbb{R}^{n\times m}$ satisfying
\begin{eqnarray*}
{\mathcal V}_{k}^{T}{\mathcal F}_{k}=I,
\end{eqnarray*}
and (\ref{aaalternative}) becomes
\begin{eqnarray}
Lu_{k+1}&=&Lx_{k}-\widetilde{G}_{k}f_{k},\nonumber\\
\widetilde{G}_{k}&=&-I+({\mathcal H}_{k}+{\mathcal F}_{k}){\mathcal V}_{k}^{T},\quad
{\mathcal H}_{k}=(\Delta u_{k-m_{k}},\ldots,\Delta u_{k-1}).\label{aaalternative2}
\end{eqnarray}
The resulting methods are collected in the so-called Anderson's family, \cite{walkern,fangs}. Two particular members are emphasized: the Type-I method (denoted by AA-I from now on), which corresponds to ${\mathcal V}_{k}=({\mathcal H}_{k}^{T}{\mathcal F}_{k})^{-1}{\mathcal H}_{k}$ in (\ref{aaalternative2}) and Type-II method (or AA-II), which is the original Anderson mixing (\ref{aaalternative}).

To our knowledge, some convergence results can be seen in \cite{walkern,PotraE2013,TothK2015}. In \cite{walkern} the authors identify some Anderson methods for linear problems and in some sense with the GMRES method and the Arnoldi (FOM) method; some convergence results can be derived from this identification. In \cite{PotraE2013} the equivalence with GMRES for linear problems is completely characterized. Finally, \cite{TothK2015} gives some proofs of convergence of the Anderson acceleration when applied to contractive mappings: $q$-linear convergence of the residual for linear problems under certain conditions when $nw=1$ and local $r$-linear convergence in the nonlinear case. (These types of convergence are defined in the paper.)

On the other hand, as
observed in \cite{walkern}, the implementation of AAM should be carried out by attending to three main points: a convenient formulation of the minimization problem (\ref{minimize}), a numerical method for its efficient resolution and, finally, the parameter $nw$, which plays a similar role to that of the extrapolation width $mw$ in the VEM. In our computations below, we have followed the treatment described in \cite{walkern}. This is based on the use of the unconstrained form (\ref{minip}) and its numerical resolution with $QR$ decomposition. For other alternatives in both problems, see the discussion in \cite{walkern,fangs} and the references cited there. (According to our results below, the use of alternative preconditioning techniques might be recommendable in some cases.) As far as the choice of $nw$ is concerned, a similar strategy to that of $mw$ will be used, since our experiments, \cite{walkern}, suggest that $nw$ (as $mw$) strongly depends on the problem under study and large values are not recommended.
Finally the codes are implemented by retaining the definition of $n_{k}=\min\{nw,k\}$ since other alternatives, \cite{yangmlw}, did not improve the results in a relevant way.

\section{Numerical comparisons}
\label{se3}

Presented here is a comparative study on the use of VEM and AAM as acceleration techniques from the \PM type methods (\ref{mm2}) in traveling wave computations. The comparison is organized according to two main points: the type of traveling wave to be generated and the elements of each family of methods to be used for the generation. The first point includes the following case studies:
\begin{enumerate}
\item Classical solitary waves, generalized solitary waves and periodic traveling waves of the four-parameter Boussinesq system, \cite{Bona_Chen_Saut_1,Bona_Chen_Saut_2}.
\item Localized ground state solutions of NLS type equations, \cite{Yang2012,yang2}.
\item Highly oscillatory solitary waves of the one- and two-dimensional Benjamin equation, \cite{ben0,ben1,ben2,kim,kima1,kima2}.
\end{enumerate}
This plethora of waves attempts to discuss and overcome different computational difficulties and with the aim of establishing as more general conclusions as possible. As for the second point, each family of techniques has been represented by the following methods:
\begin{itemize}
\item MPE and RRE standing for polynomial extrapolation methods.
\item VEA and TEA standing for $\epsilon$ algorithms.
\item AA-I and AA-II standing for the AAM.
\end{itemize}
For simplicity, the acceleration will be applied to the \PM method (\ref{mm2}), (\ref{mm3c}) with $\gamma=p/(p-1)$. Due to the similar behaviour of the methods of the family (\ref{mm2}), illustrated in \cite{alvarezd}, the conclusions from the corresponding results can reasonable serve when the \PM method is substituted by any of (\ref{mm2}).


In all the cases considered, the traveling wave profiles appeared as solutions of initial value problems of ode's. Their discretization to generate approximations to the profiles, was carried out in a common and standard way. The corresponding initial periodic boundary value problem (on a sufficiently long interval) was discretized by using Fourier collocation techniques, \cite{Boyd,Canutohqz}; the discretization leads to a  nonlinear system of algebraic equations for the approximate values of the profile at the collocation points or for the discrete Fourier coefficients of the approximation. This system is iteratively solved with the classical Petviashvili method (\ref{mm2}), (\ref{mm3c}) along with the selected acceleration technique. This will be described in each equation considered below.

Several stopping criteria for the iterations are implemented:
\begin{itemize}
\item A maximum number of iterations.
\item The iteration stops when one of the following quantities are below a prefixed, small tolerance $TOL$:
\begin{itemize}
\item[(i)] The difference in Euclidean norm between two consecutive iterations
\begin{eqnarray*}
\label{fsec31}
E_{n}=||u_{n+1}-u_{n}||,\quad n=0,1,\ldots
\end{eqnarray*}
\item[(ii)] The residual error (also in Euclidean norm)
\begin{eqnarray}
\label{fsec32}
RES_{n}=||Lu_{n}-N(u_{n})||,\quad n=0,1,\ldots
\end{eqnarray}
\item[(iii)]
The discrepancy between the stabilizing factor and (in case of convergence) its limit one
\begin{eqnarray}
\label{fsec33}
SFE_{n}=|s(u_{n})-1|,\quad n=0,1,\ldots
\end{eqnarray}
\end{itemize} 
\end{itemize}

The numerical experiments that form the comparative study are of different type:
\begin{itemize}
\item For several values of $\kappa$ (in the case of VEM) and $nw$ (in the case of AAM) we have computed the number of iterations required by each method to achieve a residual error below the prefixed tolerance. This allows to compare some performance of the methods between the two families, among different techniques within a same family and indeed with the \PM method without acceleration.
\item Some eigenvalues of the iteration matrices for the classical fixed point algorithm and the \PM method have been computed (with the corresponding standard MATLAB function) in order to explain the behaviour of the second one, \cite{alvarezd} and how the acceleration eventually changes it.
\item The form of the approximate profiles and some experiments to check their accuracy are also displayed.
\end{itemize}

\subsection{Traveling wave solutions of Boussinesq systems}
\label{se31}
In this first example we study the numerical generation of traveling wave solutions of the four-parameter family of Boussinesq system
  \begin{eqnarray}
\eta_{t}+u_{x}+(\eta u)_{x}+au_{xxx}-b\eta_{xxt}&=&0,\label{fsec311}\\
u_{t}+\eta_{x}+uu_{x}+c\eta_{xxx}-du_{xxt}&=&0,\nonumber
\end{eqnarray}
where $\eta=\eta(x,t), u=u(x,t), x\in\mathbb{R}, t\geq 0$ and the four parameters $a, b, c, d$ satisfy
\begin{eqnarray}
a+b=\frac{1}{2}(\theta^{2}-\frac{1}{3}),\quad c+d=\frac{1}{2}(1-\theta^{2}),\label{fsec312}
\end{eqnarray}
with some constant $\theta^{2}\in [0,1]$, \cite{Boussinesq,Bona_Chen_Saut_1,Bona_Chen_Saut_2}. System (\ref{fsec311}) appears as one of the alternatives to model the bidirectional propagation of the irrotational free surface flow of an incompressible, inviscid fluid in a uniform horizontal channel under the effects of gravity when the surface tension and cross-channel variations of the fluid are assumed to be negligible. If $h_{0}$ denotes the undisturbed water depth, then $\eta(x,t)$ stands for the deviation of the free surface from $h_{0}$ at the point $x$ and time $t$, while $u(x,t)$ is the horizontal velocity of the fluid at $x$ and at the height $y=\theta h_{0}$ (where $y=0$ corresponds to the channel bottom) at time $t$. For the derivation of (\ref{fsec311}) from the two-dimensional Euler equations and the mathematical theory see \cite{Bona_Chen_Saut_1,Bona_Chen_Saut_2}. For the modification of (\ref{fsec312}) to include the influence of surface tension see \cite{DaripaD2003,ChenNS2009,ChenNS2011}.

The Boussinesq system (\ref{fsec311}) admits different types of traveling wave solutions. First, being an approximation to the corresponding two-dimensional Euler equations in the theory of surface waves, it is expected to have classical  solitary wave solutions. They are solutions of the initial value problem of (\ref{fsec311}) of smooth traveling wave form $\eta=\eta(x-c_{s}t), u=u(x-c_{s}t)$ with some speed $c_{s}>0$ and decaying to zero as $X=x-c_{s}t\rightarrow\pm\infty$. Substitution into (\ref{fsec311}) and after one integration the profiles $u=u(X), \eta=\eta(X)$ must satisfy the ode system
\begin{eqnarray}
\label{fsec313} \begin{pmatrix}
    c_{s}(1-b\partial_{XX})&-(1+a\partial_{XX}) \\
    -(1+c\partial_{XX})&c_{s}(1-d\partial_{XX})\\
    \end{pmatrix}\begin{pmatrix}
    \eta \\
    u\\
    \end{pmatrix}=\begin{pmatrix}
    u\eta \\
    \frac{u^{2}}{2}\\
    \end{pmatrix}.
\end{eqnarray}
The problems of existence, asymptotic decay and stability of solutions of (\ref{fsec313}) have been analyzed in many references and for particular values of $a, b, c, d$, see \cite{DougalisM2008} and references therein. Furthermore, in the same reference, linearly well-posed systems (\ref{fsec313}) may be studied as a first order ode system
 and, based on normal form theory, a discussion on the values of the parameters leading to Boussinesq systems admitting solitary wave solutions with speed $c_{s}>1$ is established. According to it, two classes of systems can be distinguished. The first one admits classical (in the sense above defined) solitary wave solutions. This group contains the Bona-Smith system ($a=0, b=d=(3\theta^{2}-1)/6, c=(2-3\theta^{2})/3, 2/3<\theta^{2}<1$), \cite{BonaS1976}, or the BBM-BBM system ($a=c=0, b=d=1/6$), \cite{BonaC1988}. (The classical Boussinesq system, which corresponds to $a=b=c=0, d=1/3$, is also in this group, although it is out of the general discussion of \cite{DougalisM2008} and has been studied separately, \cite{pegow}.) A second class of Boussinesq systems admits generalized solitary wave solutions, that is traveling wave profiles which are not homoclinic to zero at infinity but to small amplitude periodic waves, \cite{Lombardi2000}. The KdV-KdV system ($a=c>0, b=d=0$), \cite{BonaDM2007,BonaDM2008}, is an example of this second group. 

Finally, the existence of periodic traveling wave solutions of (\ref{fsec311}) is studied in \cite{Chen_Chen_Nguyen} by applying topological degree theory for positive operator to the corresponding periodic initial value problem posed on an interval $(-l,l)$ and some cnoidal wave solutions of the BBM-BBM system are computed. A smooth periodic traveling wave solution $\eta=\eta(x.c_{s}t), u=u(x.c_{s}t)$ with some speed $c_{s}>0$ must satisfy the ode system
\begin{eqnarray}
\label{fsec313b} \begin{pmatrix}
    c_{s}(1-b\partial_{xx})&-(1+a\partial_{xx}) \\
    -(1+c\partial_{xx})&c_{s}(1-d\partial_{xx})\\
    \end{pmatrix}\begin{pmatrix}
    \eta \\
    u\\
    \end{pmatrix}=\begin{pmatrix}
    u\eta \\
    \frac{u^{2}}{2}\\
    \end{pmatrix}+\begin{pmatrix}
    K_{1} \\
    K_{2}\\
    \end{pmatrix},
\end{eqnarray}
for some real constants $K_{1}, K_{2}$. These are related to the period parameter $l$. The resolution of (\ref{fsec313b}) involves a modified system for  which these constants of integration are set to zero, \cite{Chen_Chen_Nguyen}. This is briefly described as follows. One first searches for constant solutions $\eta=C_{1}, u=C_{2}$ of (\ref{fsec313b}). This leads to the system
\begin{eqnarray*}
\begin{pmatrix}
    c_{s}&-1\\
    -1&c_{s}\\
    \end{pmatrix}\begin{pmatrix}
    C_{1}\\
    C_{2}\\
    \end{pmatrix}=\begin{pmatrix}
    C_{1}C_{2} \\
    \frac{C_{2}^{2}}{2}\\
    \end{pmatrix}+\begin{pmatrix}
    K_{1} \\
    K_{2}\\
    \end{pmatrix},
\end{eqnarray*}
which can be solved as a cubic equation for $C_{2}$:
\begin{eqnarray}
&&\frac{C_{2}^{3}}{2}-\frac{3}{2}c_{s}C_{2}^{2}+(c_{s}^{2}-1+K_{1})C_{2}-K_{2}-c_{s}K_{1}=0,\label{fsec313c}\\
&&C_{1}=c_{s}C_{2}-C_{2}^{2}/2-K_{2}.\label{fsec313c1}
\end{eqnarray}
Once $(C_{1},C_{2})$ is obtained (there may be more than one solution indeed) the differences $\widetilde{\eta}=\eta-C_{1}, \widetilde{u}=u-C_{2}$ must satisfy
\begin{eqnarray}
\label{fsec313d} \begin{pmatrix}
    c_{s}(1-b\partial_{xx})-C_{2}&-(1+a\partial_{xx}) -C_{1}\\
    -(1+c\partial_{xx})&c_{s}(1-d\partial_{xx})-C_{2}\\
    \end{pmatrix}\begin{pmatrix}
    \widetilde{\eta} \\
    \widetilde{u}\\
    \end{pmatrix}=\begin{pmatrix}
    \widetilde{u}\widetilde{\eta}\\
    \frac{\widetilde{u}^{2}}{2}\\
    \end{pmatrix}.
\end{eqnarray}
This strategy will be considered in the numerical generation of the profiles in (\ref{fsec313b}): the system (\ref{fsec313d}) wil be discretized to compute approximations to the variables $\widetilde{\eta}, \widetilde{u}$ and to the variables $\eta=\widetilde{\eta}+C_{1}, u=\widetilde{u}+C_{2}$ from them.

In order to generate numerically classical and generalized solitary wave solutions of (\ref{fsec311}) the corresponding periodic value problem of (\ref{fsec313}) on a long enough interval $(-l,l)$ is discretized  with a Fourier collocation method leading to a discrete system of the form
\begin{eqnarray}
\label{fsec314} \underbrace{ \begin{pmatrix}
    c_{s}(I_{m}-bD^{2})&-(I_{m}+aD^{2}) \\
    -(I_{m}+cD^{2})&c_{s}(I_{m}-dD^{2})\\
    \end{pmatrix}}_{L}\begin{pmatrix}
    \eta_{h} \\
    u_{h}\\
    \end{pmatrix}
    =\underbrace{\begin{pmatrix}
    u_{h}.\eta_{h} \\
    \frac{u_{h}.^{2}}{2}\\
    \end{pmatrix}
    }_{N(\eta_{h},u_{h})}.
\end{eqnarray}
where
$\eta_{h}, u_{h}\in \mathbb{R}^{m}$ are approximations
$\eta_{h,j}\approx \eta(x_{j}), u_{h,j}\approx u(x_{j})$ to the
values of a solution of (\ref{fsec313b}) at the grid points
$x_{j}=-l+jh, h=2l/m, j=0,\ldots m-1$, $D$ is the pseudospectral differentiation matrix,
\cite{Boyd,Canutohqz}, $I_{m}$ is the
$m\times m$ identity matrix and the nonlinear term $N$, which is homogeneous of degree $p=2$, involves Hadamard products. In the case of periodic traveling waves and as was mentioned above, system (\ref{fsec314}) will be substituted in the implementation by
\begin{eqnarray}
\label{fsec314b}\begin{pmatrix}
    c_{s}(I_{m}-bD^{2})-C_{2}I_{m}&-(I_{m}+aD^{2})-C_{1}I_{m} \\
    -(I_{m}+cD^{2})&c_{s}(I_{m}-dD^{2})-C_{2}I_{m}\\
    \end{pmatrix}\begin{pmatrix}
    \widetilde{\eta}_{h} \\
    \widetilde{u}_{h}\\
    \end{pmatrix}
    =\begin{pmatrix}
    \widetilde{u}_{h}.\widetilde{\eta}_{h} \\
    \frac{\widetilde{u}_{h}.^{2}}{2}\\
    \end{pmatrix}  
\end{eqnarray}
for the approximations $\widetilde{\eta}_{h},
    \widetilde{u}_{h}$ to the $\widetilde{\eta},
    \widetilde{u}$ variables at the grid points and where $C_{1}, C_{2}$ are previously known from the resolution of (\ref{fsec313c}), (\ref{fsec313c1}). Then $\eta_{h}=\widetilde{\eta}_{h}+C_{1}, u_{h}=\widetilde{u}_{h}+C_{2}$.

The methods (\ref{mm2})  along with the corresponding acceleration technique are then applied to the discrete systems (\ref{fsec314}) and (\ref{fsec314b}). The implementation is performed in the Fourier space; for example (\ref{fsec314}) becomes
\begin{eqnarray*}
&&  \begin{pmatrix}
    c_{s}(1+b\left(\frac{p\pi}{l}\right)^{2})&-(1-a\left(\frac{p\pi}{l}\right)^{2}) \\
    -(1-c\left(\frac{p\pi}{l}\right)^{2})&c_{s}(1+d\left(\frac{p\pi}{l}\right)^{2})\\
    \end{pmatrix}\begin{pmatrix}
    \left(\widehat{\eta_{h}}\right)_{p} \\
    \left(\widehat{u_{h}}\right)_{p}\\
    \end{pmatrix}
    =\begin{pmatrix}
    \left(\widehat{u_{h}.\eta_{h}}\right)_{p}\\
    \frac{1}{2}\left(\widehat{u_{h}.u_{h}}\right)_{p}\\
    \end{pmatrix},\\
&&-\frac{m}{2}\leq p\leq \frac{m}{2}.
\end{eqnarray*}
Thus the $2m\times 2m$ system (\ref{fsec314}) is divided into $m$ blocks of $2\times 2$ systems for the corresponding $p$-th discrete  Fourier coefficients $\left(\widehat{\eta_{h}}\right)_{p},
    \left(\widehat{u_{h}}\right)_{p}, -m/2\leq p\leq m/2.$ (For simplicity, we assume that $m=2^{s}$ for some $s>1$.) Alternatively, (\ref{fsec314}) can be written in the form
\begin{eqnarray*}
\label{fsec315}
\begin{pmatrix}
    \eta_{h} \\
    u_{h}\\
    \end{pmatrix}=T_{h}\begin{pmatrix}
    \eta_{h} \\
    u_{h}\
    \end{pmatrix}=\begin{pmatrix}
    A_{h}\ast (u_{h}.\eta_{h})+B_{h}\ast
    \frac{u_{h}.^{2}}{2}\\
    B_{h}\ast (u_{h}.\eta_{h})+C_{h}\ast
    \frac{u_{h}.^{2}}{2}\\
    \end{pmatrix},
\end{eqnarray*}
where $\ast$ denotes periodic convolution and if
$\omega=\exp(-2\pi i/m)$, the vectors $A_{h}, B_{h}, C_{h}$ have
discrete Fourier coefficients
\begin{eqnarray*}
&&(\widehat{A_{h}})_{p}=\frac{1-a\left(\frac{p\pi}{l}\right)^{2}}{\Delta(p)},\quad
(\widehat{B_{h}})_{p}=\frac{c_{s}(1+b\left(\frac{p\pi}{l}\right)^{2})}{\Delta(p)},\\
&&(\widehat{C_{h}})_{p}=\frac{1-c\left(\frac{p\pi}{l}\right)^{2}}{\Delta(p)},\quad
(\widehat{D_{h}})_{p}=\frac{c_{s}(1+d\left(\frac{p\pi}{l}\right)^{2})}{\Delta(p)},\\
&&\Delta(p)=c_{s}^{2}(1+b\left(\frac{p\pi}{l}\right)^{2})(1+d\left(\frac{p\pi}{l}\right)^{2})-(1-a\left(\frac{p\pi}{l}\right)^{2})(1-c\left(\frac{p\pi}{l}\right)^{2}),\\
&&-\frac{m}{2}\leq p\leq \frac{m}{2}.
\end{eqnarray*}
In order to explain the behaviour of the iteration, the size of the eigenvalues of the iteration matrix will be relevant in the numerical study. In this case,
the corresponding iteration matrix of the classical fixed point iteration at a solution $u^{*}=(\eta_{h}^{*},u_{h}^{*})$ has the form

\begin{eqnarray*}
\label{fsec316}
S=L^{-1}\begin{pmatrix}
    {\rm diag}(u_{h}^{*})&{\rm diag}(\eta_{h}^{*}) \\
    0&{\rm diag}(u_{h}^{*})\\
    \end{pmatrix},
\end{eqnarray*}
(where ${\rm diag}(v)$ stands for the diagonal matrix with diagonal entries given by the components of $v\in\mathbb{R}^{m}$). Some information on the spectrum of $S$ is known.
We
already have the eigenvalue $\lambda=2$, corresponding to the
degree of homogeneity of the nonlinear part, with
$u^{*}=(\eta_{h}^{*},u_{h}^{*})$ as an eigenvector. Also, the application of $D$ to (\ref{fsec314}) leads to
\begin{eqnarray*}
&&c_{s}D\eta_{h}^{*}-Du_{h}^{*}=D(\eta_{h}^{*}.u_{h}^{*})=u_{h}^{*}.D\eta_{h}^{*}+\eta_{h}^{*}.Du_{h}^{*}\\
&&-D\eta_{h}^{*}+c_{s}\left(I-\frac{1}{3}D^{2}\right)Du_{h}^{*}=D\left(\frac{u_{h}.^{2}}{2}\right)=u_{h}^{*}Du_{h}^{*},
\end{eqnarray*}
which means that $\lambda=1$ is an eigenvalue of $S$ and
$(D\eta_{h}, Du_{h})^{T}$ is an associated eigenvector. This
corresponds to the \lq translational\rq\ invariance of (\ref{fsec311}).

Three particular systems of (\ref{fsec311}) will be taken to illustrate the numerical generation of traveling waves. The first one is the classical Boussinesq system ($a=b=c=0, d=1/3$),
\cite{Boussinesq,Bona_Chen_Saut_1,Bona_Chen_Saut_2}
\begin{eqnarray}
\eta_{t}+u_{x}+(\eta u)_{x}&=&0,\nonumber\\
u_{t}+\eta_{x}+uu_{x}-\frac{1}{3}u_{xxt}&=&0.\label{fsec317}
\end{eqnarray}
which is known to have classical solitary wave solutions, \cite{pegow}. The second one is the so-called KdV-KdV system ($a=c=1/6, b=d=0$)
\begin{eqnarray}
\eta_{t}+u_{x}+(\eta u)_{x}+\frac{1}{6}u_{xxx}&=&0,\nonumber\\
u_{t}+\eta_{x}+uu_{x}+\frac{1}{6}\eta_{xxx}&=&0,\label{fsec318}
\end{eqnarray}
that admits generalized solitary wave solutions, \cite{BonaDM2007,BonaDM2008}. Finally, in order to illustrate the numerical generation of periodic traveling waves, the BBM-BBM system ($a=c=0, b=d=1/6$),
\begin{eqnarray}
\eta_{t}+u_{x}+(\eta u)_{x}-\frac{1}{6}\eta_{xxt}&=&0,\nonumber\\
u_{t}+\eta_{x}+uu_{x}-\frac{1}{6}u_{xxt}&=&0,\label{fsec319}
\end{eqnarray}
will be taken, \cite{Chen_Chen_Nguyen}.

\subsubsection{Numerical generation of classical solitary waves of (\ref{fsec317})}
\begin{figure}[htbp]
\centering \subfigure[]{
\includegraphics[width=6.6cm]{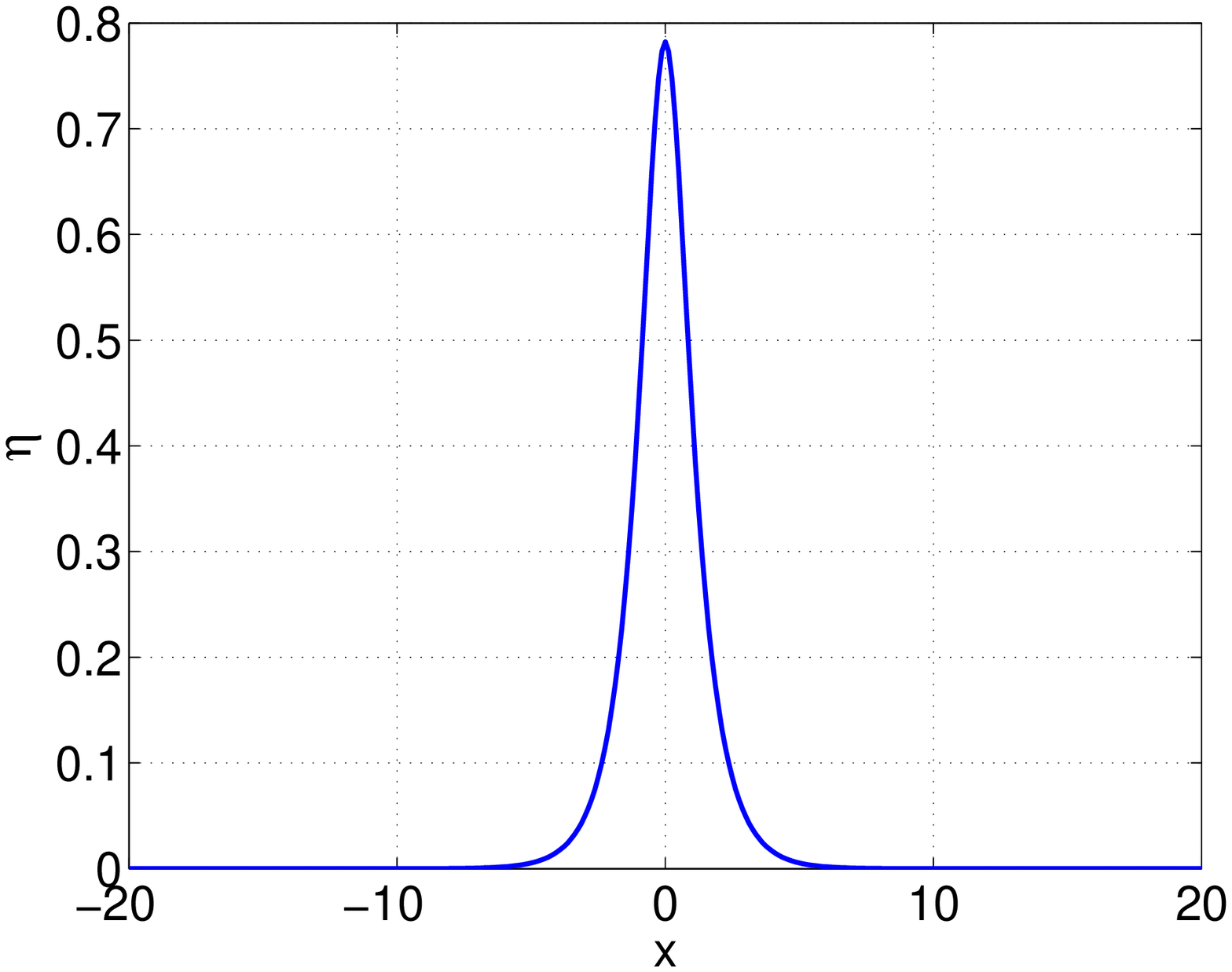} }
\subfigure[]{
\includegraphics[width=6.6cm]{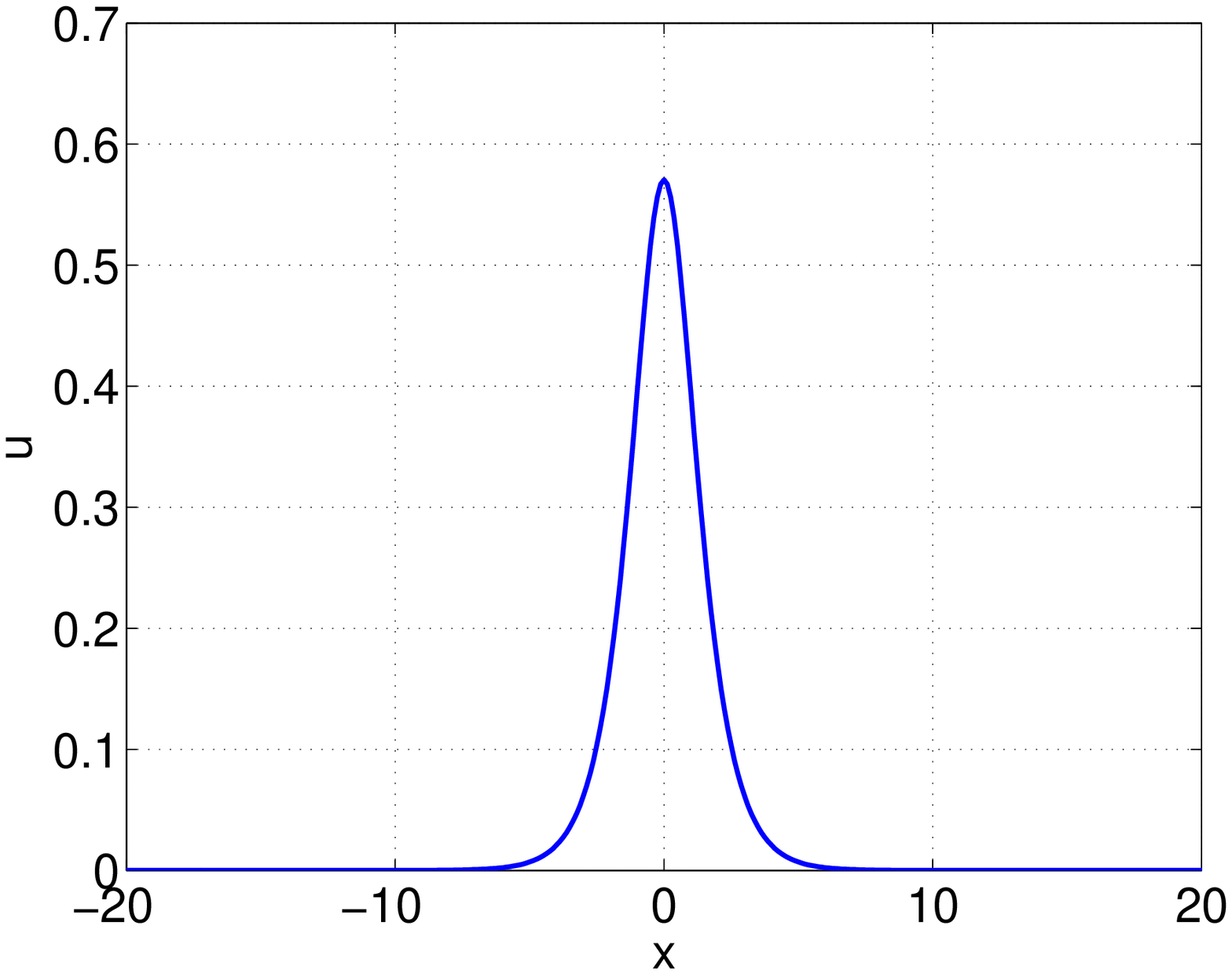} }
\subfigure[]{
\includegraphics[width=6.6cm]{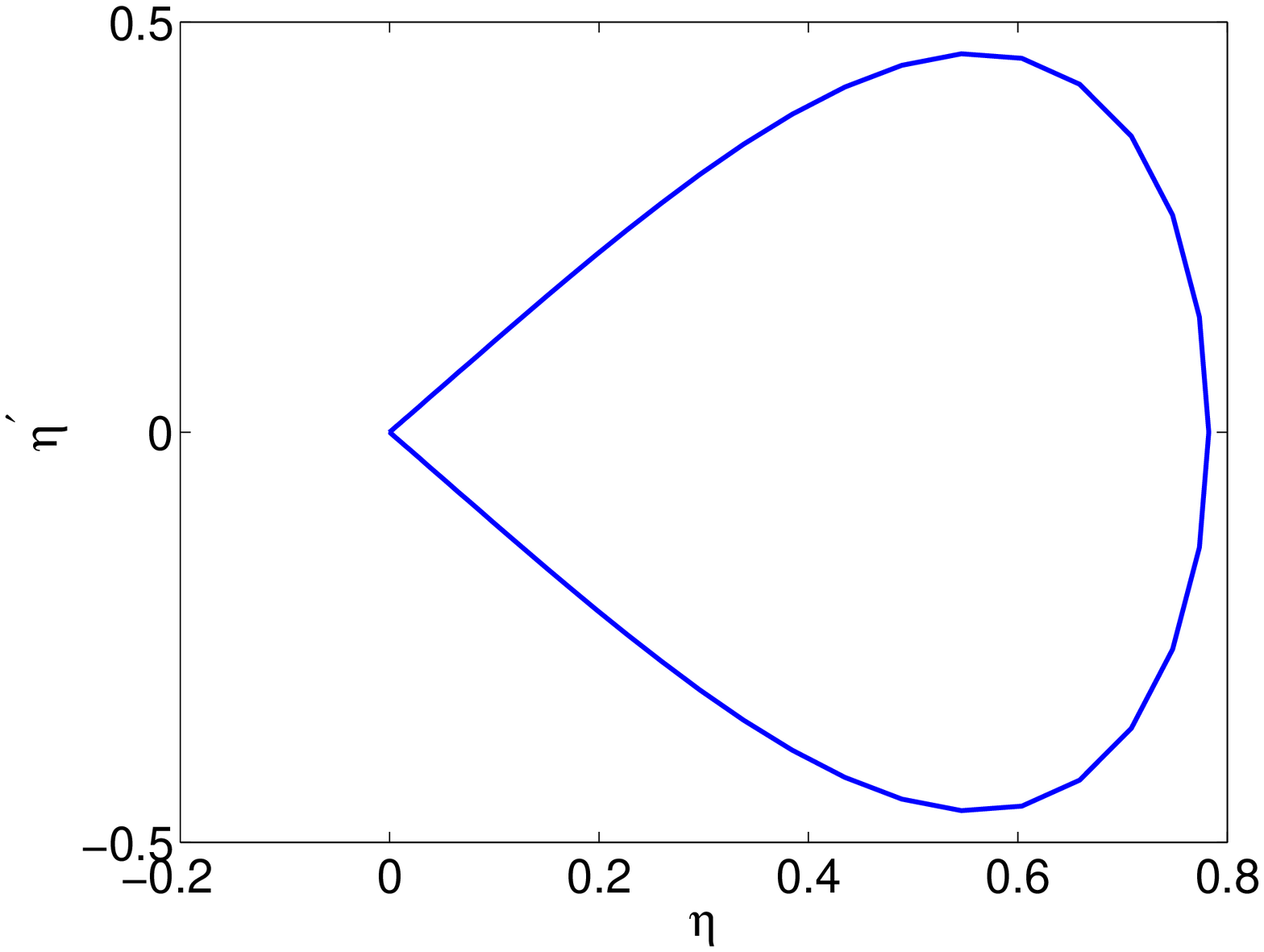} }
\subfigure[]{
\includegraphics[width=6.6cm]{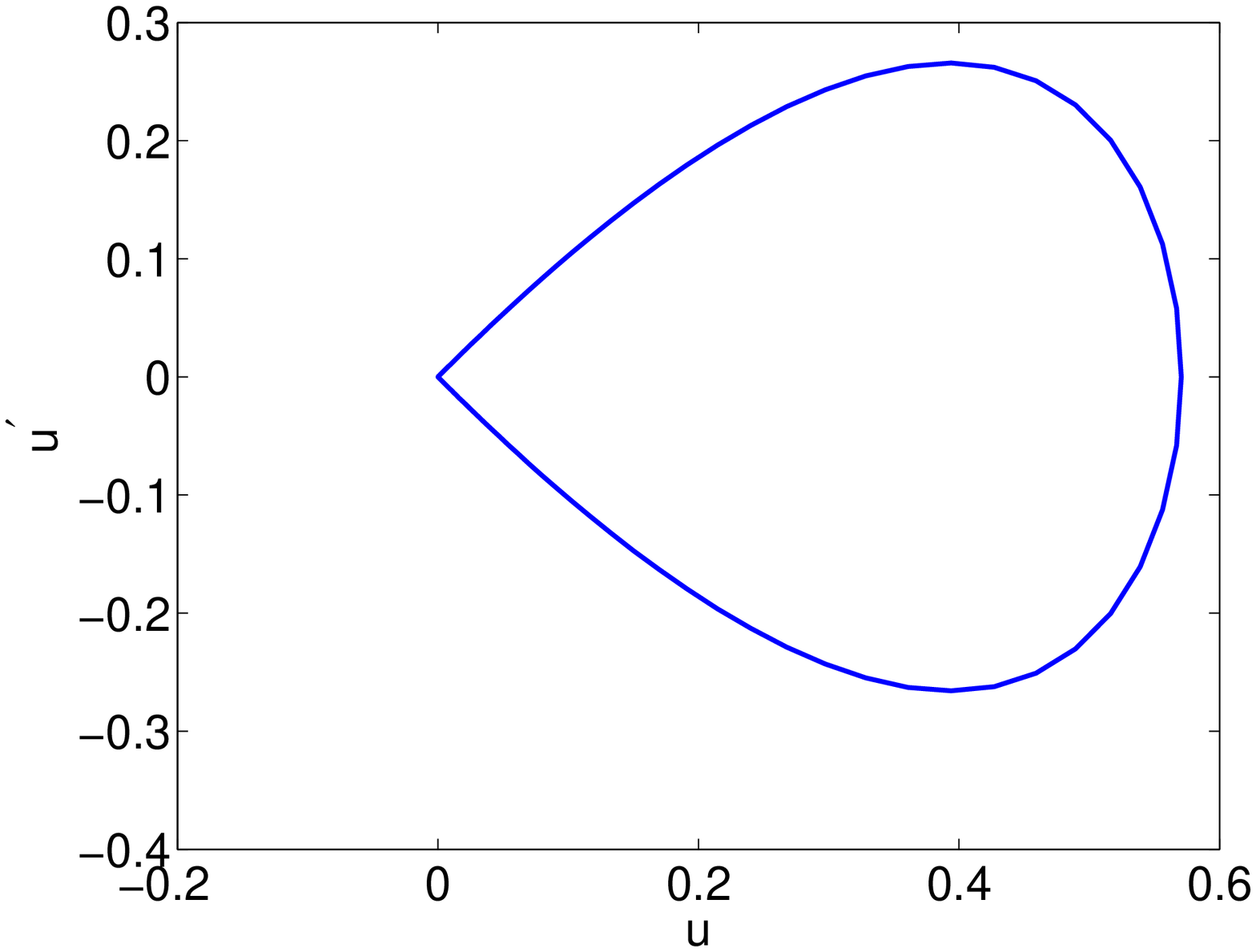} }
\caption{Approximate profiles generated by the \PM method (\ref{mm2}), (\ref{mm3c}) for
(\ref{fsec317}) with $c_{s}=1.3$: (a) $\eta$, (b) $u$; (c) Phase portrait of $\eta$, (d) Phase portrait of $u$.} \label{figuresw0}
\end{figure}

\begin{figure}[htbp]
\centering \subfigure[]{
\includegraphics[width=6.6cm]{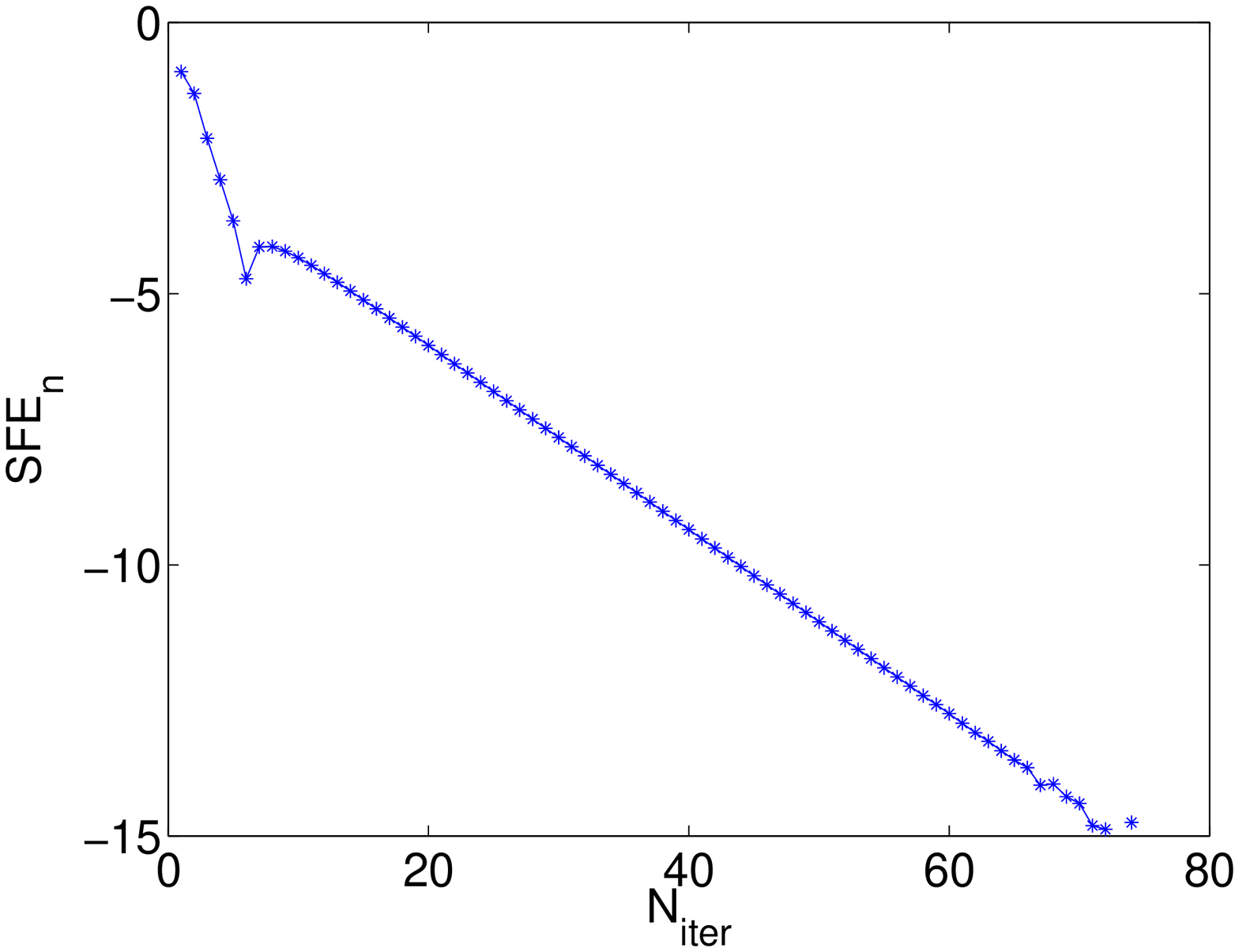} }
\subfigure[]{
\includegraphics[width=6.6cm]{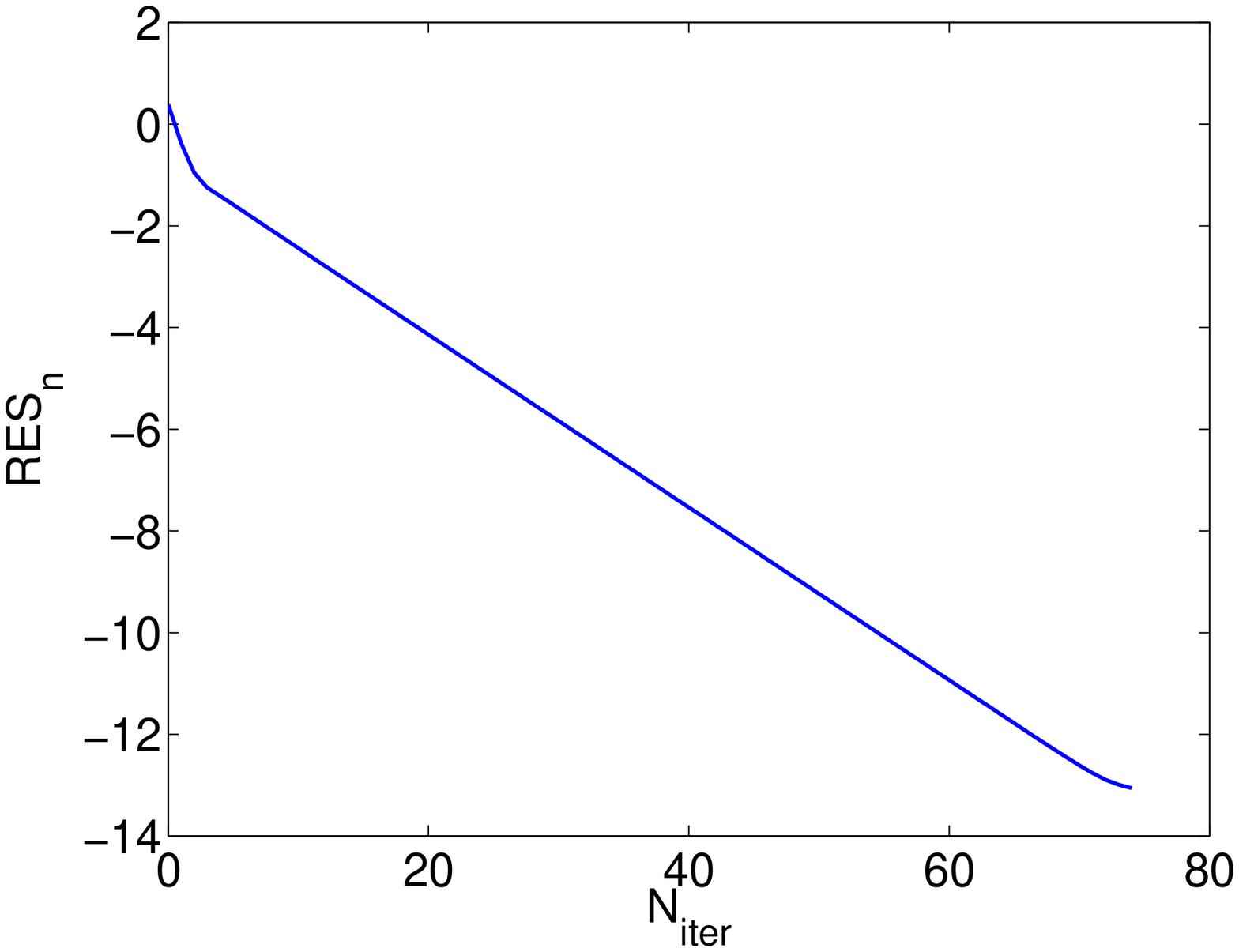} }
\caption{Convergence results of the \PM method for
(\ref{fsec317}): (a) Discrepancy (\ref{fsec33}) for the stabilizing
factor $s_{n}=s(\eta_{n},u_{n})$ vs number of iterations: (b)
Residual error (\ref{fsec32}) vs number of iterations. (Semi-logarithm scale in both cases.)}
\label{figuresw1}
\end{figure}

In the case of system (\ref{fsec317}) a first experiment of comparison of the acceleration techniques has been made by taking $c_{s}=1.3$ and a hyperbolic secant profile as initial iteration with $l=64$ and $m=1024$ collocation points. The Petviashvili method (\ref{mm2}), (\ref{mm3c}) with $\gamma=2$ was first run, generating approximate $\eta$ and $u$ profiles as shown in Figures \ref{figuresw0}(a) and (b) while Figures \ref{figuresw0}(c) and (d) stand for the corresponding phase portraits of the approximate profiles in (a) and (b). (They show the classical character of the solitary waves, represented as homoclinic to zero orbits with exponential decay, \cite{pegow}. )
The accuracy of the iteration is checked in Figure \ref{figuresw1}. Figure \ref{figuresw1}(a) illustrates the convergence of the sequence $s_{n}=s(\eta_{n},u_{n})$ of stabilizing factors, computed with the corresponding formula (\ref{mm3c}) and the optimal value $\gamma=2$. The discrepancy (\ref{fsec33}) is below the tolerance $TOL=10^{-13}$ in $n=62$ iterations, while the first residual error below $TOL$ is $9.092489E-14$ at $n=76$. (This also happens in the rest of the experiments: when the procedure is convergent, the error $|1-s_{n}|$ achieves the tolerance before the residual error; therefore, the control on this last one is a harder test and will be adopted as the main one to stop the iteration.) 
Convergence is also confirmed by Table \ref{tav1b}. This shows the six largest magnitude eigenvalues of the iteration matrix (\ref{fsec316}) (first column) and of the iteration matrix (at the same iterate) of the \PM procedure (second column) both at the last computed iterate  $(\eta_{f},u_{f})$. The first column reveals the dominant eigenvalues $\lambda_{1}=2, \lambda_{2}=1$, both simple, while the rest is below one. The filtering effect of the \PM method, \cite{alvarezd}, is observed in the second column; the dominant eigenvalue is filtered to zero (recall that $\gamma=2$) and the rest is preserved. Since $\lambda_{2}=1$ corresponds to the translational symmetry of (\ref{fsec311}), this guarantees the local convergence of the method (also in the orbital sense mentioned above).

\begin{table}
\begin{center}
\begin{tabular}{|c|c|}
\hline\hline  Iteration matrix $S(\eta_{f},u_{f})$&
Iteration matrix $F^{\prime}(\eta_{f},u_{f})$\\\hline
$1.999999E+00$&$9.999999E-01$\\
$9.999999E-01$&$6.763242E-01$\\
$6.763242E-01$&$5.411229E-01$\\
$5.411229E-01$&$4.820667E-01$\\
$4.820667E-01$&$4.567337E-01$\\
$4.567337E-01$&$4.465122E-01$\\
\hline\hline
\end{tabular}
\end{center}
\caption{Classical solitary wave generation of (\ref{fsec317}) . Six largest magnitude eigenvalues of the
approximated iteration matrix $S=L^{-1}N^{\prime}(\eta_{f},u_{f})$ (first column) and of the iteration matrix $F^{\prime}(\eta_{f},u_{f})$,  generated by the \PM
method (\ref{mm2}), (\ref{mm3c}) with $\gamma=2$, both evaluated
at the last computed iterate $(\eta_{f},u_{f})$.}\label{tav1b}
\end{table}

\begin{figure}[htbp]
\centering \subfigure[]{
\includegraphics[width=8.6cm]{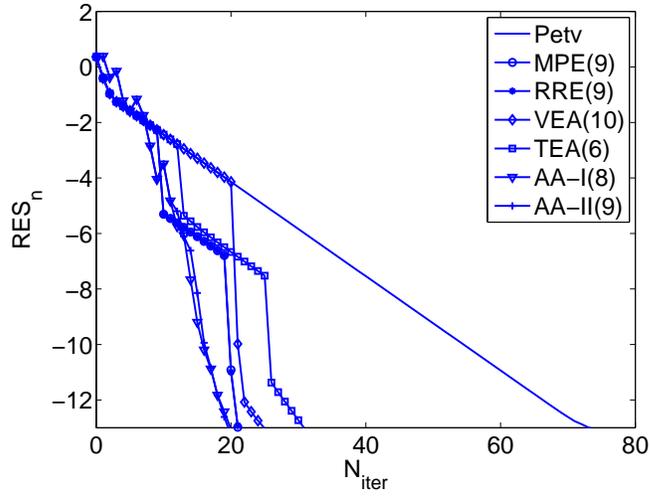} }
\subfigure[]{
\includegraphics[width=8.6cm]{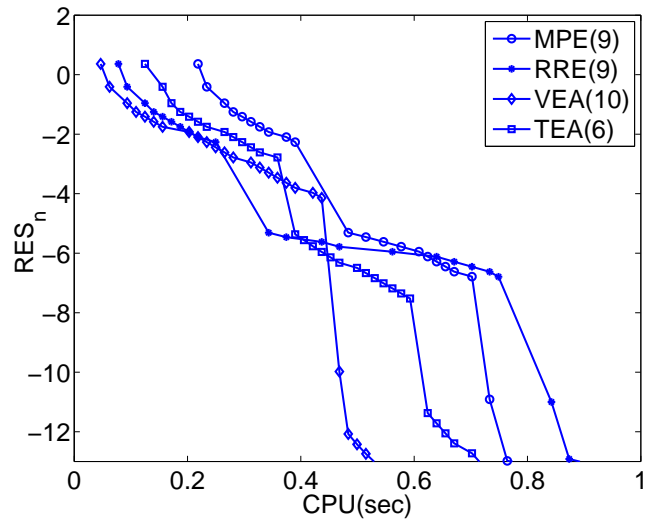} }
\subfigure[]{
\includegraphics[width=8.6cm]{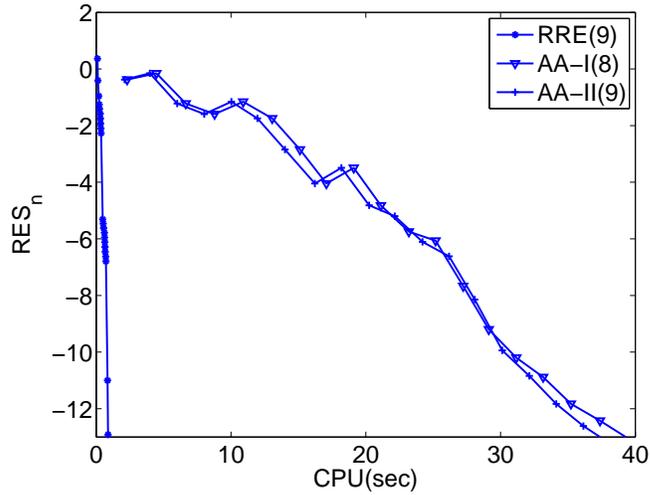} }
\caption{Convergence results of the \PM method for
(\ref{fsec317}):  (a) Residual error (\ref{fsec32}) as function of number of iterations for the \PM method without acceleration and for six acceleration techniques with the best parameters $\kappa$ and $nw$ (according to Tables \ref{tav5b} and \ref{tav6b}). (b) Residual error (\ref{fsec32}) as function of CPU time (in seconds) for six acceleration techniques with the best parameters $\kappa$ and $nw$ (according to Tables \ref{tav5b} and \ref{tav6b}). (c) Comparison  of residual error (against CPU time) between the most efficient VEM and the AAM.
}
\label{figuresw2}
\end{figure}

The improvement of the performance of the \PM method with several acceleration techniques is now computationally analyzed. A first point to study is the choice of the parameters $\kappa$ (for the VEM) and $nw$ (for the AAM). Table \ref{tav2b} shows, for values of $\kappa$ between one and ten, the number of iterations required by MPE, RRE, VEA and TEA to achieve a residual error below $TOL=10^{-13}$. (The residual error, corresponding to the last iteration is in parenthesis for each computation.) From these results, the following comments can be made:
\begin{table}
\begin{center}
\begin{tabular}{|c|c|c|c|c|}
\hline\hline  $\kappa$&MPE($\kappa$)&RRE($\kappa$)&VEA($\kappa$)&TEA($\kappa$)\\\hline
$1$&$269$&$99$&$631$&$408$\\
&($8.9136E-14$)&($9.1312E-14$)&($9.7920E-14$)&($7.3356E-14$)\\
$2$&$64$&$48$&$43$&$43$\\
&($7.4794E-14$)&($8.3903E-14$)&($7.7841E-14$)&($8.1979E-14$)\\
$3$&$43$&$43$&$38$&$42$\\
&($7.5001E-14$)&($8.0682E-14$)&($7.7395E-14$)&($9.6255E-14$)\\
$4$&$33$&$33$&$33$&$37$\\
&($7.9601E-14$)&($8.2823E-14$)&($7.7824E-14$)&($8.0527E-14$)\\
$5$&$28$&$28$&$31$&$39$\\
&($8.5285E-14$)&($9.8557E-14$)&($8.7444E-14$)&($9.3163E-14$)\\
$6$&$26$&$26$&$29$&${\bf 35}$\\
&($8.3589E-14$)&($7.7189E-14$)&($7.3281E-14$)&($7.9462E-14$)\\
$7$&$27$&$27$&$33$&$35$\\
&($8.6403E-14$)&($8.1151E-14$)&($7.2215E-14$)&($8.4479E-14$)\\
$8$&$27$&$25$&$29$&$37$\\
&($9.7379E-14$)&($8.0842E-14$)&($7.3955E-14$)&($7.4142E-14$)\\
$9$&${\bf 23}$&${\bf 24}$&$30$&$41$\\
&($9.5276E-14$)&($8.3798E-14$)&($8.7013E-14$)&($7.2068E-14$)\\
$10$&$25$&$25$&${\bf 27}$&$35$\\
&($7.6433E-14$)&($8.3980E-14$)&($9.3658E-14$)&($8.5795E-14$)\\
\hline\hline
\end{tabular}
\end{center}
\caption{Classical solitary wave generation of (\ref{fsec317}) . Number of iterations required by MPE, RRE, VEA and TEA as function of $\kappa$ to achieve a residual error below $TOL=10^{-13}$. The residual error at the last computed iterate is in parenthesis. Without acceleration,  the \PM
method (\ref{mm2}), (\ref{mm3c}) with $\gamma=2$ requires $n=76$ iterations with a residual error $9.0925E-14$.}\label{tav2b}
\end{table}
\begin{enumerate}
\item[(a)] For $\kappa\geq 2$, all the methods improve the performance of the \PM method without acceleration (cf. Figure \ref{figuresw1}(b)). The reduction in the number of iterations varies in a range  $50-70\%$.
\item[(b)] In general, polynomial methods (MPE and RRE, which essentially behaves in an equivalent way) are more efficient than $\epsilon$-algorithms (with VEA slightly better than TEA). In the best cases, the improvement is about $70\%$ in the case of MPE and RRE and RRE, about $65\%$ with respect to VEA and about $54\%$ in the case of TEA. (However, one has to take into account that the cycle in the case of polynomial methods is $mw=\kappa+1$ and in the case of $\epsilon$-algorithms is $mw=2\kappa$; cf. Figure \ref{figuresw2}.)
\end{enumerate}

In the case of the AAM, the corresponding results are in Table \ref{tav3b}. Now, the role of the parameter $\kappa$ (or $mw$) is played by $nw$.
\begin{table}
\begin{center}
\begin{tabular}{|c|c|c|}
\hline\hline  $nw$&AA-I($nw$)&AA-II($nw$)\\\hline
$1$&$38$($8.4014E-14$)&$35$($4.8504E-14$)\\
$2$&$28$($4.9835E-14$)&$26$($5.6978E-14$)\\
$3$&$28$($5.5678E-14$)&$25$($6.2897E-14$)\\
$4$&$27$($3.8773E-14$)&$22$($1.4624E-14$)\\
$5$&$22$($4.4004E-14$)&$20$($6.5530E-14$)\\
$6$&$21$($7.9925E-14$)&$21$($2.4615E-14$)\\
$7$&$21$($2.3111E-14$)&$20$($5.3227E-14$)\\
$8$&$20$($8.0666E-14$)&$20$($2.7701E-14$)\\
$9$&$20$($4.5873E-14$)&$19$($9.6556E-14$)\\
$10$&$20$($2.7208E-14$)&$19$($6.8255E-14$)\\
\hline\hline
\end{tabular}
\end{center}
\caption{Classical solitary wave generation of (\ref{fsec317}) . Number of iterations required by AA-I and AA-II as function of $nw$ to achieve a residual error below $TOL=10^{-13}$. The residual error at the last computed iterate is in parenthesis. Without acceleration,  the \PM
method (\ref{mm2}), (\ref{mm3c}) with $\gamma=2$ requires $n=76$ iterations with a residual error $9.0925E-14$.}\label{tav3b}
\end{table}
The results show that the performance of the methods is essentially the same. The best results are obtained with $nw=8$ in the case of AA-I and $nw=9$ for the AA-II. On the other hand, as mentioned in \cite{ni,walkern}, the value of $nw$ cannot be too large, because of ill-conditioning. In this example, this was observed for AA-I when $nw=9, 10$. (The corresponding results in Table \ref{tav3b} were obtained by using standard preconditioning.) Finally, compared to the \PM method without acceleration, the reduction in the number of iterations is in range of $50-80\%$.

Since the implementation of the methods is different, a comparison between VEM and AAM should take into account several efficiency indicators. In our example, we have measured the performance by computing the residual error as function of the number of iterations (i.~e. comparing the best results of Tables \ref{tav2b} and \ref{tav3b}) and as function of the computational time. 
The comparison of the methods in terms of the number of iterations is illustrated in Figure \ref{figuresw2}(a). This shows, in semilogarithmic scale, the residual error as function of the number of iterations for the \PM method without acceleration (solid line)  and accelerated with the six selected techniques, implemented with the values of $\kappa$ and $nw$ that, according to Tables \ref{tav2b} and \ref{tav3b}, lead to the best number of iterations. For this example, the AA-I(8) and AA-II(9) give, for a tolerance of $TOL=10^{-13}$ in the residual error, a slightly smaller number of iterations than the (mostly equivalent) RRE(9) and MPE(9). The initially worse performance of VEA(10) and TEA(6) is corrected after the first cycle. For example, in the case of VEA(10), after this first cycle, Figure \ref{figuresw2}(a) shows that the reduction in the residual error is the fastest. 

A second comment concerns the computational efficiency. Figure \ref{figuresw2}(b) shows (again in semi-log scale) the residual error as function of the CPU time in seconds for the four VEM. According to this, VEA(10) is the most efficient, followed by MPE(9), TEA(6) and RRE(9). The comparison in CPU time of this last one (the worst one among the VEM in efficiency) with AA-I(8) and AA-II(9) is shown in Figure \ref{figuresw2}(c) and reveals the poor performance in computational time as the main drawback of the AAM for this case. (see the formulation and implementation described in Section \ref{se2} to attempt to give an explanation of it.)

\subsubsection{Numerical generation of generalized solitary waves of (\ref{fsec318})}
\begin{figure}[htbp]
\centering \subfigure[]{
\includegraphics[width=6.6cm]{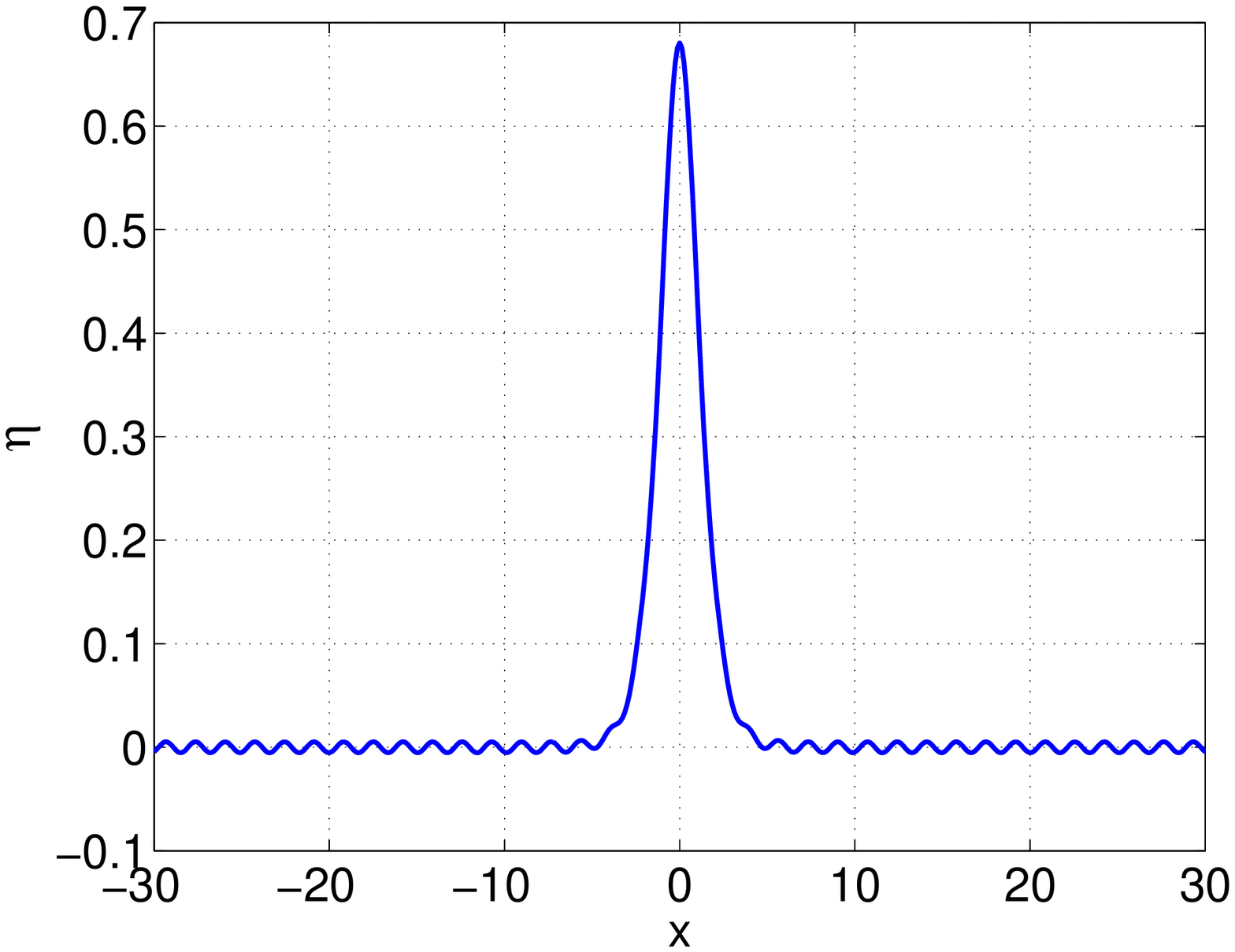} }
\subfigure[]{
\includegraphics[width=6.6cm]{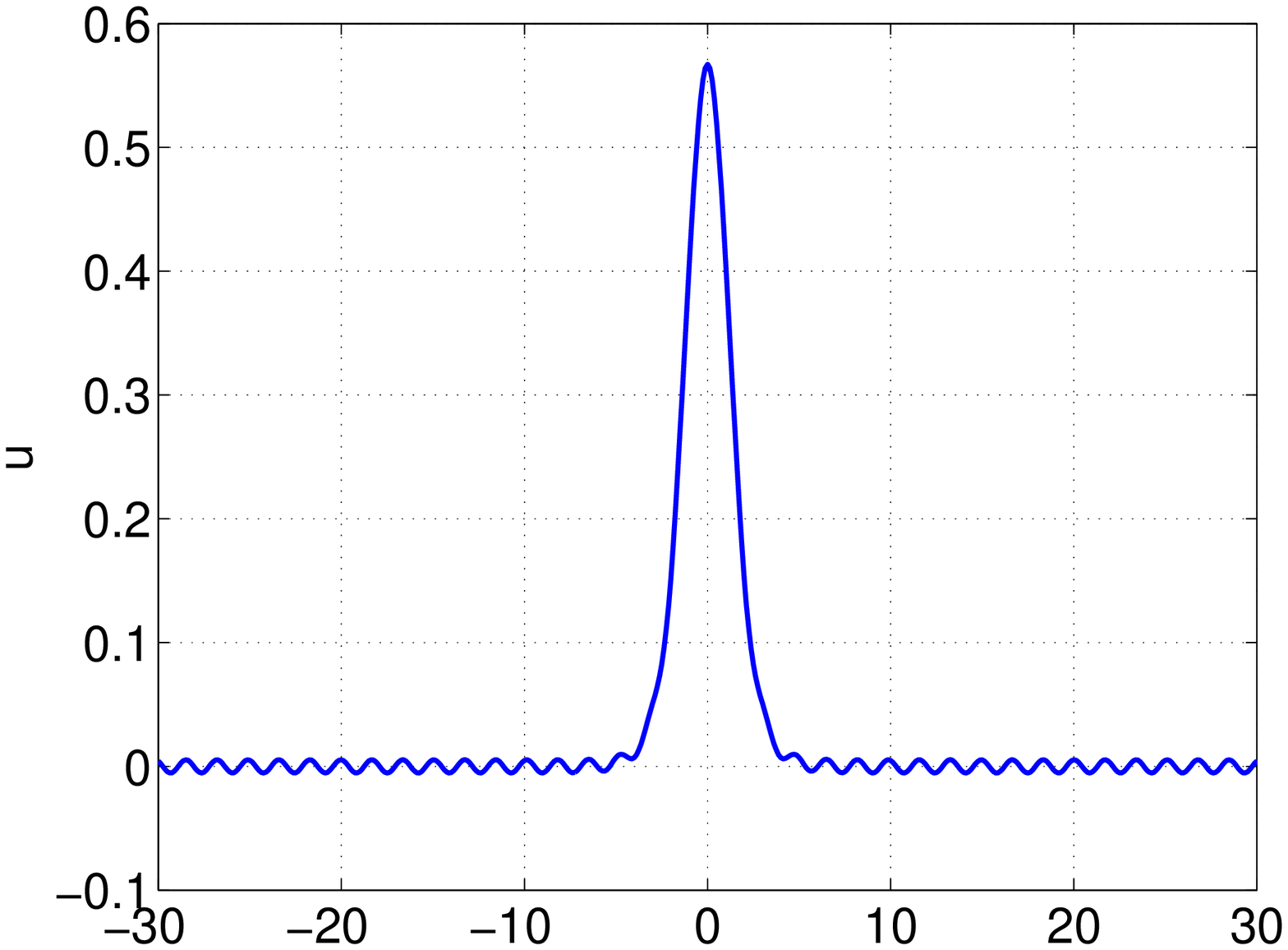} }
 \subfigure[]{
\includegraphics[width=6.6cm]{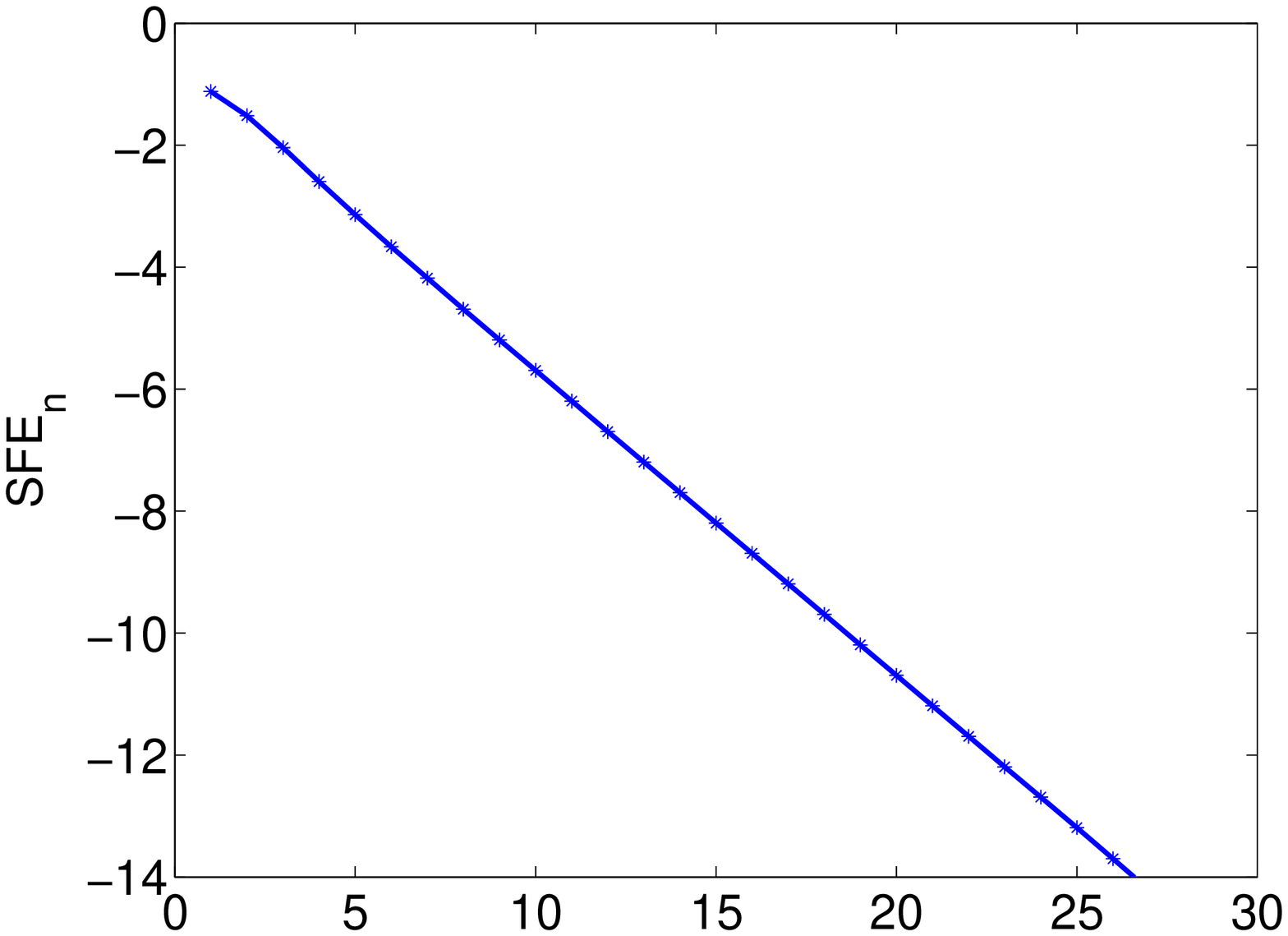} }
\subfigure[]{
\includegraphics[width=6.6cm]{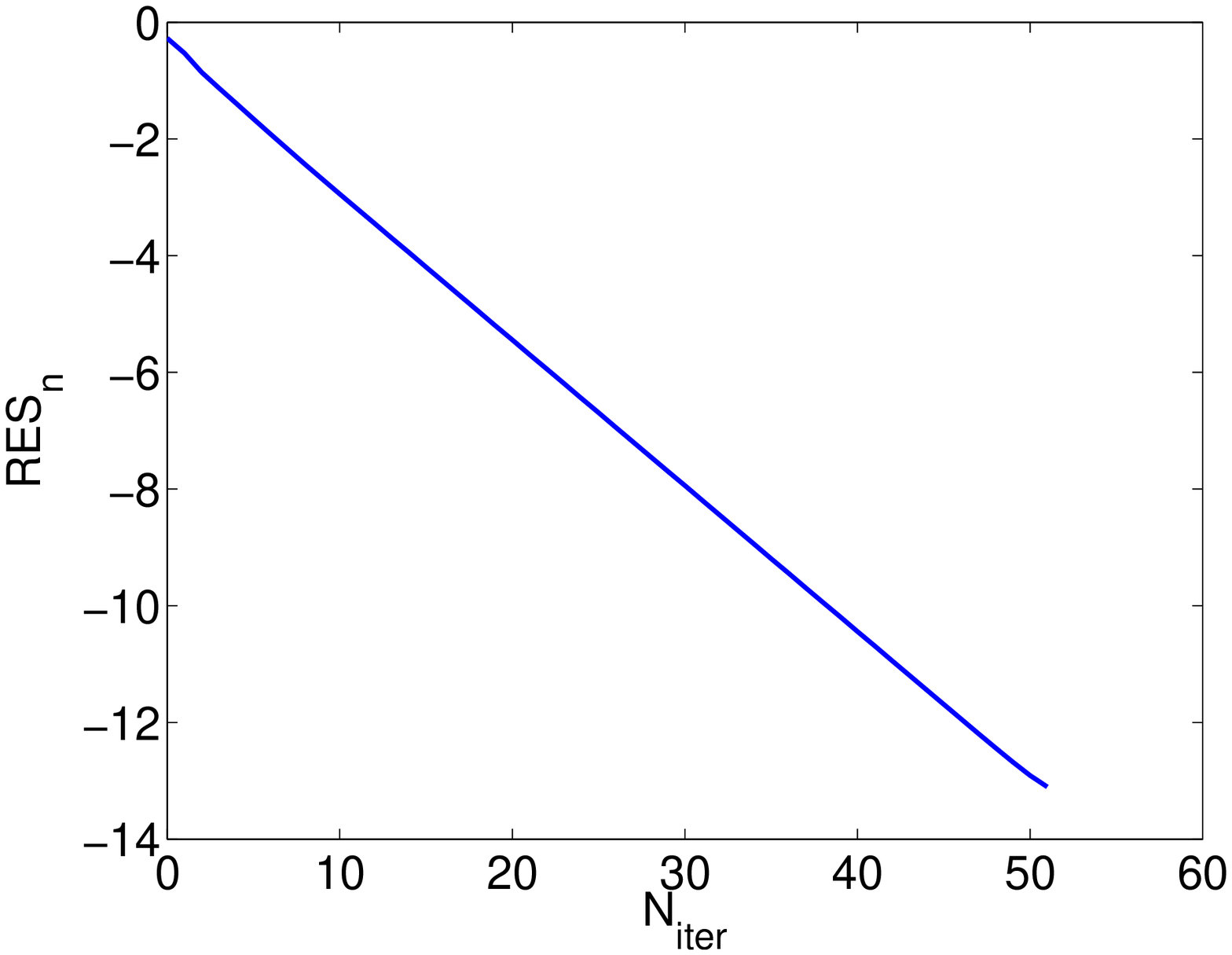} }
\caption{Generalized solitary wave generation of (\ref{fsec318}). (a)-(b) Approximate $\eta$ and $u$ profiles generated by the \PM method (\ref{mm2}), (\ref{mm3c}) for
(\ref{fsec318}); (c) Discrepancy (\ref{fsec33}) for the stabilizing
factor $s_{n}=s(\eta_{n},u_{n})$ vs number of iterations; (d)
Residual error (\ref{fsec32}) vs number of iterations. (Semi-logarithm scale in both cases.)} \label{figuresw3}
\end{figure}

\begin{figure}[htbp]
\centering
\subfigure[]{
\includegraphics[width=6.6cm]{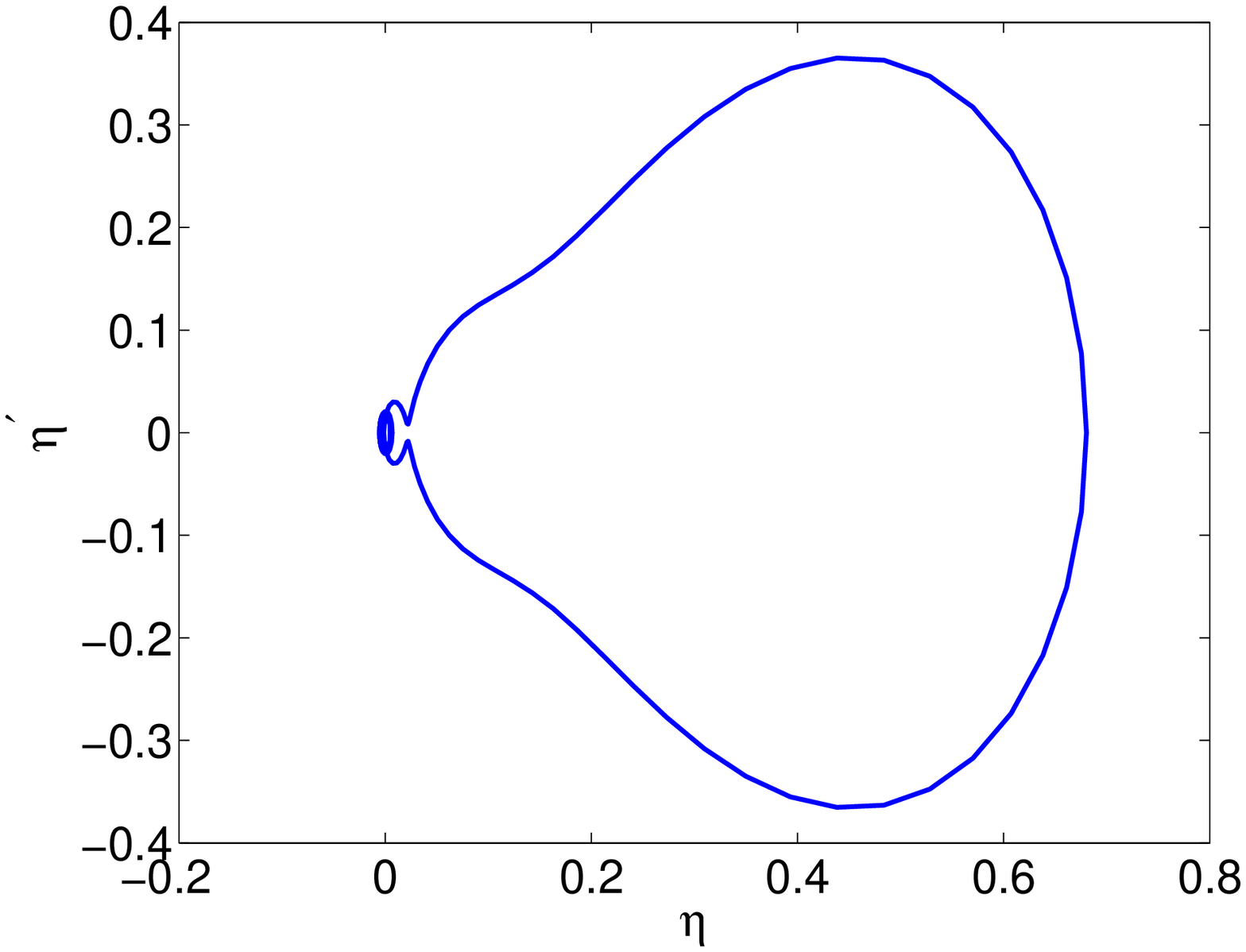} }
\subfigure[]{
\includegraphics[width=6.6cm]{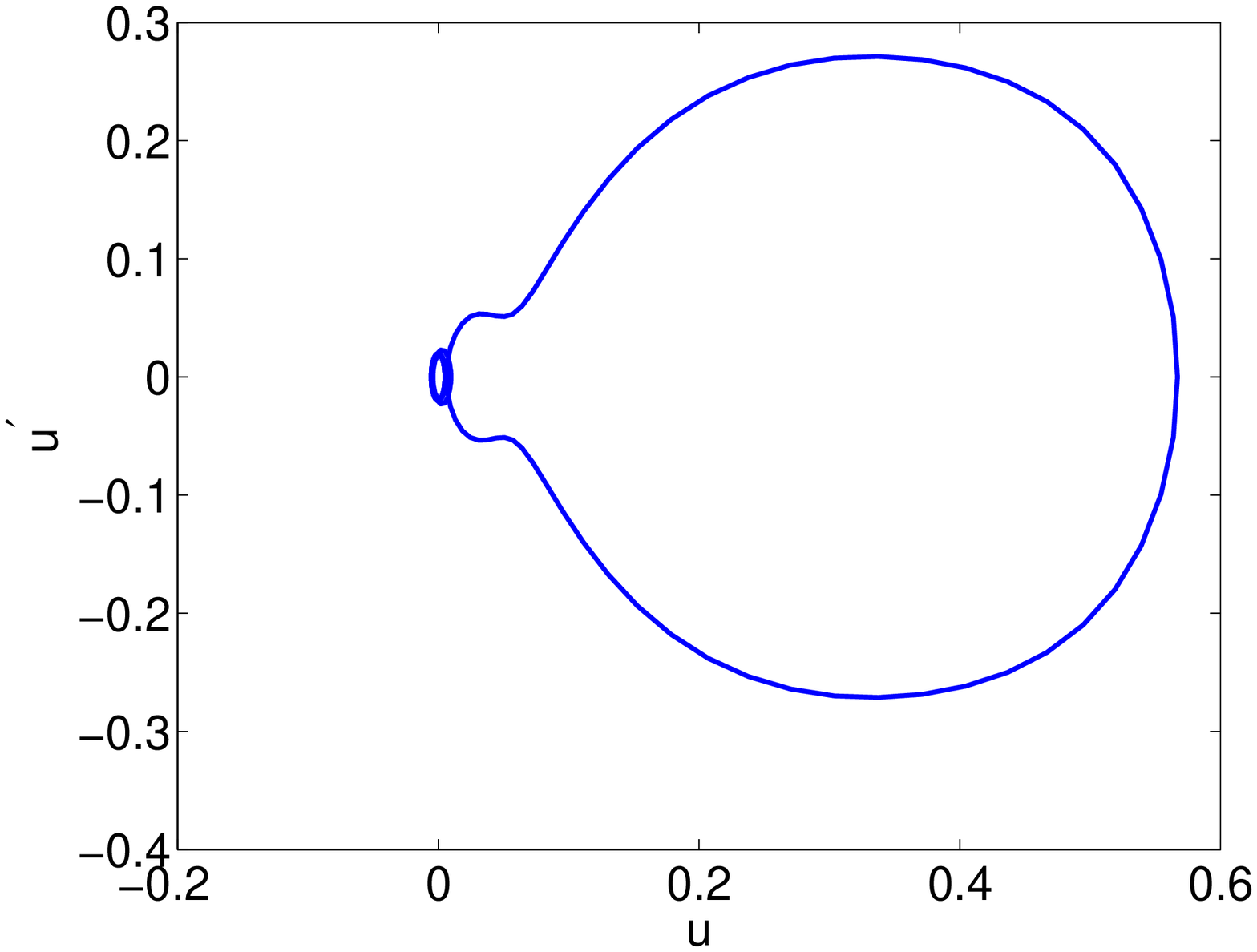} }
\caption{Generalized solitary wave generation of (\ref{fsec318}). (a) Phase portrait of $\eta$, (b) Phase portrait of $u$.} \label{figuresw3a}
\end{figure}

Here we show the results concerning the generation of approximate generalized solitary waves of the KdV-KdV system (\ref{fsec318}). In this case we have considered a speed $c_{s}=1.3$ and a Gaussian-type profile as initial guess for $\eta$ and $u$, with $l=64$ and $m=1024$ Fourier collocation points. The approximate $\eta$ and $u$ profiles generated by the \PM method (without acceleration) are displayed in Figures \ref{figuresw3}(a) and (b) respectively (observe the oscillatory ripples to the left and right of the main pulse), and the performance of the method (measured in terms of the convergence of the stabilizing factor and the behaviour of the residual error as function of number of iterations) is shown in Figures \ref{figuresw3}(c) and (d) respectively. The method achieves a residual error of $1.150546E-12$ in $n=47$ iterations and $7.859422E-14$ in $n=52$ iterations. In this case, the corresponding phase portraits in Figures \ref{figuresw3a}(a) and (b) show the generalized character of the waves, with orbits that are homoclinic to small amplitude periodic oscillations at infinity.
The last computed iterate, corresponding to this residual error, is used to evaluate the iteration matrices $S(\eta_{f},u_{f})$ and $F^{\prime}(\eta_{f},u_{f})$ of the classical fixed point and \PM method, respectively. The associated six largest magnitude eigenvalues are shown in Table \ref{tav4b}. (The generalized character of the computed solitary wave is also noticed by the presence of conjugate complex eigenvalues in the linearization matrix at the wave, cf. Table \ref{tav1b}.)

\begin{table}
\begin{center}
\begin{tabular}{|c|c|}
\hline\hline  Iteration matrix $S(\eta_{f},u_{f})$&
Iteration matrix $F^{\prime}(\eta_{f},u_{f})$\\\hline
$1.999999E+00$&$1.000000E+00$\\
$1.000000E+00$&$5.625613E-01$\\
$5.625613E-01$&$-3.525656E-01$\\
$-3.525656E-01$&$-3.521308E-01$\\
$-3.521308E-01$&$3.069304E-01+i 6.906434E-02$\\
$3.069304E-01+i 6.906434E-02$&$3.069304E-01-i 6.906434E-02$\\
\hline\hline
\end{tabular}
\end{center}
\caption{Generalized solitary wave generation of (\ref{fsec318}) . Six largest magnitude eigenvalues of the
approximated iteration matrix $S=L^{-1}N^{\prime}(\eta_{f},u_{f})$ (first column) and of the iteration matrix $F^{\prime}(\eta_{f},u_{f})$,  generated by the \PM
method (\ref{mm2}), (\ref{mm3c}) with $\gamma=2$, both evaluated
at the last computed iterate $(\eta_{f},u_{f})$.}\label{tav4b}
\end{table}

\begin{figure}[htbp]
\centering \subfigure[]{
\includegraphics[width=8.6cm]{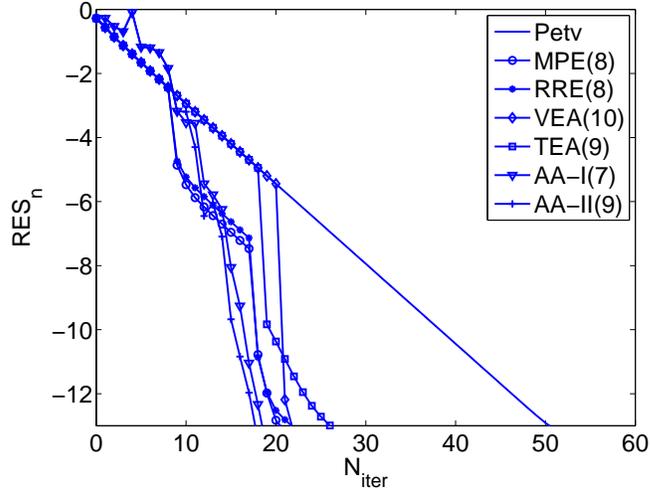} }
\subfigure[]{
\includegraphics[width=8.6cm]{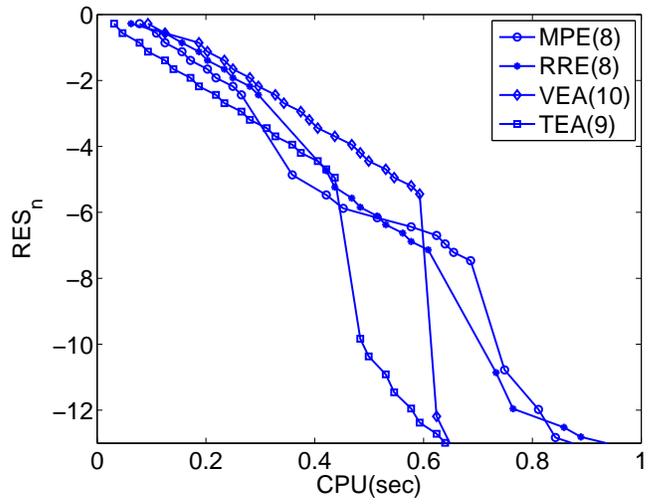} }
\subfigure[]{
\includegraphics[width=8.6cm]{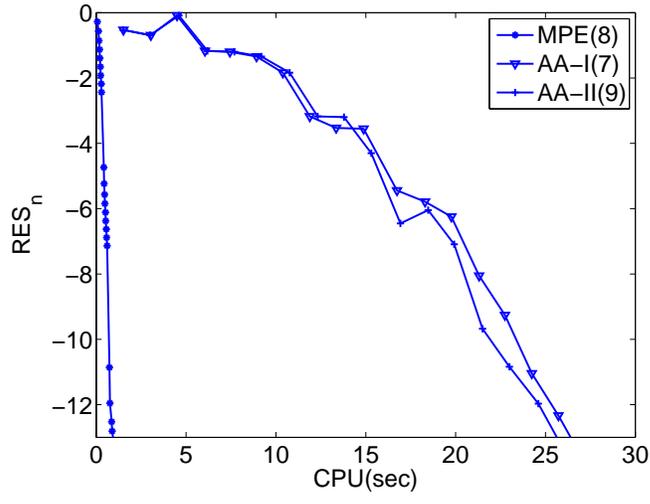} }
\caption{Convergence results of the \PM method for
(\ref{fsec318}): (a) Residual error (\ref{fsec32}) as function of number of iterations for the \PM method without acceleration and for six acceleration techniques with the best parameters $mw$ and $nw$ (according to Tables \ref{tav5b} and \ref{tav6b}). (b) Residual error (\ref{fsec32}) as function of CPU time (in seconds) for six acceleration techniques with the best parameters $mw$ and $nw$ (according to Tables \ref{tav5b} and \ref{tav6b}). (c) Comparison  of residual error (against CPU time) between the most efficient VEM and the AAM.
}
\label{figuresw4}
\end{figure}

The performance of the acceleration techniques is first checked in Tables \ref{tav5b} and \ref{tav6b} (respectively) which are the analogous to Tables  \ref{tav2b} and \ref{tav3b} respectively for the generalized solitary wave generation. The conclusions are the same as those of the generation of approximate classical solitary wave profiles for system (\ref{fsec317}): in terms of the number of iterations, AAM give the best performance and amongst the VEM, the extrapolation methods MPE and RRE are (in this case slightly) more efficient than the vector $\epsilon$-algorithms VEA and TEA, see Figure \ref{figuresw4}(a). The ranking is the opposite when residual error is measured in terms of the computational time. Figure \ref{figuresw4}(b) shows that VEA is the fastest and MPE the slowest. Even though, this is much faster than any of the Anderson algorithms, as shown in Figure \ref{figuresw4}(c). 
\begin{table}
\begin{center}
\begin{tabular}{|c|c|c|c|c|}
\hline\hline  $\kappa$&MPE($\kappa$)&RRE($\kappa$)&VEA($\kappa$)&TEA($\kappa$)\\\hline
$1$&$93$&$64$&$283$&$67$\\
&($8.2739E-14$)&($9.4330E-14$)&($6.7393E-14$)&($6.5919E-14$)\\
$2$&$42$&$43$&$81$&$49$\\
&($9.3934E-14$)&($8.6837E-14$)&($6.1237E-14$)&($8.7121E-14$)\\
$3$&$37$&$37$&$38$&$37$\\
&($9.8740E-14$)&($5.6549E-14$)&($6.7020E-14$)&($7.5320E-14$)\\
$4$&$33$&$32$&$39$&$30$\\
&($7.5531E-14$)&($6.9363E-14$)&($3.6501E-14$)&($4.2863E-14$)\\
$5$&$28$&$28$&$32$&$30$\\
&($6.6301E-14$)&($7.0496E-14$)&($5.9777E-14$)&($8.5244E-14$)\\
$6$&$28$&$28$&$29$&$28$\\
&($7.8419E-14$)&($7.4939E-14$)&($3.2861E-14$)&($7.7335E-14$)\\
$7$&$27$&$27$&$32$&$32$\\
&($3.0803E-14$)&($3.2406E-14$)&($6.5803E-14$)&($4.9439E-14$)\\
$8$&$23$&$23$&$28$&$29$\\
&($9.3872E-14$)&($8.8601E-14$)&($8.0707E-14$)&($9.1792E-14$)\\
$9$&$23$&$23$&$26$&$27$\\
&($3.0380E-14$)&($3.1758E-14$)&($6.4506E-14$)&($8.0379E-14$)\\
$10$&$24$&$24$&$24$&$28$\\
&($5.8031E-14$)&($5.0834E-14$)&($7.1291E-14$)&($5.8663E-14$)\\
\hline\hline
\end{tabular}
\end{center}
\caption{Generalized solitary wave generation of (\ref{fsec318}) . Number of iterations required by MPE, RRE, VEA and TEA as function of $\kappa$ to achieve a residual error below $TOL=10^{-13}$. The residual error at the last computed iterate is in parenthesis. Without acceleration,  the \PM
method (\ref{mm2}), (\ref{mm3c}) with $\gamma=2$ requires $n=52$ iterations with a residual error $7.8594E-14$.}\label{tav5b}
\end{table}

\begin{table}
\begin{center}
\begin{tabular}{|c|c|c|}
\hline\hline  $nw$&AA-I($nw$)&AA-II($nw$)\\\hline
$1$&$28$($7.6888E-14$)&$30$($3.6918E-14$)\\
$2$&$21$($3.2119E-14$)&$21$($3.4864E-14$)\\
$3$&$19$($4.9065E-14$)&$20$($4.6961E-14$)\\
$4$&$19$($2.5310E-14$)&$18$($1.9550E-14$)\\
$5$&$18$($5.9489E-14$)&$18$($4.2896E-14$)\\
$6$&$18$($2.0285E-14$)&$18$($1.1674E-14$)\\
$7$&$17$($5.1055E-14$)&$17$($5.5586E-14$)\\
$8$&$17$($2.8378E-14$)&$17$($1.9837E-14$)\\
$9$&$17$($7.2231E-14$)&$16$($7.1040E-14$)\\
$10$&$17$($2.8512E-14$)&$16$($5.7507E-14$)\\
\hline\hline
\end{tabular}
\end{center}
\caption{Generalized solitary wave generation of (\ref{fsec318}) . Number of iterations required by AA-I and AA-II as function of $nw$ to achieve a residual error below $TOL=10^{-13}$. The residual error at the last computed iterate is in parenthesis. Without acceleration,  the \PM
method (\ref{mm2}), (\ref{mm3c}) with $\gamma=2$ requires $n=53$ iterations with a residual error $7.8594E-14$.}\label{tav6b}
\end{table}

\subsubsection{Numerical generation of periodic traveling waves of (\ref{fsec319})}
\begin{figure}[htbp]
\centering 
\subfigure[]{
\includegraphics[width=6.6cm]{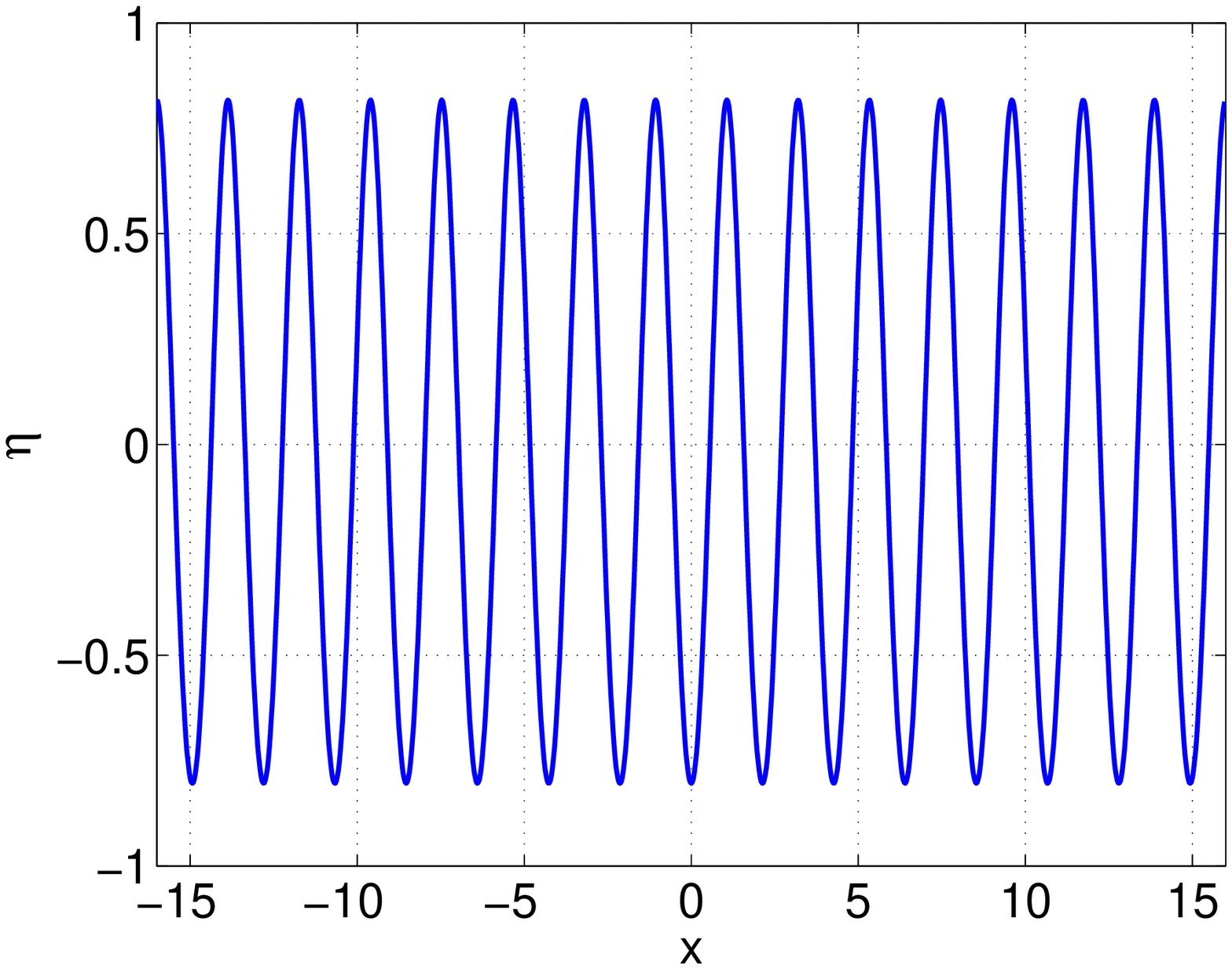} }
\subfigure[]{
\includegraphics[width=6.6cm]{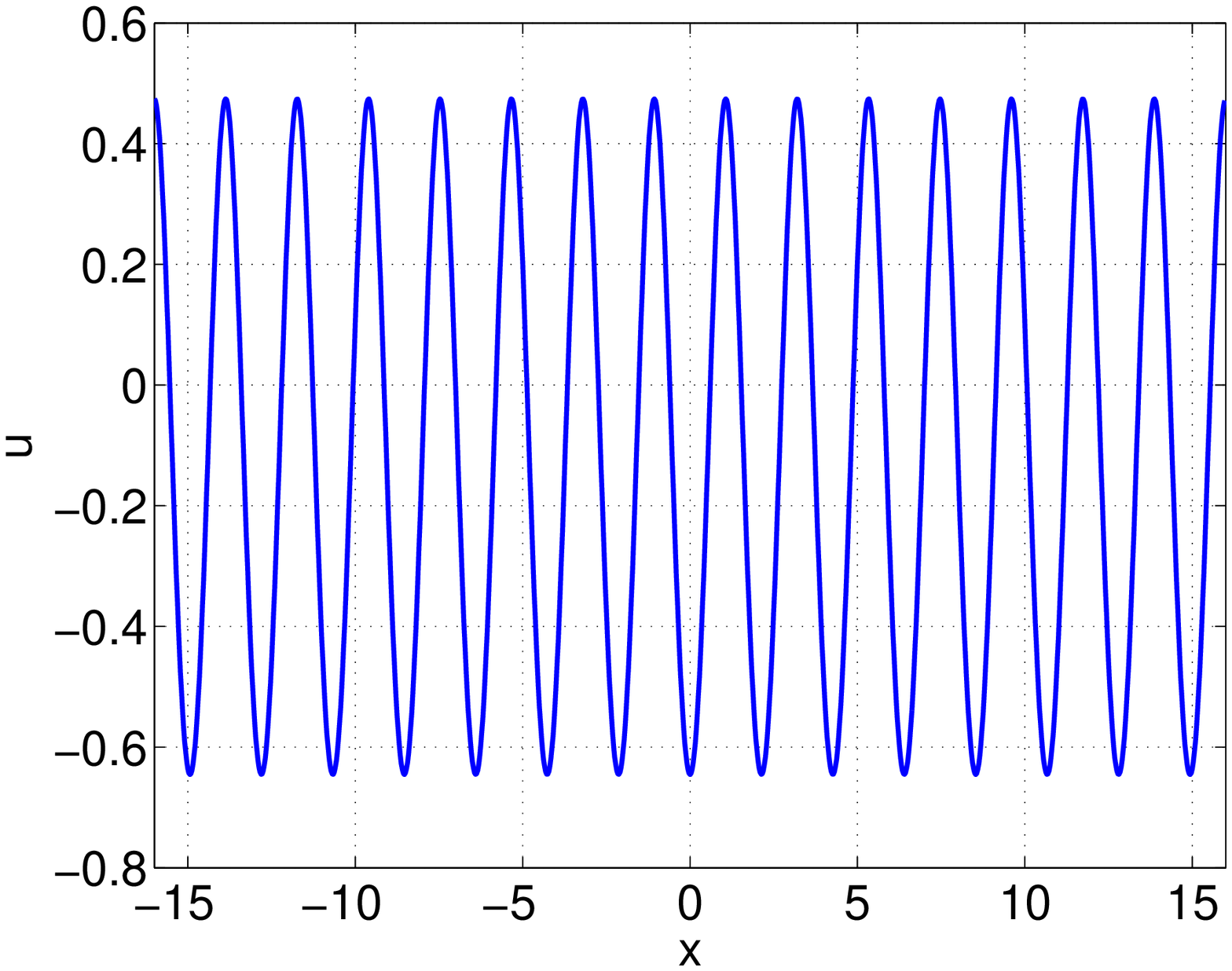} }
\subfigure[]{
\includegraphics[width=6.6cm]{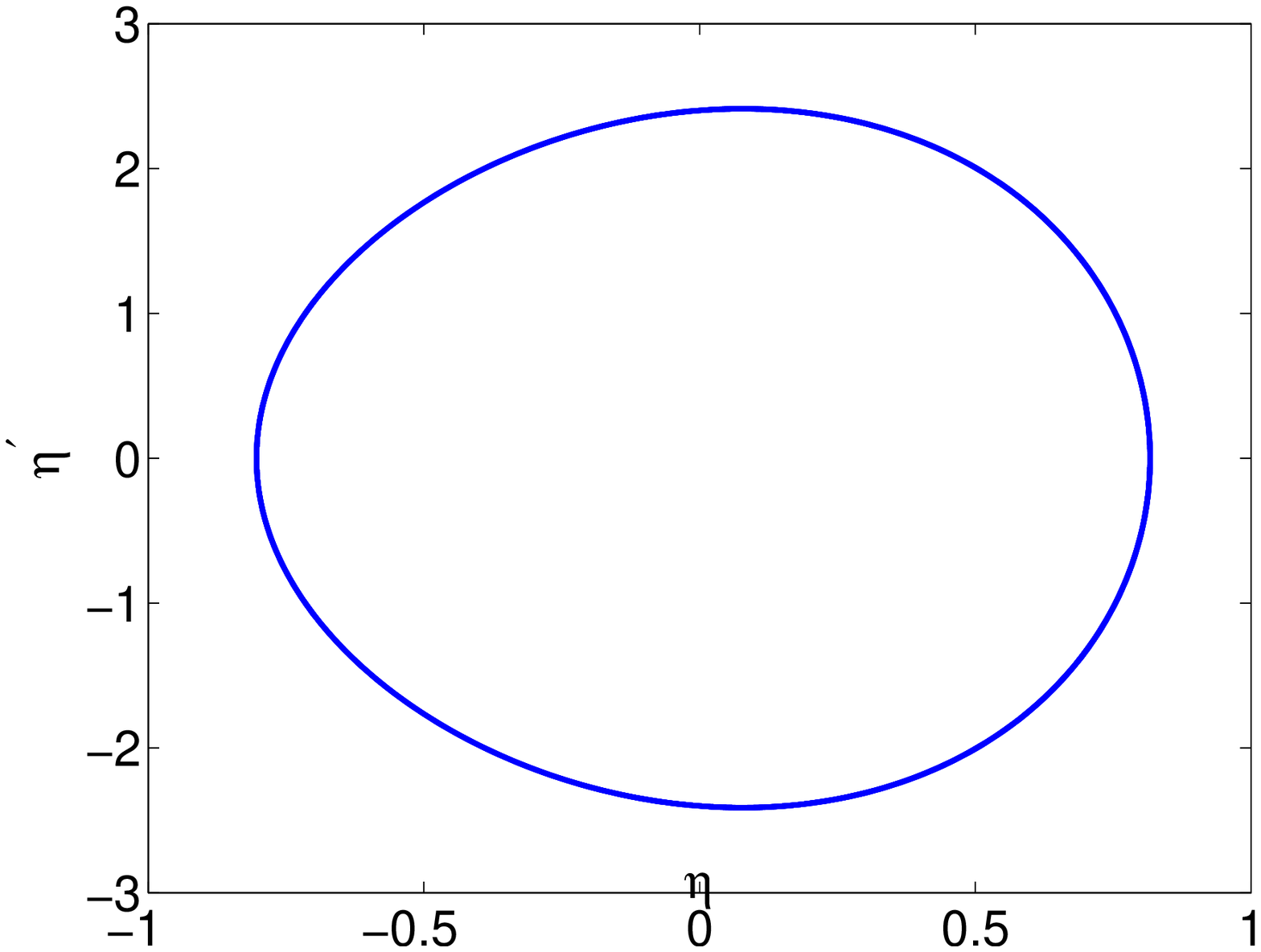} }
\subfigure[]{
\includegraphics[width=6.6cm]{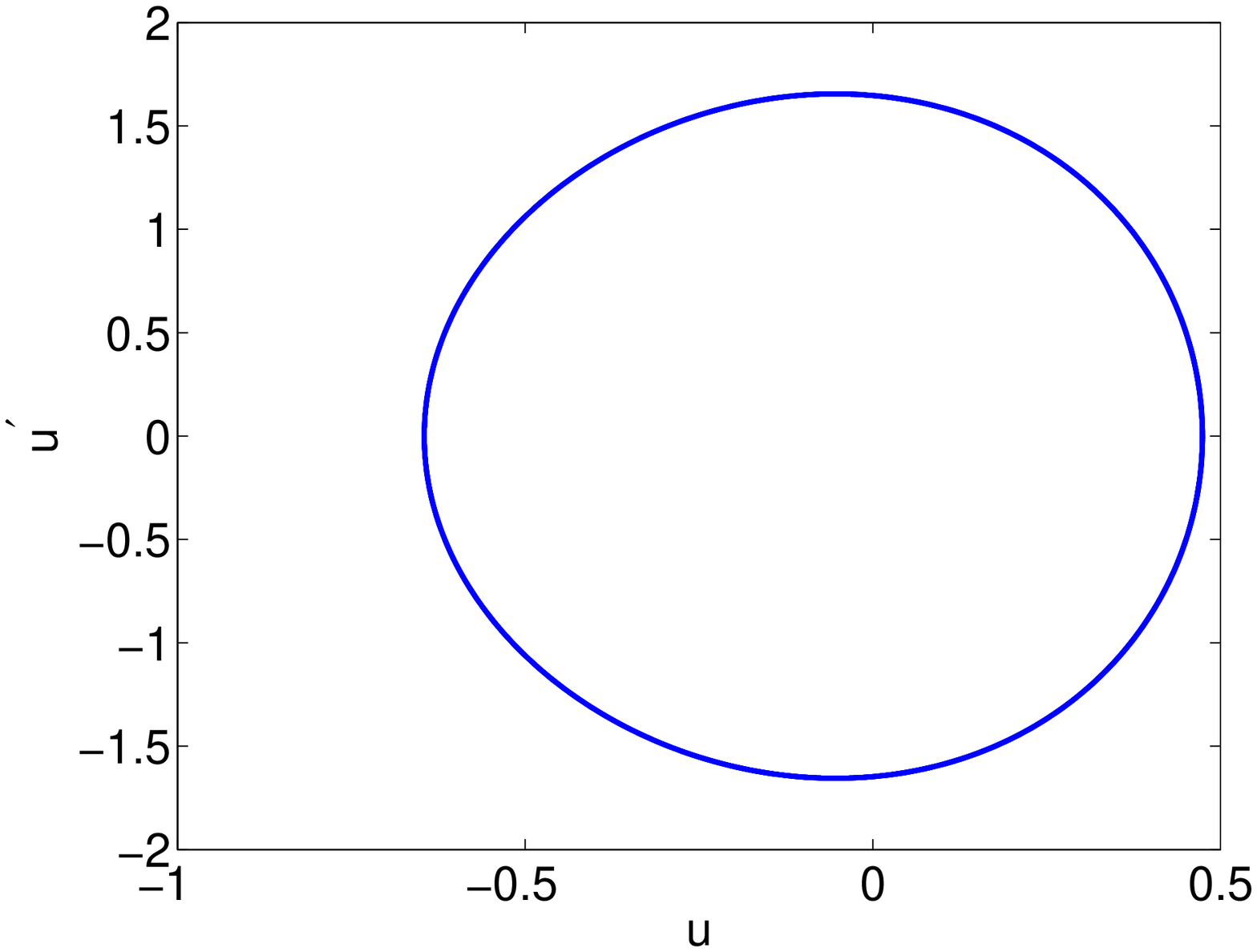} }
\caption{Numerical generation of periodic traveling waves of (\ref{fsec319}). Approximate profiles for $K_{1}=0.75, K_{2}=1$. (a) $\eta$ profile; (b) $u$ profile; (c) $\eta$ phase portrait. (d) $u$ phase portrait. } \label{Fig_ptw1}
\end{figure}

\begin{figure}[htbp]
\centering 
\includegraphics[width=9.6cm]{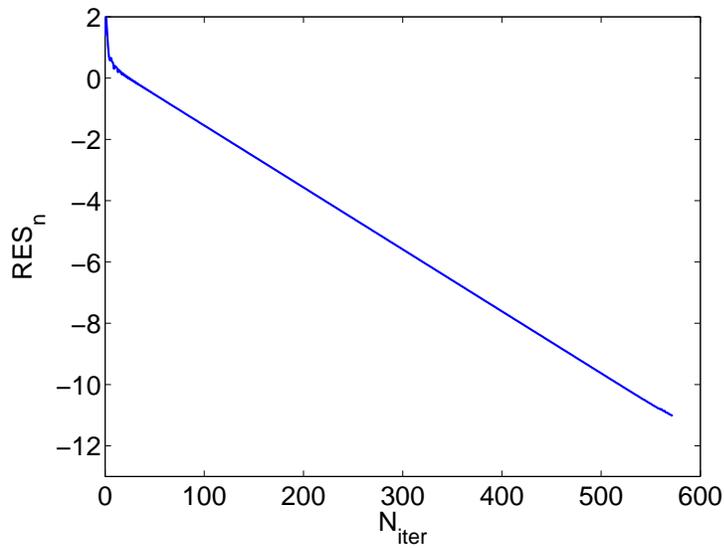}
\caption{Numerical generation of periodic traveling waves of (\ref{fsec319}). Approximate profiles for $K_{1}=0.75, K_{2}=1$. Residual error (\ref{fsec32}) vs number of iterations.} \label{Fig_ptw2}
\end{figure}

The numerical generation of periodic traveling wave solutions of the BBM-BBM system (\ref{fsec319}) will complete the study about traveling wave generation of Boussinesq systems (\ref{fsec311}). Here the initial data are similar to those of the previous cases, although now $l=16$ is taken. Once system (\ref{fsec313c}) is solved, the application of the \PM type method to (\ref{fsec314b}) generates, for $K_{1}=0.75, K_{2}=1$ (taken as an example) the computed profiles shown in Figure \ref{Fig_ptw1}(a)-(b).
The periodic behaviour is also observed in the corresponding phase plots, shown in Figure \ref{Fig_ptw1}(c), (d), while the performance is illustrated in Figure \ref{Fig_ptw2}, which corresponds to the behaviour of the residual error as function of the number of iterations. The method attains a residual error of $9.335366E-12$ in $n=572$  iterations, showing the need of some acceleration technique.

This slow performance is justified by the corresponding table of eigenvalues of the linearization operators, Table \ref{tavptw1} in this case.
\begin{table}
\begin{center}
\begin{tabular}{|c|c|}
\hline\hline  Iteration matrix $S(\eta_{f},u_{f})$&
Iteration matrix $F^{\prime}(\eta_{f},u_{f})$\\\hline
$2.000000E+00$&$1.000000E+00$\\
$1.000000E+00$&$-9.545242E-01$\\
$-9.545242E-01$&$-5.353103E-01-6.459204E-01i$\\
$-5.353103E-01-6.459204E-01i$&$-5.353103E-01-6.459204E-01i$\\
$-5.353103E-01-6.459204E-01i$&$-5.353103E-01+6.459204E-01i$\\
$-5.353103E-01+6.459204E-01i$&$-5.353103E-01+6.459204E-01i$\\
\hline\hline
\end{tabular}
\end{center}
\caption{Periodic traveling wave generation of (\ref{fsec319}) with $K_{1}=0.75, K_{2}=1$. Six largest magnitude eigenvalues of the
approximated iteration matrix $S=L^{-1}N^{\prime}(\eta_{f},u_{f})$ (first column) and of the iteration matrix $F^{\prime}(\eta_{f},u_{f})$,  generated by the \PM
method (\ref{mm2}), (\ref{mm3c}) with $\gamma=2$, both evaluated
at the last computed iterate $(\eta_{f},u_{f})$.}\label{tavptw1}
\end{table}
We observe that besides eigenvalue one (associated to the translational invariance) the next largest in magnitude eigenvalue is close to one. (As in the generalized solitary wave generation the presence of conjugate complex eigenvalues, in this case with algebraic multiplicity above one, in the spectrum of the linearization matrices is noticed.)

\begin{figure}[htbp]
\centering 
\subfigure[]{
\includegraphics[width=6.6cm]{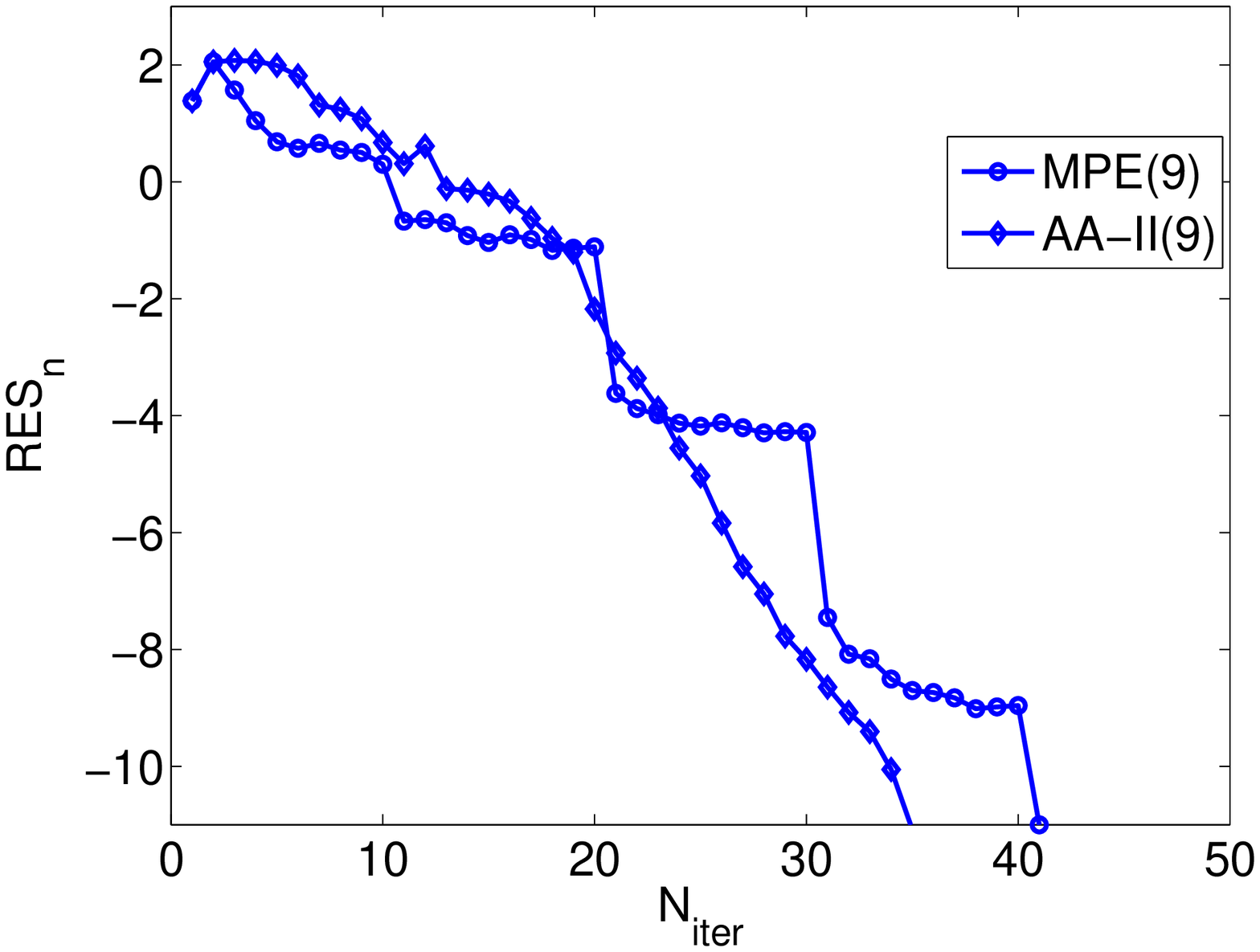} }
\subfigure[]{
\includegraphics[width=6.6cm]{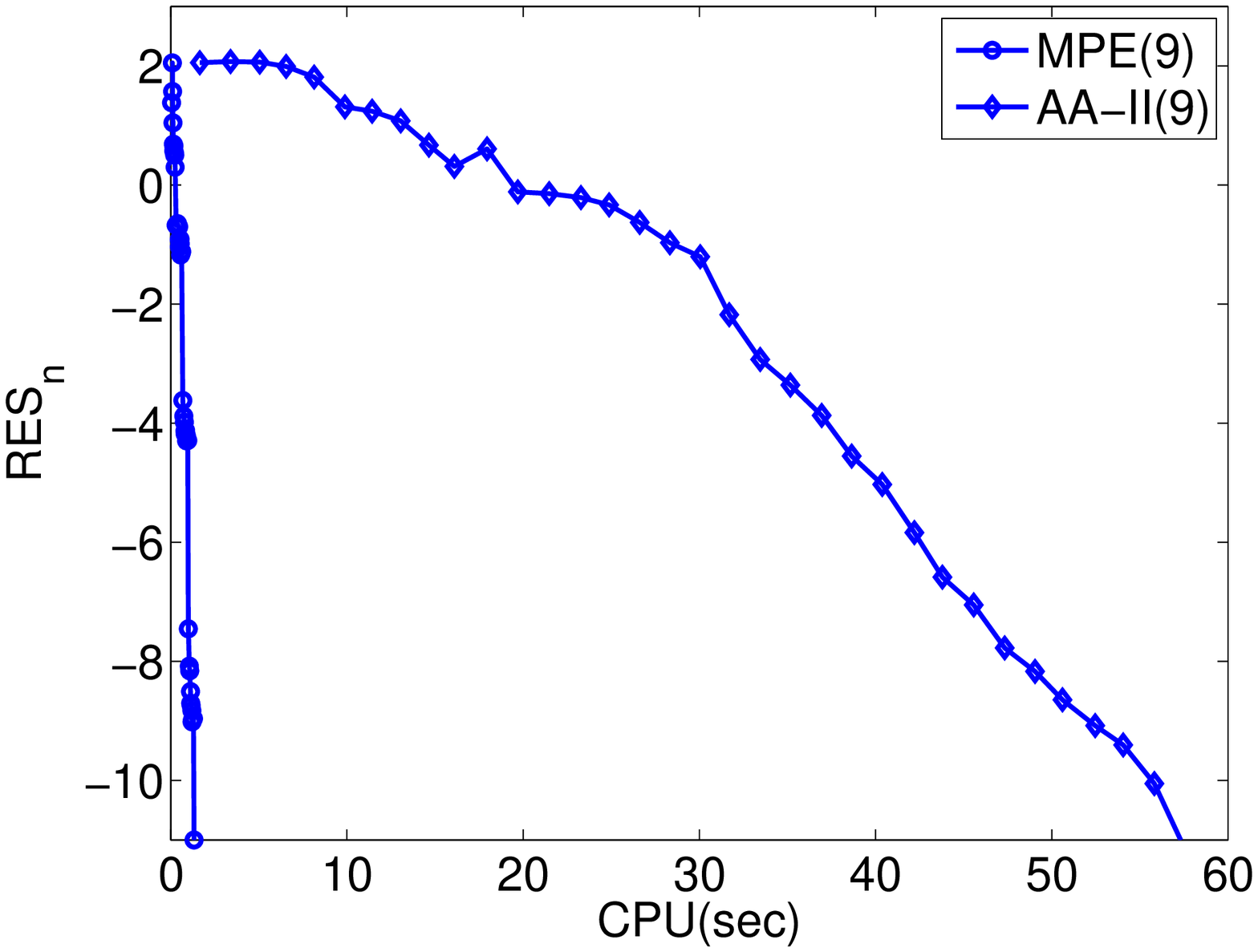} }
\caption{Numerical generation of periodic traveling waves of (\ref{fsec319}). Residual error (\ref{fsec32}) as function of: (a) number of iterations and (b) CPU time (in seconds) for MPE(9) and AA-II(9).} \label{Fig_ptw3}
\end{figure}

We now evaluate the application of VEM taking this example as illustration. The standard comparison in performance is given in Table \ref{tavptw2}. In this case the tolerance for the residual error was set as $TOL=10^{-11}$.
\begin{table}
\begin{center}
\begin{tabular}{|c|c|c|c|c|}
\hline\hline  $\kappa$&MPE($\kappa$)&RRE($\kappa$)&VEA($\kappa$)&TEA($\kappa$)\\\hline
$1$&$118$&&$278$&$88$\\
&($8.9607E-12$)&&($7.1073E-12$)&($4.9649E-12$)\\
$2$&$70$&$81$&$81$&$81$\\
&($9.2990E-12$)&($6.8950E-12$)&($6.8161E-12$)&($6.8798E-12$)\\
$3$&$54$&$53$&$78$&$64$\\
&($9.4403E-12$)&($7.2397E-12$)&($4.2025E-12$)&($9.4892E-12$)\\
$4$&$47$&$52$&$55$&$55$\\
&($7.0369E-12$)&($5.3720E-12$)&($4.0063E-12$)&($7.4587E-12$)\\
$5$&$55$&$46$&$48$&$103$\\
&($4.0773E-12$)&($8.2635E-12$)&($9.7571E-12$)&($6.2560E-12$)\\
$6$&$46$&$44$&$53$&$79$\\
&($9.9878E-12$)&($5.3034E-12$)&($5.7135E-12$)&($5.4689E-12$)\\
$7$&$45$&$49$&$51$&$73$\\
&($7.4645E-12$)&($3.6062E-12$)&($9.3703E-12$)&($8.3118E-12$)\\
$8$&$45$&$46$&$53$&$69$\\
&($9.1401E-12$)&($4.3655E-12$)&($9.0789E-12$)&($4.4378E-12$)\\
$9$&$41$&$41$&$58$&$77$\\
&($9.9555E-12$)&($4.7100E-12$)&($9.8141E-12$)&($8.4583E-12$)\\
$10$&$45$&$45$&$64$&64\\
&($3.9218E-12$)&($4.4096E-12$)&($6.2339E-12$)&($5.0514E-12$)\\
\hline\hline
\end{tabular}
\end{center}
\caption{Periodic traveling wave generation of (\ref{fsec319}). Number of iterations required by MPE, RRE, VEA and TEA as function of $\kappa$ to achieve a residual error (\ref{fsec32}) below $TOL=10^{-11}$. The residual error at the last computed iterate is in parenthesis. Without acceleration,  the \PM
method (\ref{mm2}), (\ref{mm3c}) with $\gamma=2$ requires $n=572$ iterations with a residual error $9.3354E-12$.}\label{tavptw2}
\end{table}
Some conclusions from it are the following:
\begin{enumerate}
\item Better performance of polynomial methods compared to $\epsilon$-algorithms.
\item MPE and RRE are virtually equivalent, especially when $\kappa$ grows. There are more differences between VEA and TEA, but they decrease when $\kappa$ grows.
\item For polynomial methods, the best results are obtained for large $\kappa$ (around $\kappa=9$) while for $\epsilon$ algorithms, it is better to take small $\kappa$ (around $mw=4, 5$). This implies a similar length of each cycle (width of extrapolation).
\end{enumerate}
We now analyze the results corresponding to AAM by using Table \ref{tavptw3}, which evaluates the performance of AA-I and AA-II for the same example.
\begin{table}
\begin{center}
\begin{tabular}{|c|c|c|}
\hline\hline  $nw$&AA-I($nw$)&AA-II($nw$)\\\hline
$1$&$57$($2.2757E-12$)&$66$($8.3451E-12$)\\
$2$&$84$($4.7115E-12$)&$40$($6.4724E-12$)\\
$3$&$81$($6.2384E-12$)&$39$($3.1769E-12$)\\
$4$&$81$($1.4968E-12$)&$36$($3.2197E-12$)\\
$5$&$48$($9.8634E-12$)&$36$($4.5059E-12$)\\
$6$&$48$($3.2036E-12$)&$49$($3.7796E-12$)\\
$7$&Ill-conditioned&$38$($2.4939E-12$)\\
$8$&Ill-conditioned&$36$($7.2839E-12$)\\
$9$&Ill-conditioned&$35$($8.3908E-12$)\\
$10$&Ill-conditioned&$37$($2.6179E-12$)\\
\hline\hline
\end{tabular}
\end{center}
\caption{Periodic traveling wave generation of (\ref{fsec319}). Number of iterations required by AA-I and AA-II as function of $nw$ to achieve a residual error (\ref{fsec32}) below $TOL=10^{-13}$. The residual error at the last computed iterate is in parenthesis. Without acceleration,  the \PM
method (\ref{mm2}), (\ref{mm3c}) with $\gamma=2$ requires $n=572$ iterations with a residual error $9.3354E-12$.}\label{tavptw3}
\end{table}
Some conclusions from Table \ref{tavptw3}:
\begin{enumerate}
\item As in some previous cases the AAM (particularly AA-II) behave better than any VEM when measuring the performance in terms of the number of iterations. However, the polynomial methods MPE and RRE are more efficient in terms of the computational time, see Figure \ref{Fig_ptw3}.
\item The best results of AA-II are obtained with large values of $nw$. The method does not appear to be affected by ill-conditioning, contrary to AA-I, which becomes useless from $nw=7$.
\end{enumerate}

\subsection{Example 2. Localized ground state generation}
A second group of experiments illustrates the generation of localized ground states in nonlinear Schr\"{o}dinger (NLS) type models with potentials. In particular, the equation
\begin{eqnarray}
\label{doub_well11}
    iu_{t}+\partial_{xx} u+V(x)u+|u|^{2}u=0,
    \end{eqnarray}
with potential $V(x)=6{\rm sech}^{2}(x)$, is considered as an example,
\cite{lakobay,yang2}. A localized ground state
solution of (\ref{doub_well11}) has the form $u(x,t)=e^{i\mu t}U(x)$, where $\mu\in
\mathbb{R}$ and $U(x)$ is assumed to be real and localized ($U\rightarrow 0,\; |x|\rightarrow\infty$). Substitution into (\ref{doub_well11}) leads to
\begin{eqnarray}
\label{doub_well12}
     U^{\prime\prime}(x)+V(x)U(x)-\mu U(x)+|U(x)|^{2}U(x)=0.
\end{eqnarray}
A discretization of (\ref{doub_well12}) based on a Fourier collocation method on a sufficiently long interval $(-l,l)$ requires in this case the resolution of a system of the form (\ref{mm1}) for the approximations $U_{h}$ of $U$ at the grid points $x_{j}=-l+jh, h=2l/m, j=0,\ldots,m-1$, with
\begin{eqnarray*}
L=D^{2}+{\rm diag}(V)-\mu I_{m},\quad N(U_{h})=-U_{h}.^{3},\label{s311}
\end{eqnarray*}
where $D$ is the pseudospectral differentiation matrix, ${\rm diag}(V)$ is the diagonal matrix with
elements $V_{j}=V(x_{j}),  j=0,\ldots,m-1$ and $I_{m}$ is the
$m\times m$ identity matrix. The nonlinearity $N$ is homogeneous with degree three, where, as usual, the dot
stands for the Hadamard product.
\begin{figure}[htbp]
\centering 
\subfigure[]{
\includegraphics[width=6.6cm]{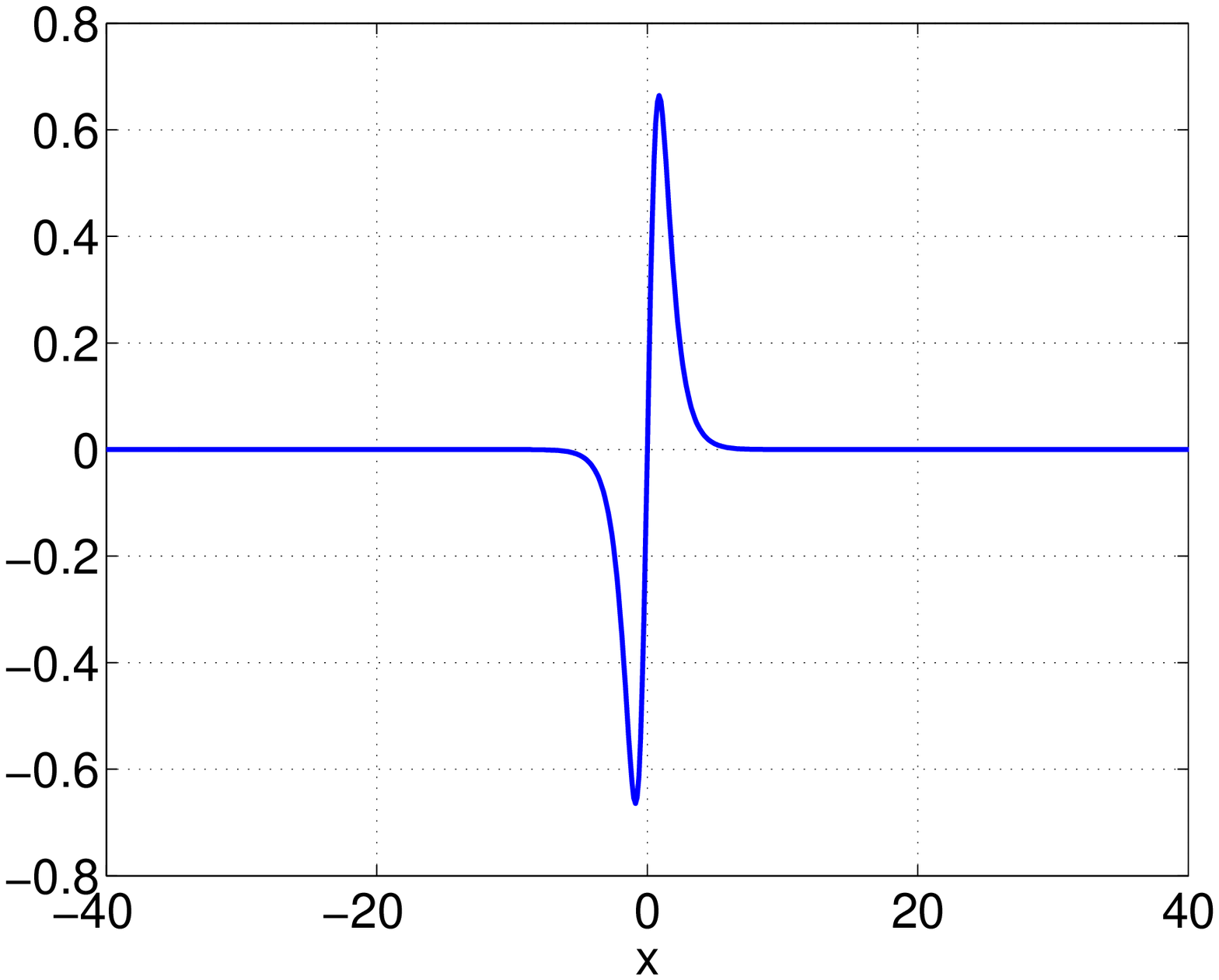} }
\subfigure[]{
\includegraphics[width=6.6cm]{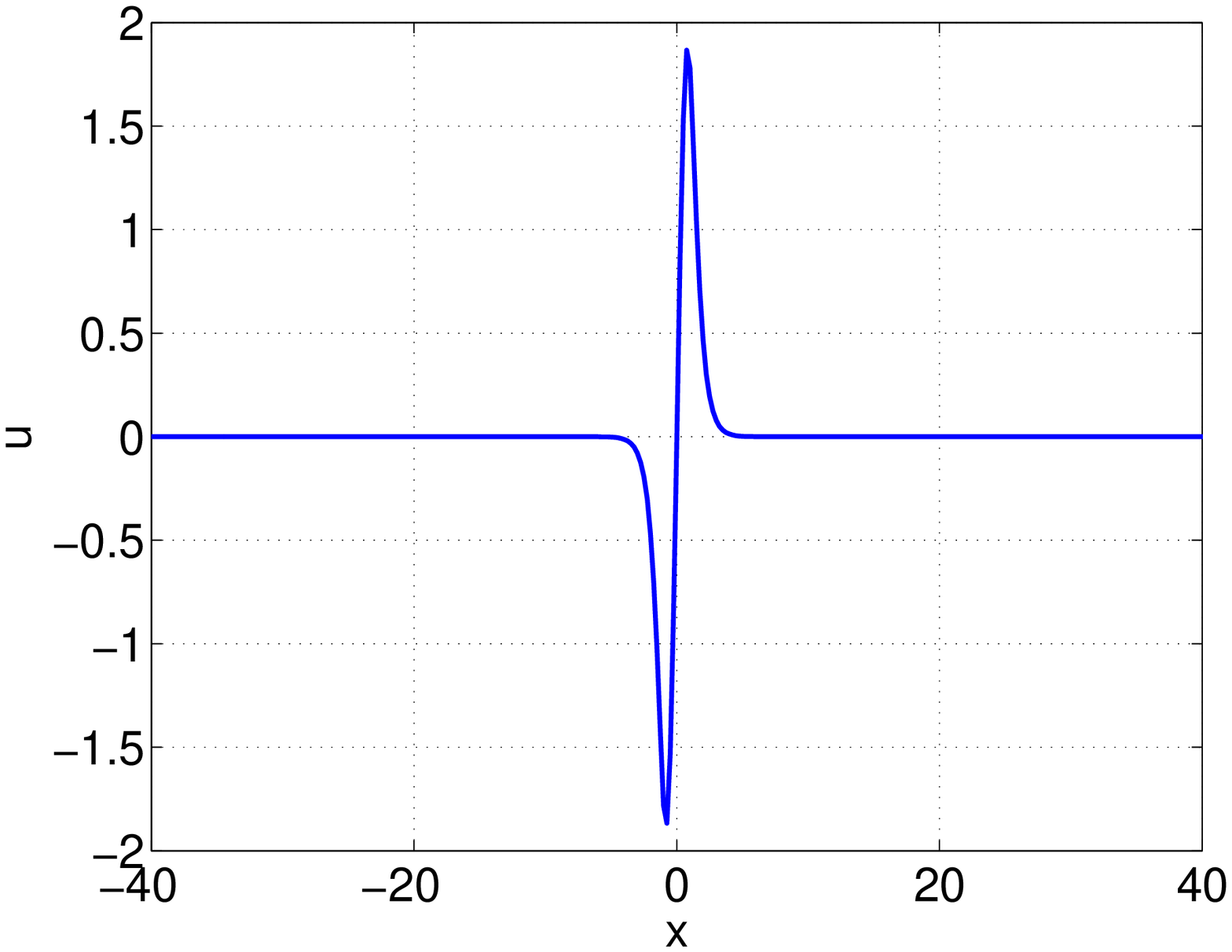} }
\subfigure[]{
\includegraphics[width=6.6cm]{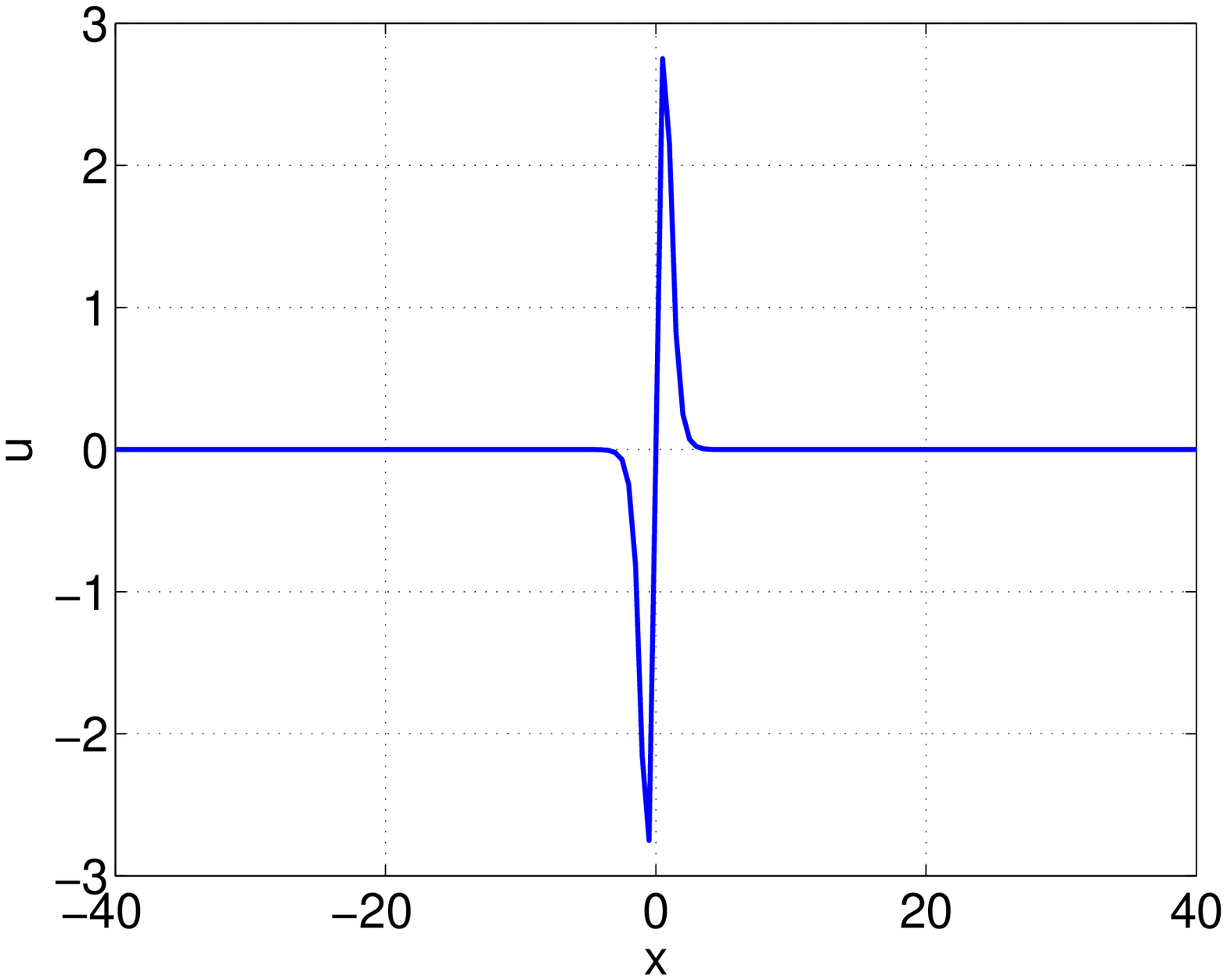} }
\subfigure[]{
\includegraphics[width=6.6cm]{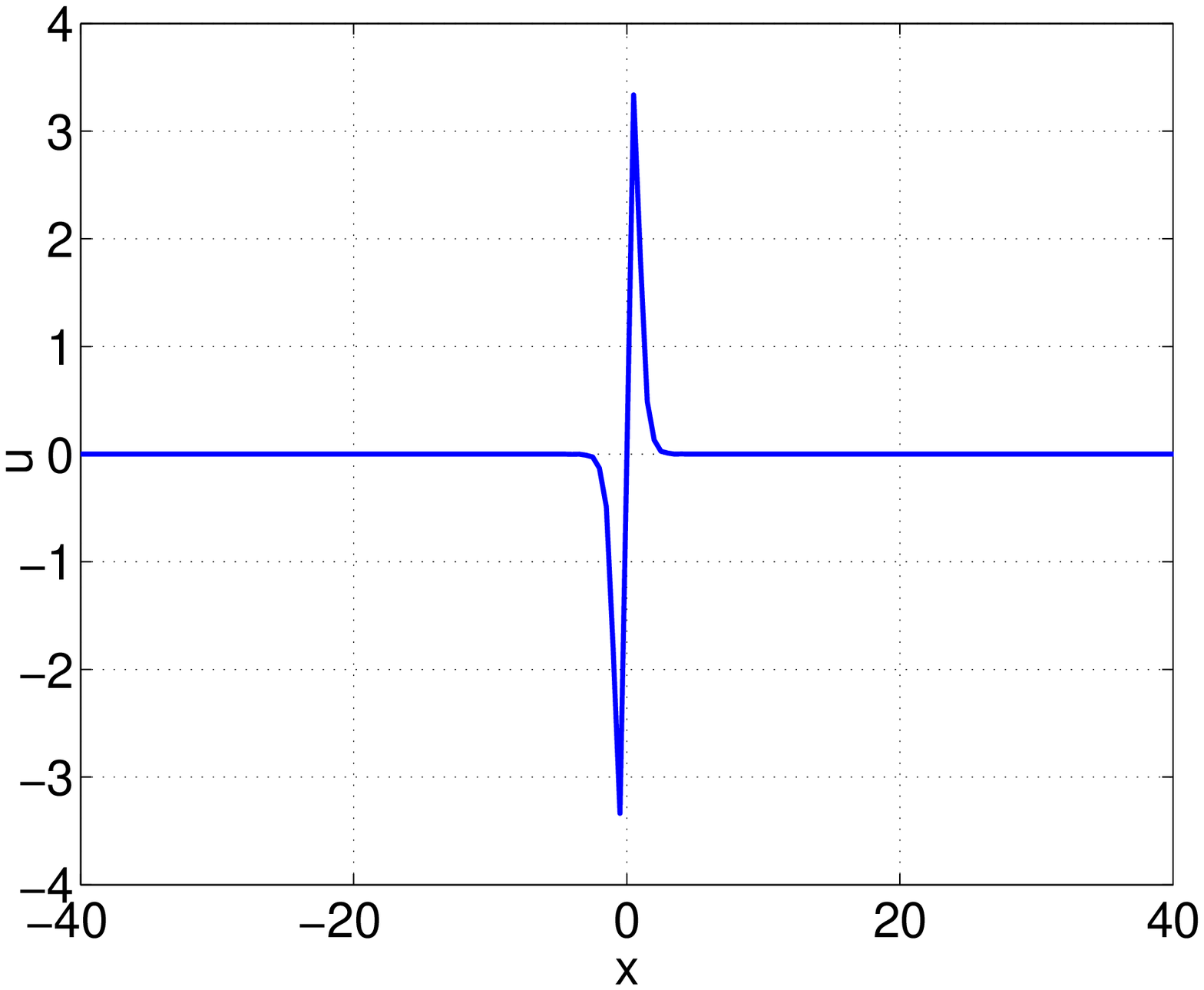} }
\caption{Numerical generation of localized ground states of (\ref{doub_well12}). Approximate asymmetric profile. (a) $\mu=1.3$; (b) $\mu=3.3$; (c) $\mu=6.3$; (d) $\mu=8.3$. The amplitude of the waves increases with $\mu$, while the shape is narrower.} \label{Fig321}
\end{figure}
The discussion below is focused on the ground state numerical generation for several values of $\mu$, which provide different challenges to the iteration. For each considered value of $\mu$, the performance of both families of acceleration techniques has been checked. 

\begin{figure}[htbp]
\centering 
\subfigure[]{
\includegraphics[width=6.6cm]{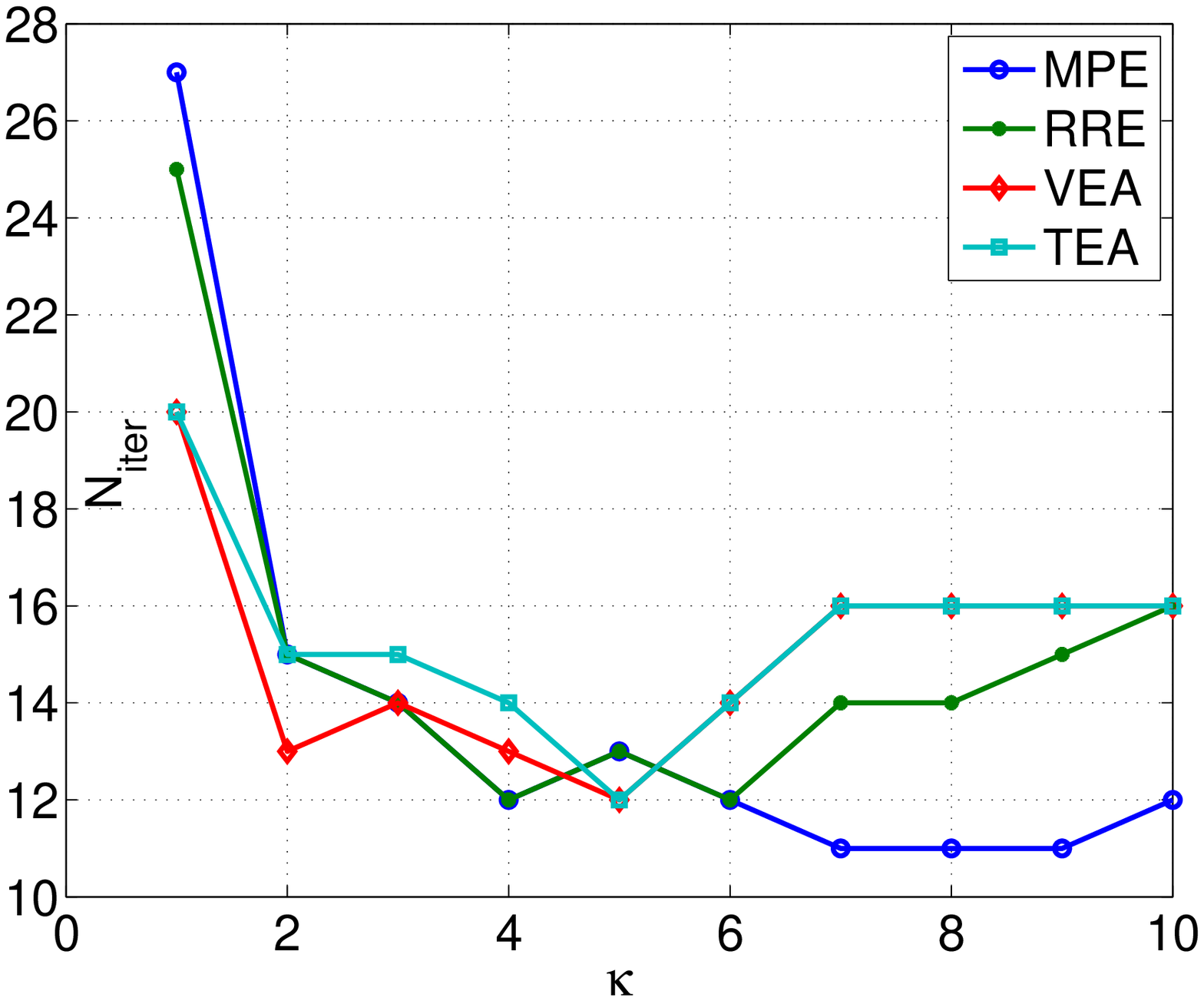} }
\subfigure[]{
\includegraphics[width=6.6cm,height=5.6cm]{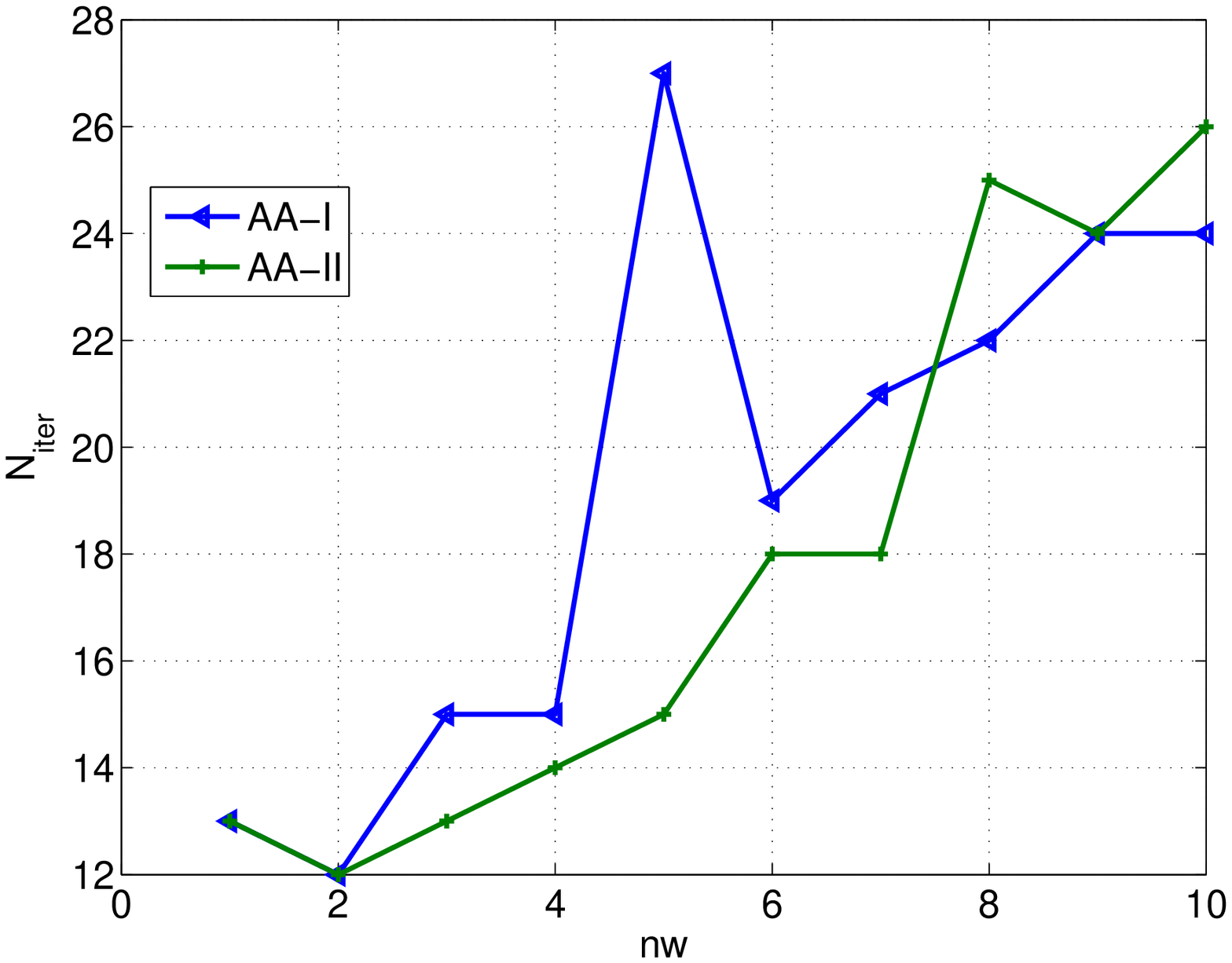} }
\caption{Numerical generation of asymmetric ground state of (\ref{doub_well12}) with $\mu=1.3$. Number of iterations required to reduce the residual error (\ref{fsec32}) below $TOL=10^{-12}$ and as function of the extrapolation width parameters $\kappa$ and $nw$. (a) VEM; (b) AAM.} \label{Fig322}
\end{figure}

\begin{figure}[htbp]
\centering 
\subfigure[]{
\includegraphics[width=8.6cm]{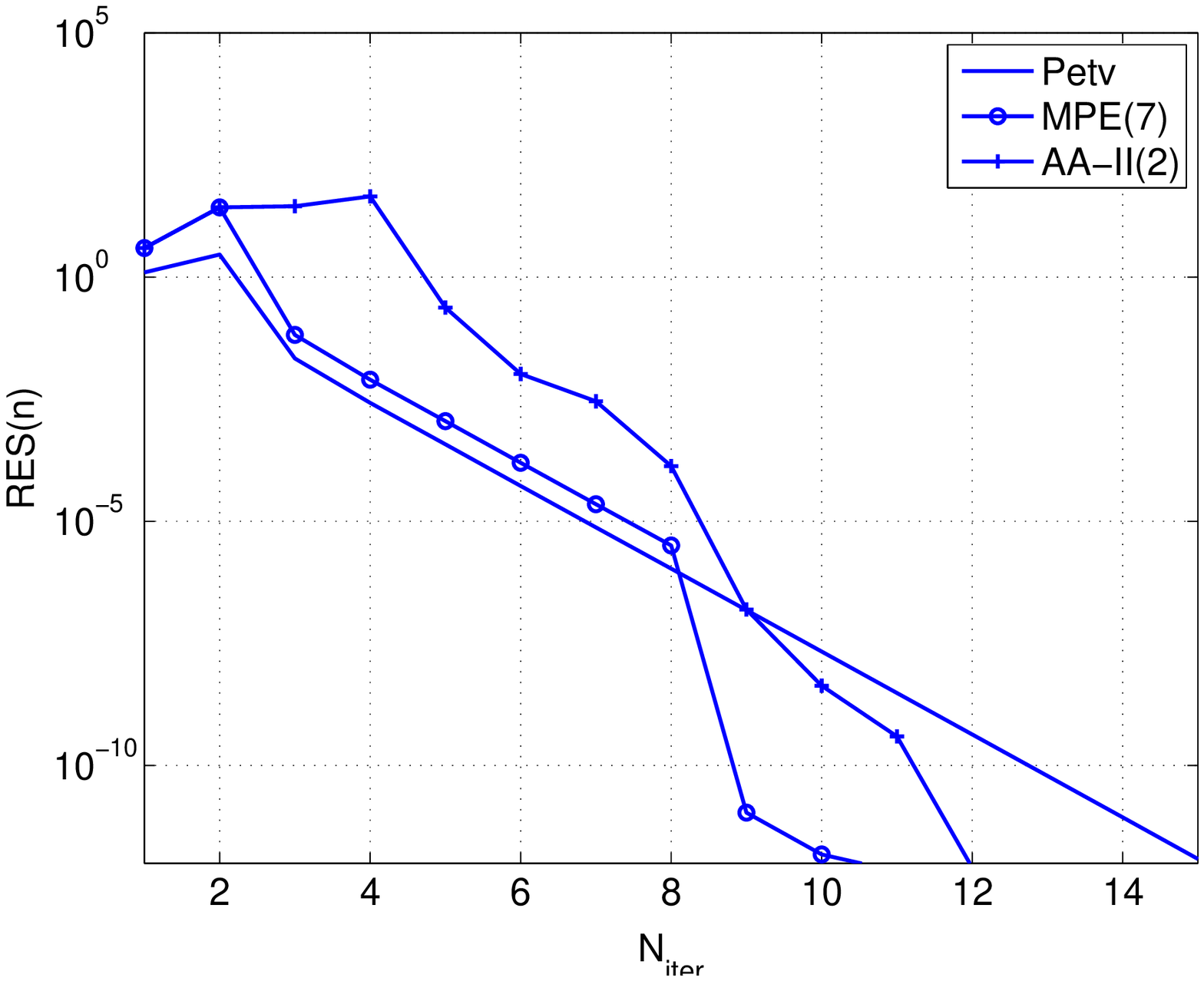} }
\subfigure[]{
\includegraphics[width=8.6cm]{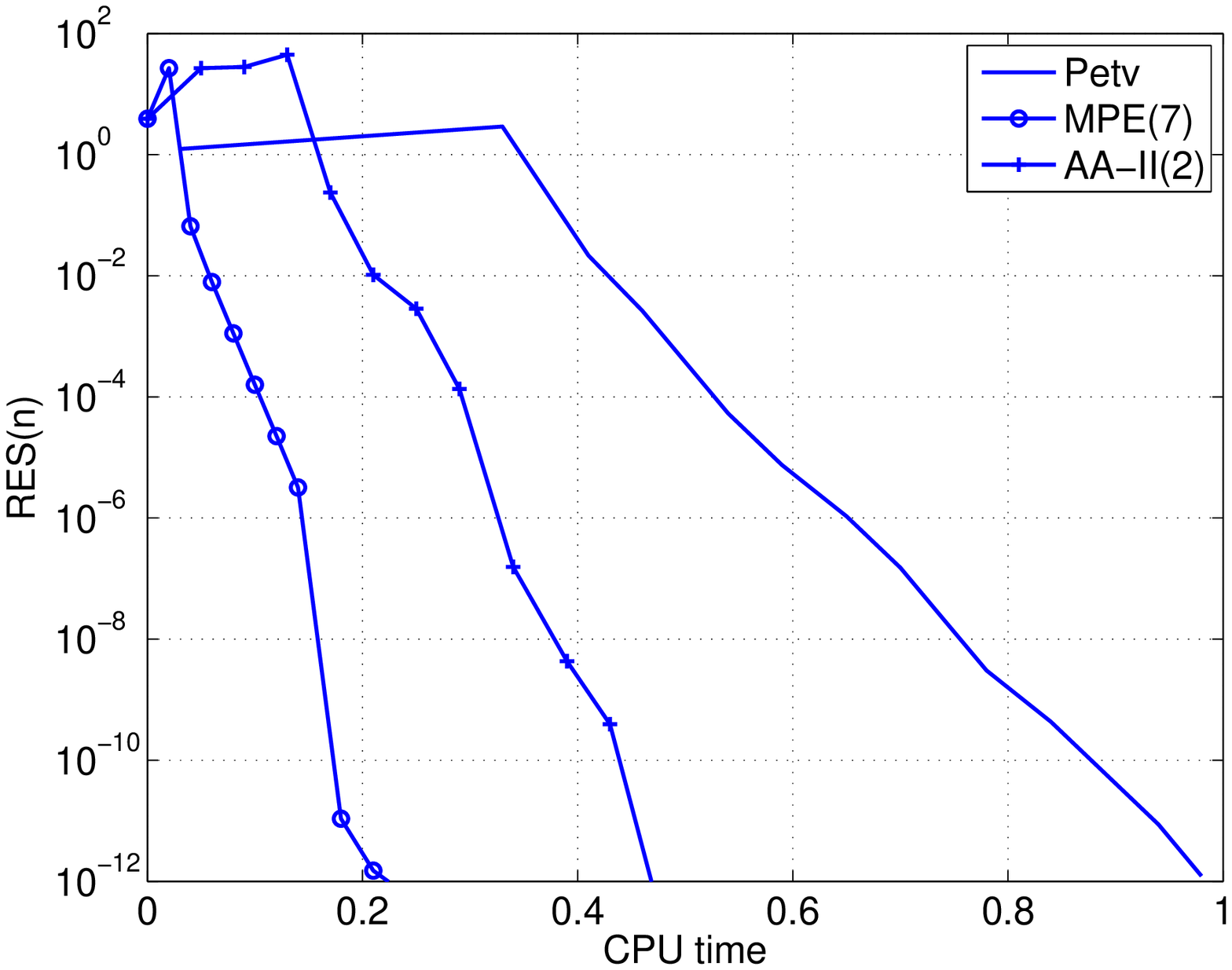} }
\caption{Numerical generation of assymetric ground state of (\ref{doub_well12}) with $\mu=1.3$. Residual error (\ref{fsec32}) as function of the number of iterations (a) and CPU time in seconds (b) for the \PM method without acceleration (solid line) and accelerated with MPE(7) (circle symbols) and AA-II(2) (plus symbols).} \label{Fig323}
\end{figure}

The first results concern the numerical generation of an asymmetric solution of (\ref{doub_well12}) for $\mu=1.3$ (Figure \ref{Fig321}(a)). Figure \ref{Fig322} compares the performance of the acceleration techniques in terms of the number of iterations required to reduce the residual error (\ref{fsec32}) below $TOL=10^{-12}$ and as function of the extrapolation width parameters $\kappa$ and $nw$. In the case of VEM, Figure \ref{Fig322}(a), all the techniques considered are comparable and the differences are not very large; MPE with $\kappa=7$ gives the minimum number of iterations. (The values $\kappa=8, 9$ also lead to the same number of iterations, but the computational effort in CPU time is higher.) The AAM, Figure \ref{Fig322}(b), are competitive with VEM for small values of $nw$ ($nw=1,2$). As $nw$ grows, the number of iterations increases (in opposite way to the behaviour of VEM with respect to $\kappa$) and the computation of the coefficients in the minimization problem becomes ill-conditioned. The comparison between the most efficient method of each family (MPE(7) and AA-II(2) respectively) is displayed in Figure \ref{Fig323}.
It shows the residual error as function of the number of iterations (a) and the CPU time (b) for the \PM method without acceleration (solid line) and accelerated with MPE(7) and AA-II(2).  In both figures, the improvement in the performance with respect to the \PM method provided by the two acceleration techniques is observed, with the best results corresponding to MPE (and, in general VEM against AAM). 
\begin{table}
\begin{center}
\begin{tabular}{|c|c||c|c|}
\hline\hline
$\mu=1.3$&$\mu=1.3$&$\mu=3.3$&$\mu=3.3$\\\hline\hline
eigs$(S)$& eigs$(F^{\prime}(u^{*}))$&eigs$(S)$& eigs$(F^{\prime}(u^{*}))$\\\hline
2.999999E+00&2.886842E-01&-6.328271E+00&-6.328271E+00\\
2.886842E-01&-1.858331E-01&3.000000E+00&8.594730E-01\\
-1.858331E-01&1.419117E-01&8.594730E-01&5.552068E-01\\
1.419117E-01&7.522396E-02&5.552068E-01&2.978699E-01\\
7.522396E-02&5.527593E-02&2.978699E-01&2.360730E-01\\
5.527593E-02&3.629863E-02&2.360730E-01&1.552434E-01\\
\hline\hline
\end{tabular}
\end{center}
\caption{Numerical generation of asymmetric profile of (\ref{doub_well12})
with $\mu=1.3$ and $\mu=3.3$. Six largest magnitude eigenvalues of the
approximated iteration matrix of the classical fixed-point method $S=L^{-1}N^{\prime}(U_{f})$ (left)  and of the \PM method, evaluated
at the last computed iterate $U_{f}$ obtained with MPE($7$).}\label{tav1}
\end{table}
Table \ref{tav1} confirms the convergence of the \PM method. It displays the six largest magnitude eigenvalues of the corresponding iteration matrix of the classical fixed-point algorithm $S=L^{-1}N^{\prime}(U_{f})$, and of the \PM method (\ref{iterop}), (\ref{mm3c}) for two values of $\mu$. Since an analytical expression for the exact profile is not known, the matrices have been evaluated at the last computed iterate given by MPE(7). In the case of $S$ (first column), the dominant eigenvalue corresponds to the degree of homogeneity $p=3$, with the rest of the eigenvalues below one. The filter action of the stabilizing factor is observed in the second column. The degree $p=3$ has been subtituted by zero (the optimal $q=\gamma (1-p)=-p$ has been taken) and the rest of the spectrum is preserved. This implies that for $\mu=1.3$, the spectral radius of $F^{\prime}(u^{*})$ is below one (second column) and this leads to the (local) convergence of the method.

\begin{figure}[htbp]
\centering 
\subfigure[]{
\includegraphics[width=6.6cm]{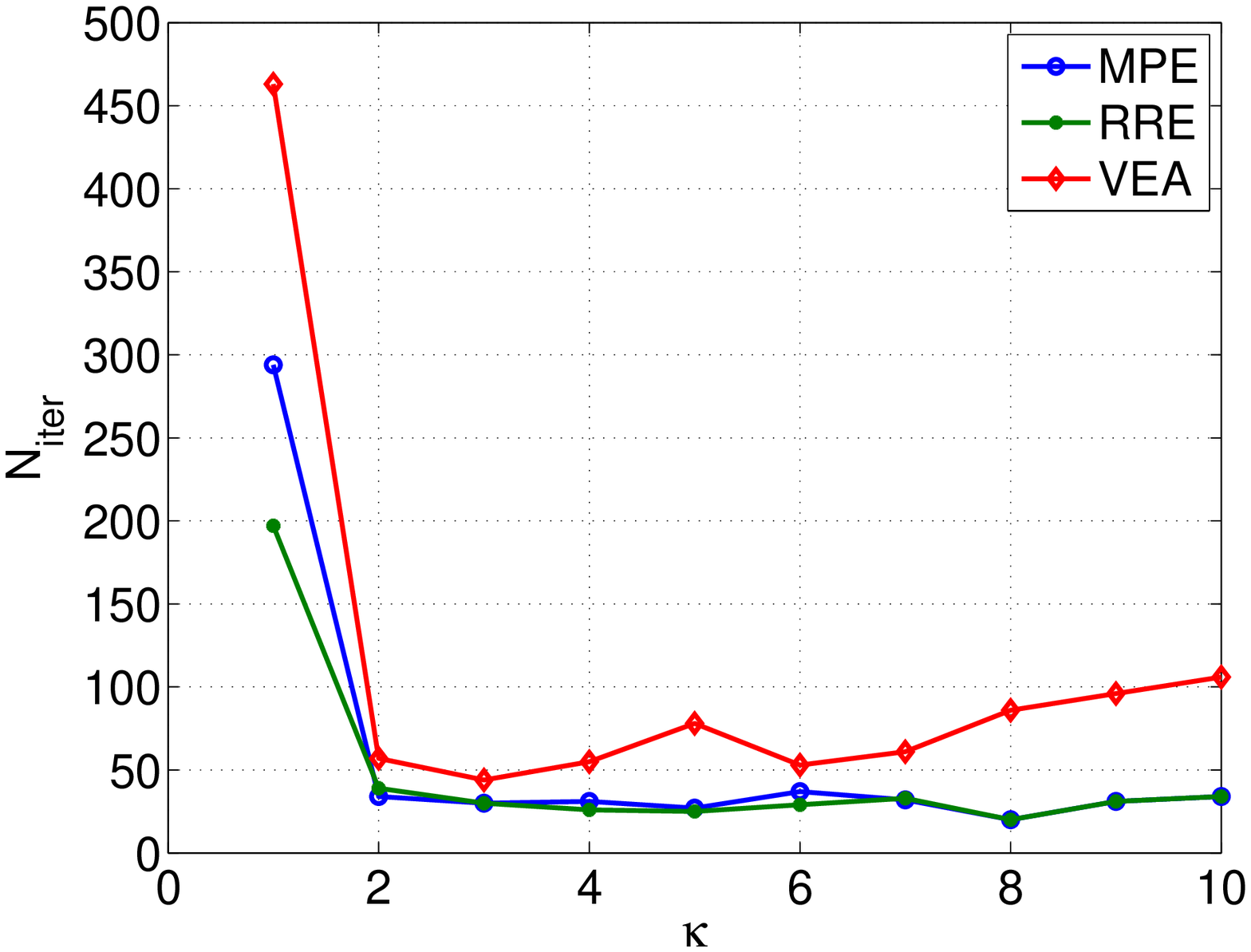} }
\subfigure[]{
\includegraphics[width=6.6cm]{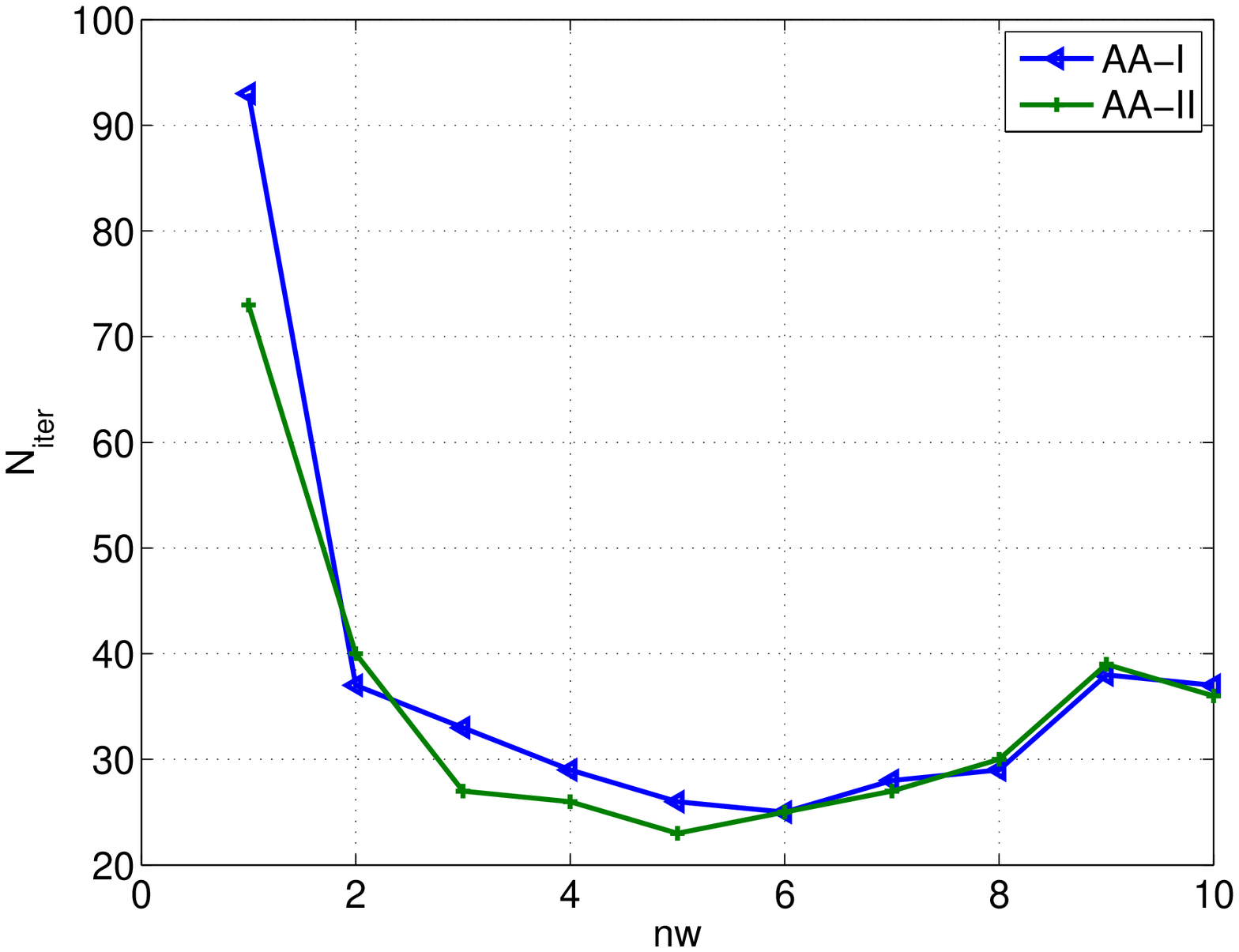} }
\caption{Numerical generation of asymmetric ground state of (\ref{doub_well12}) with $\mu=3.3$. Number of iterations required to reduce the residual error (\ref{fsec32}) below $TOL=10^{-12}$ and as function of the extrapolation width parameters $\kappa$ and $nw$. (a) VEM; (b) AAM.} \label{Fig324}
\end{figure}

\begin{figure}[htbp]
\centering 
\subfigure[]{
\includegraphics[width=6.6cm]{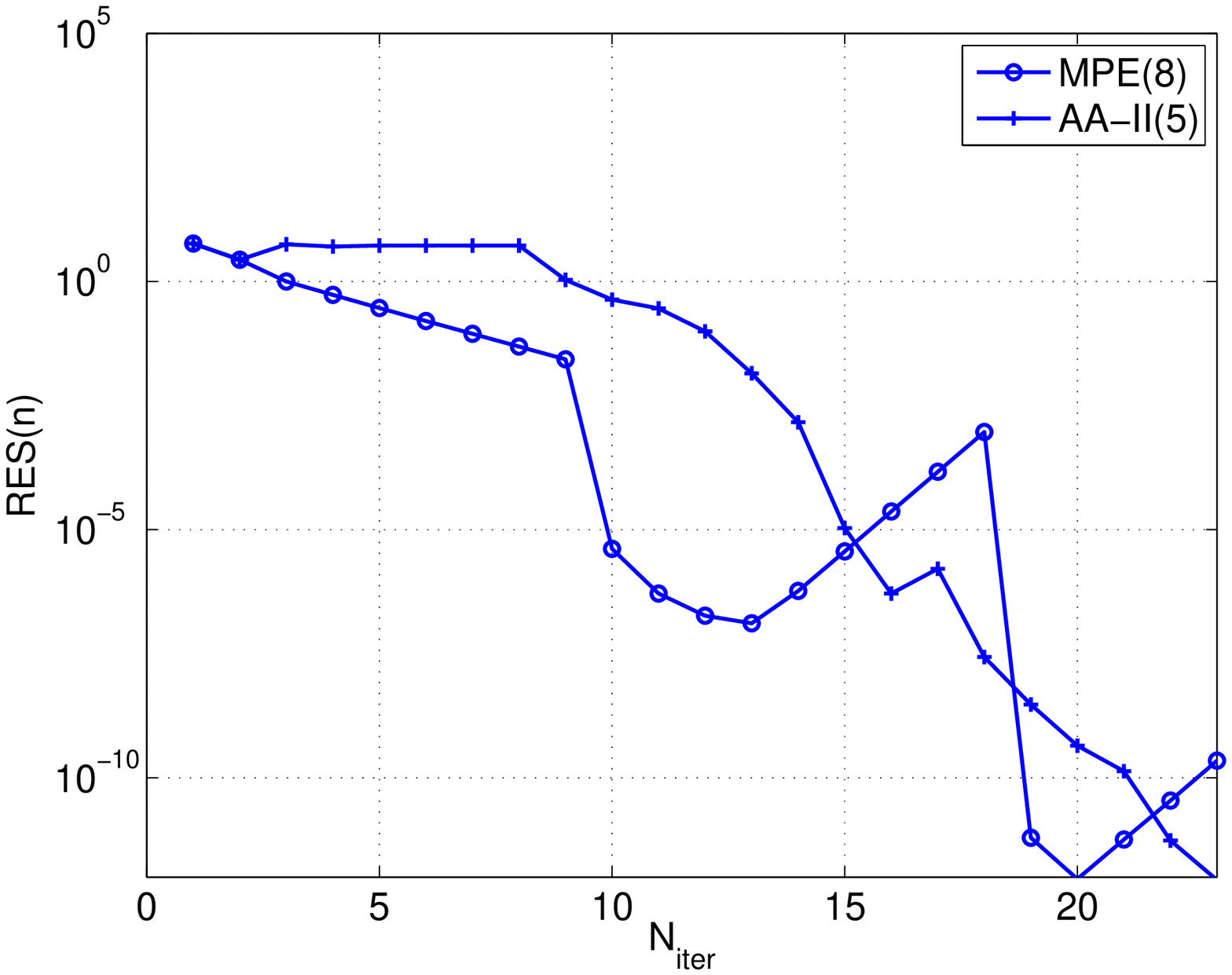} }
\subfigure[]{
\includegraphics[width=6.6cm]{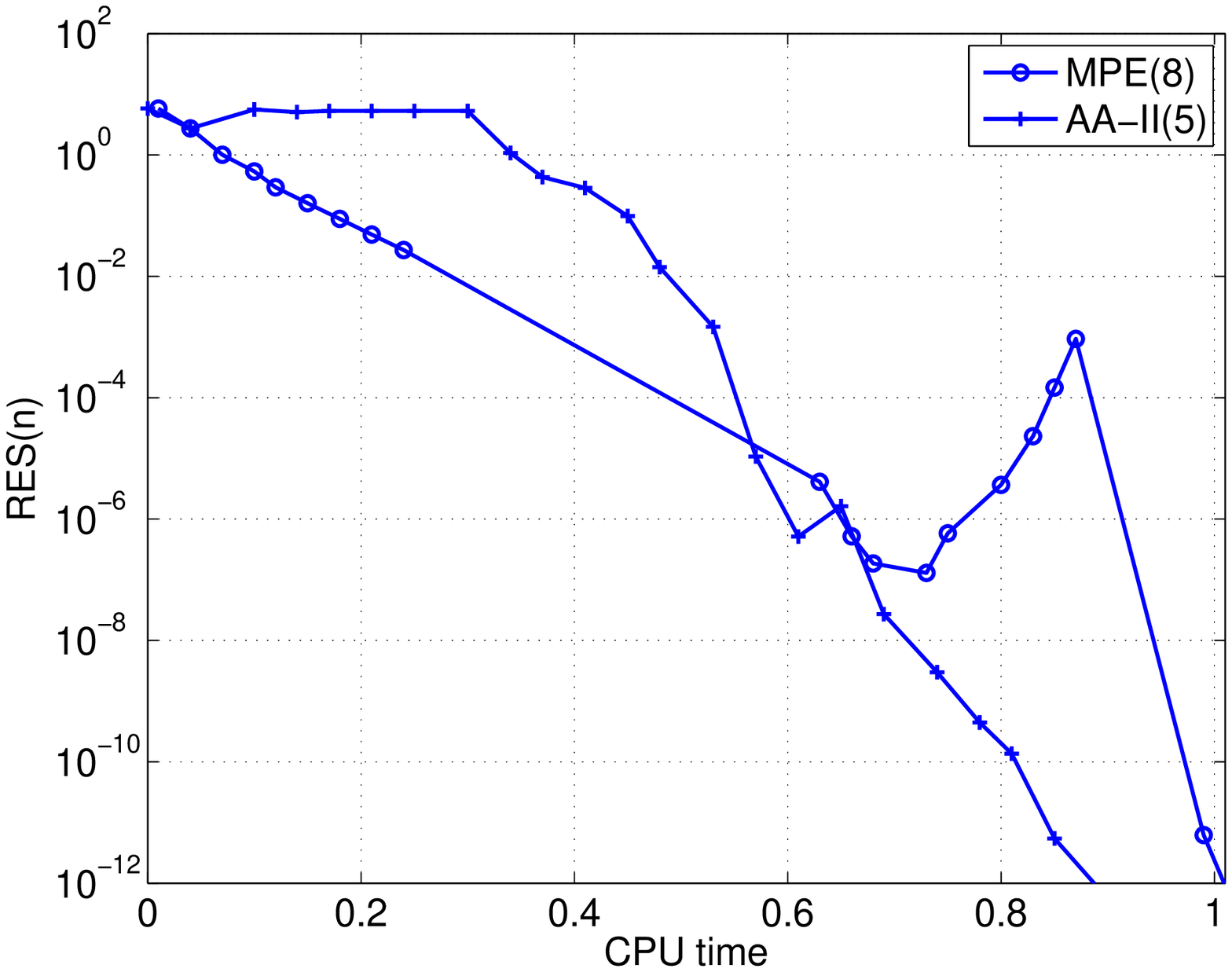} }
\caption{Numerical generation of asymmetric ground state of (\ref{doub_well12}) with $\mu=3.3$. Residual error (\ref{fsec32}) as function of the number of iterations (a) and CPU time in seconds (b) for the \PM method without acceleration (solid line) and accelerated with MPE(8) (circle symbols) and AA-II(5) (plus symbols).} \label{Fig325}
\end{figure}

For other values of $\mu$, some differences are observed. When $\mu=3.3$ the numerical generation of an asymmetric solution of (\ref{doub_well12}) (see Figure \ref{Fig321}(b))  with the \PM method without acceleration is not possible in general. Table \ref{tav1}  (third column) shows the presence of an additional eigenvalue with magnitude above one in the iteration matrix $S$ of the classical fixed point algorithm. As part of the spectrum different from the degree of homogeneity $p=3$, this eigenvalue also appears in the spectrum of the iteration matrix of the \PM method (fourth column) making thus the convergence fail. Here the use of the acceleration techniques corrects this behaviour, leading to convergence. (Both iteration matrices are in fact evaluated at the approximate profile displayed in Figure \ref{Fig321}(b).) In this case (see Figure \ref{Fig324}(a)) MPE and RRE have virtually the same performance while the $\epsilon$-algorithms start reducing their efficiency. (In this case, TEA does not always work in a reliable way and is not competitive against the other VEM.) As far as the AAM are concerned, Figure \ref{Fig324}(b), both improve the performance in a similar, relevant way. They are comparable with VEM in number of iterations (Figure \ref{Fig325}(a)) and behave better when measuring the computational time (Figure \ref{Fig325}(b)).

\begin{figure}[htbp]
\centering 
\subfigure[]{
\includegraphics[width=6.6cm]{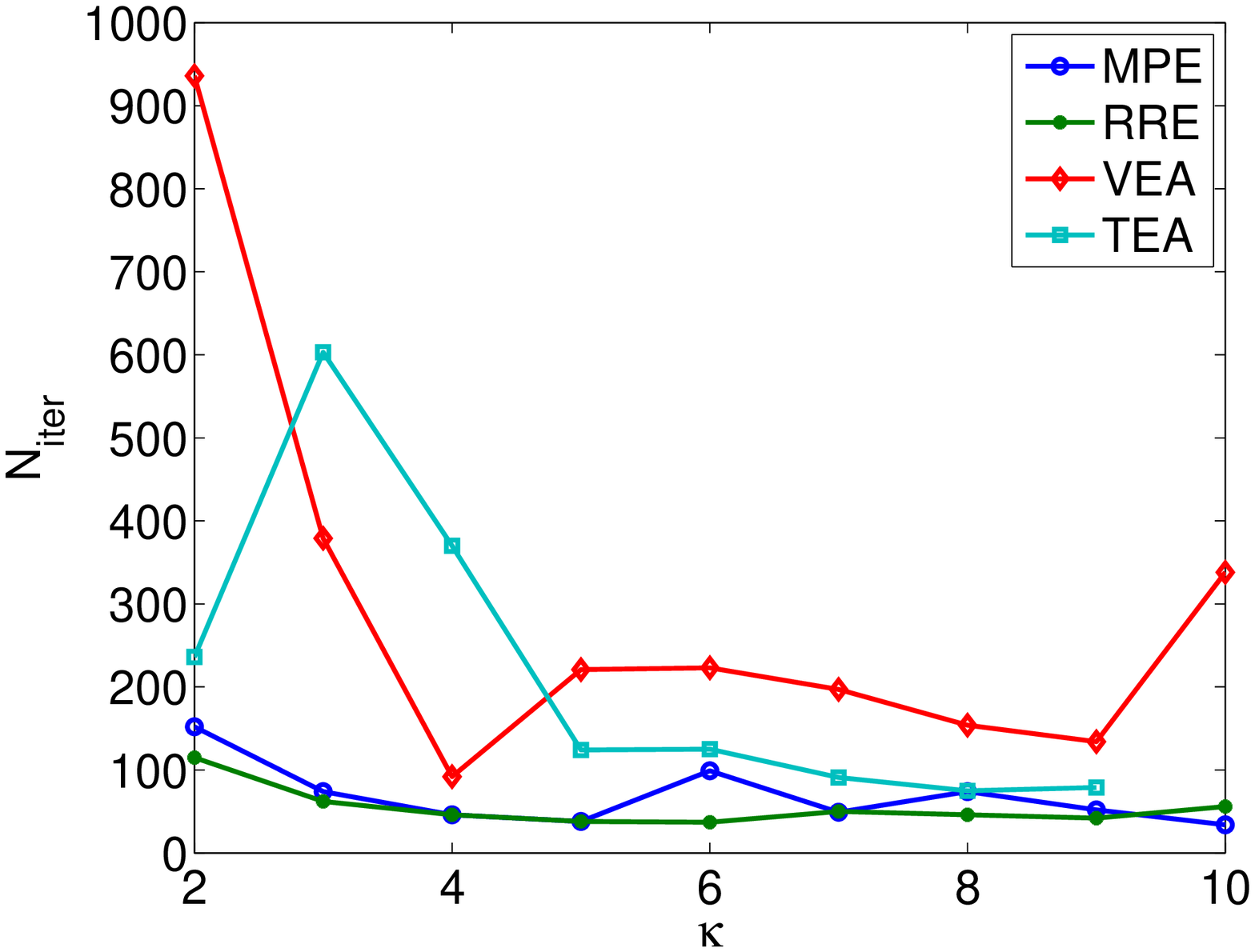} }
\subfigure[]{
\includegraphics[width=6.6cm]{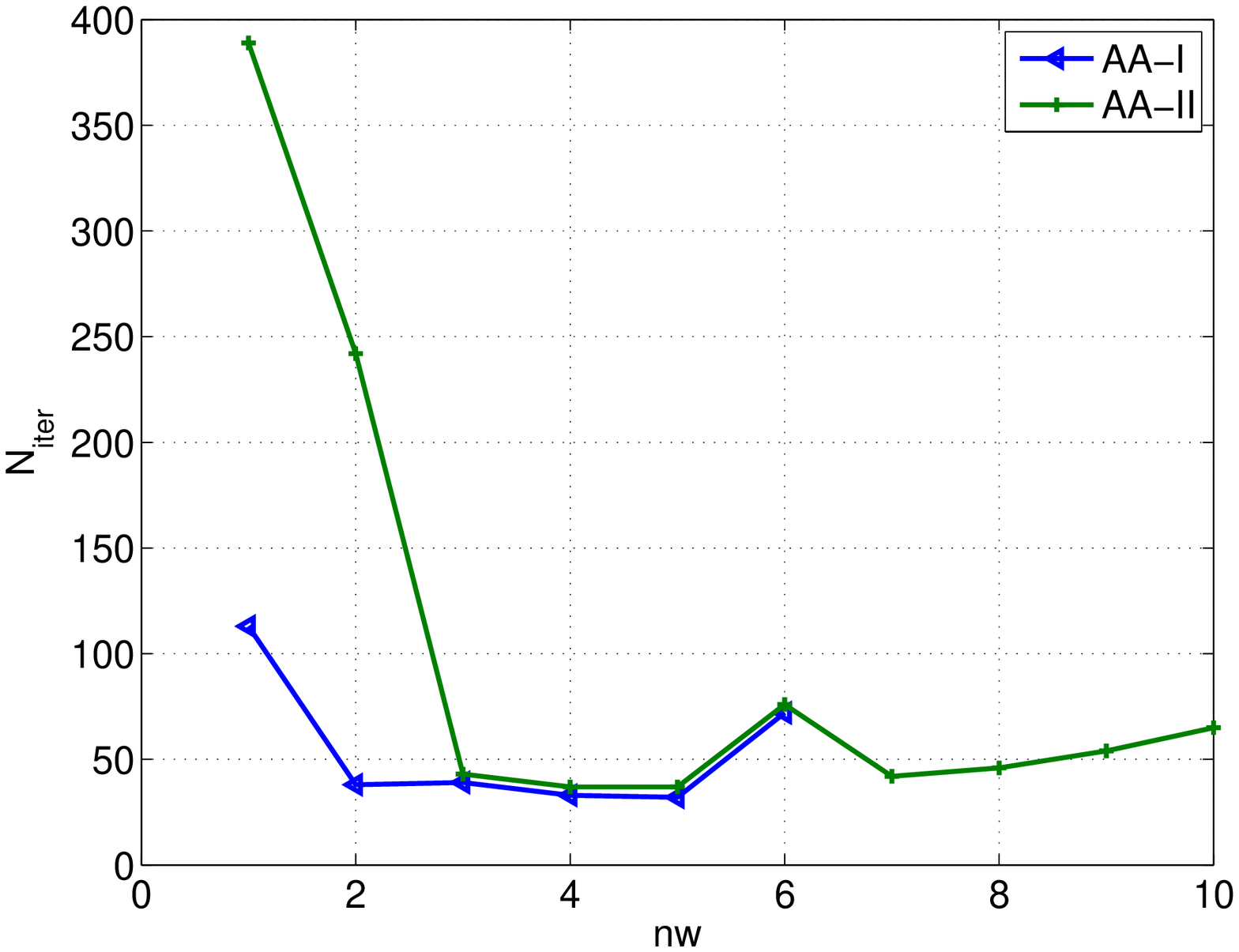} }
\caption{Numerical generation of asymmetric ground state of (\ref{doub_well12}) with $\mu=6.3$. Number of iterations required to reduce the residual error (\ref{fsec32}) below $TOL=10^{-12}$ and as function of the extrapolation width parameters $\kappa$ and $nw$. (a) VEM; (b) AAM.} \label{Fig326}
\end{figure}

\begin{figure}[htbp]
\centering 
\subfigure[]{
\includegraphics[width=6.6cm]{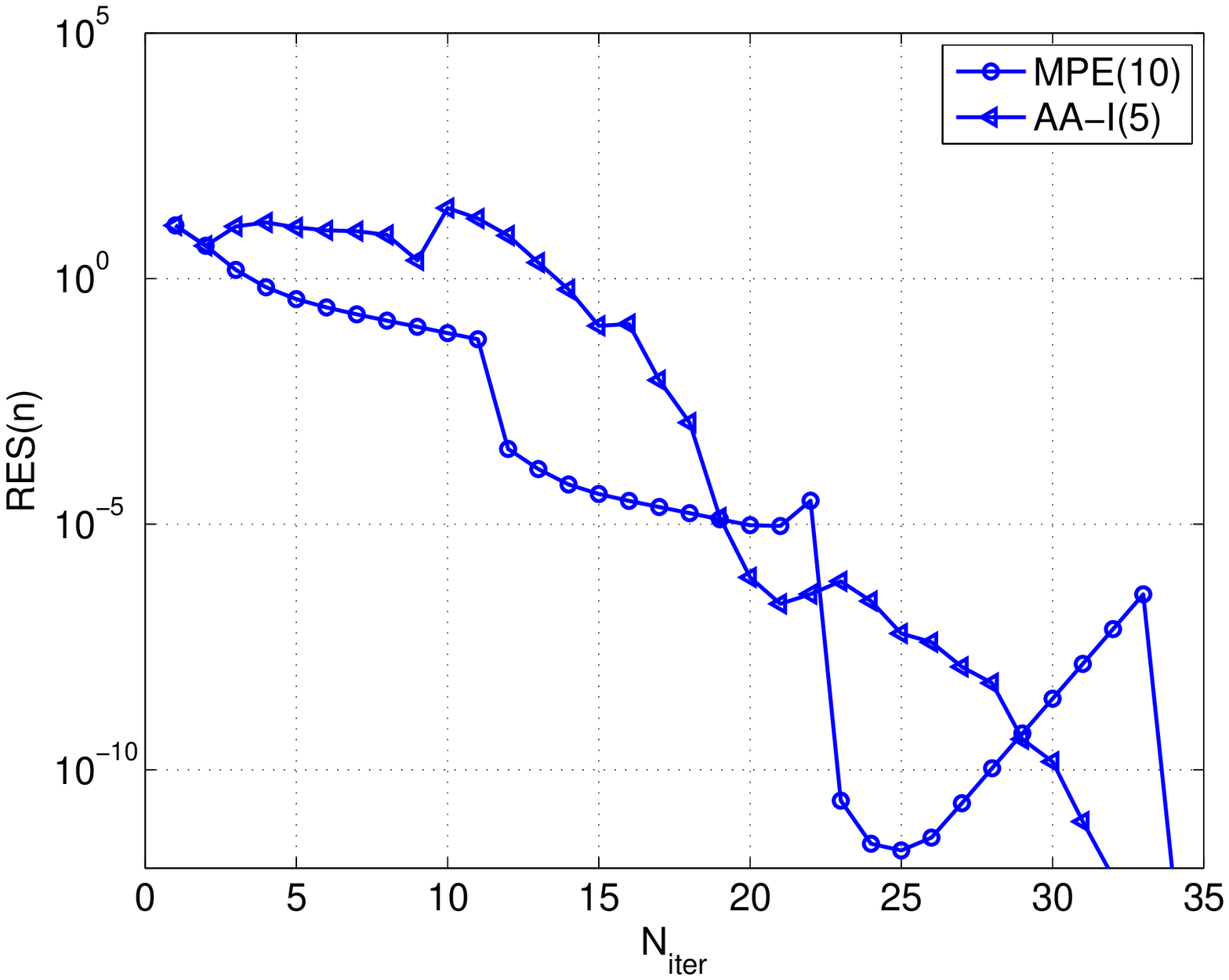} }
\subfigure[]{
\includegraphics[width=6.6cm]{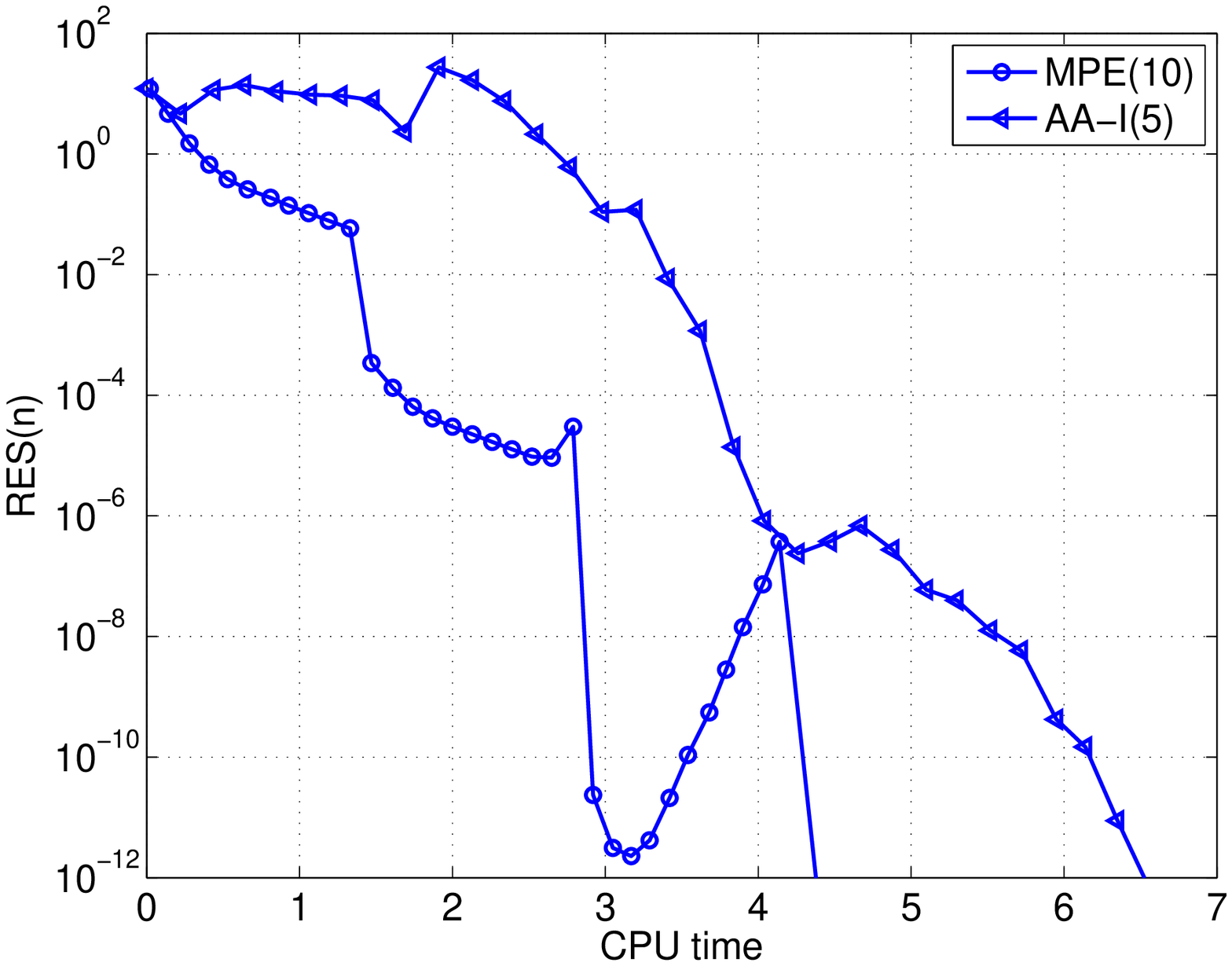} }
\caption{Numerical generation of asymmetric ground state of (\ref{doub_well12}) with $\mu=6.3$. Residual error (\ref{fsec32}) as function of the number of iterations (a) and CPU time in seconds (b) for the \PM method without acceleration (solid line) and accelerated with MPE(10) (circle symbols) and AA-I(5) (triangle symbols).} \label{Fig327}
\end{figure}

\begin{figure}[htbp]
\centering 
\subfigure[]{
\includegraphics[width=6.6cm]{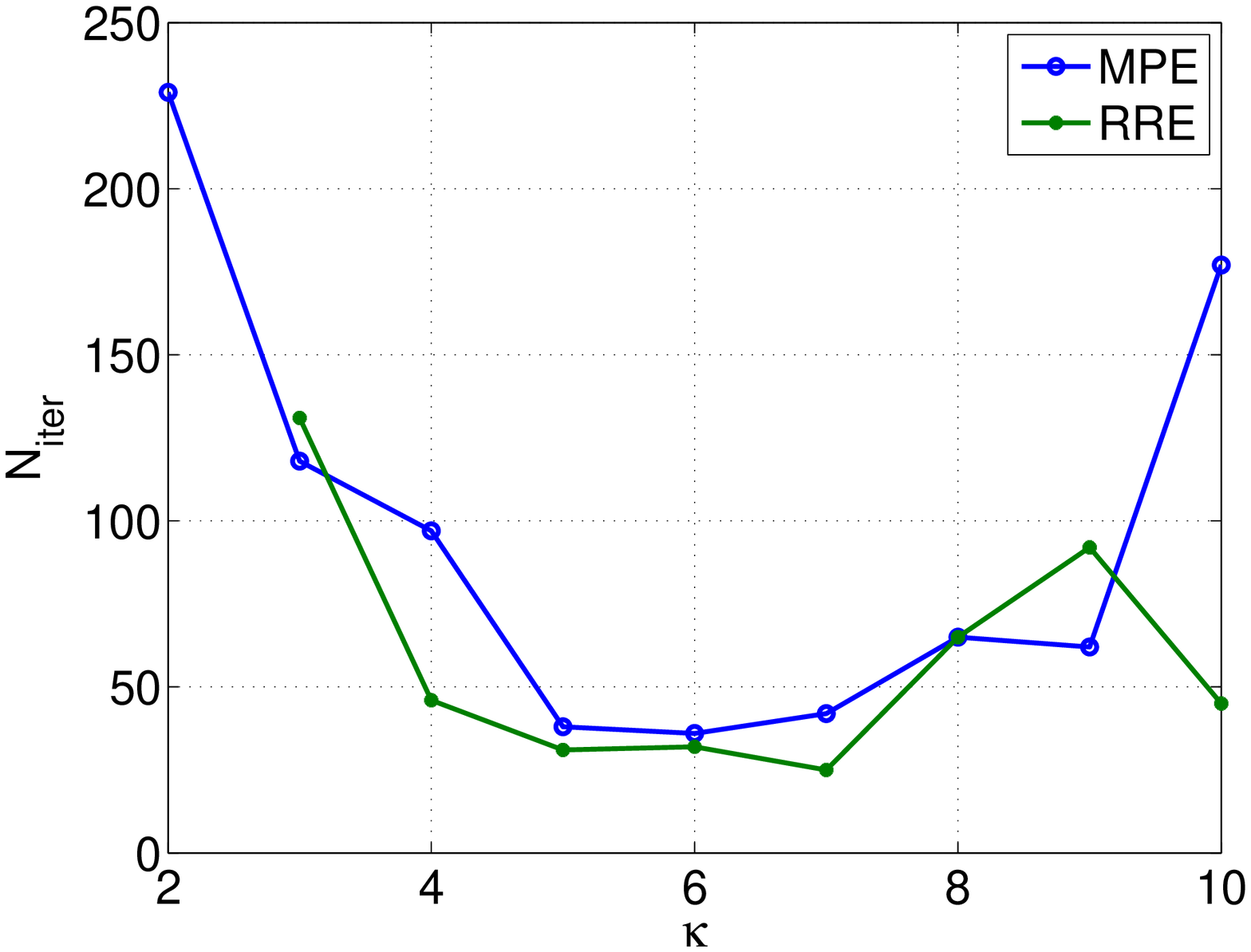} }
\subfigure[]{
\includegraphics[width=6.6cm]{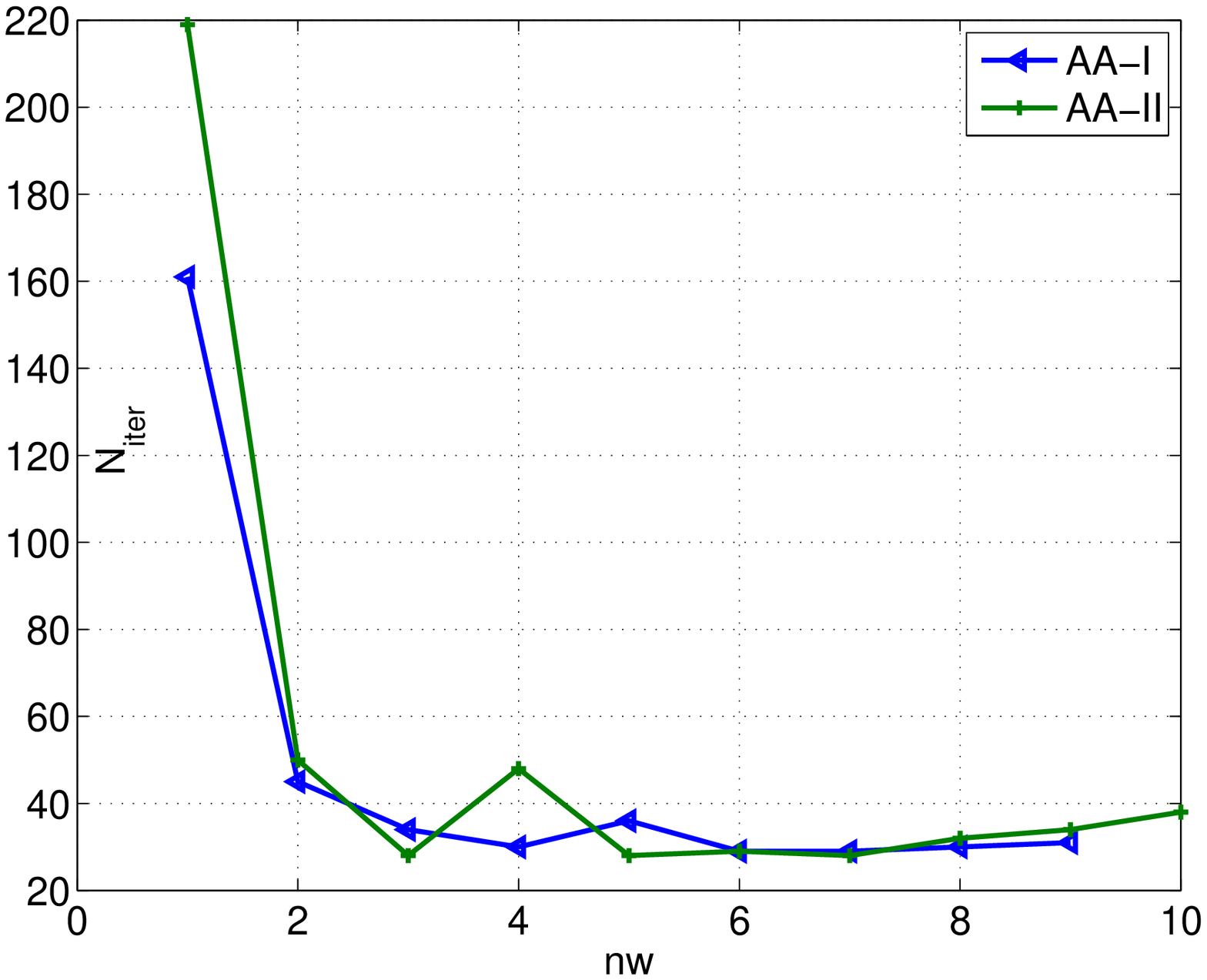} }
\caption{Numerical generation of asymmetric ground state of (\ref{doub_well12}) with $\mu=8.3$. Number of iterations required to reduce the residual error (\ref{fsec32}) below $TOL=10^{-12}$ and as function of the extrapolation width parameters $\kappa$ and $nw$. (a) VEM; (b) AAM.} \label{Fig328}
\end{figure}

\begin{figure}[htbp]
\centering 
\subfigure[]{
\includegraphics[width=6.6cm]{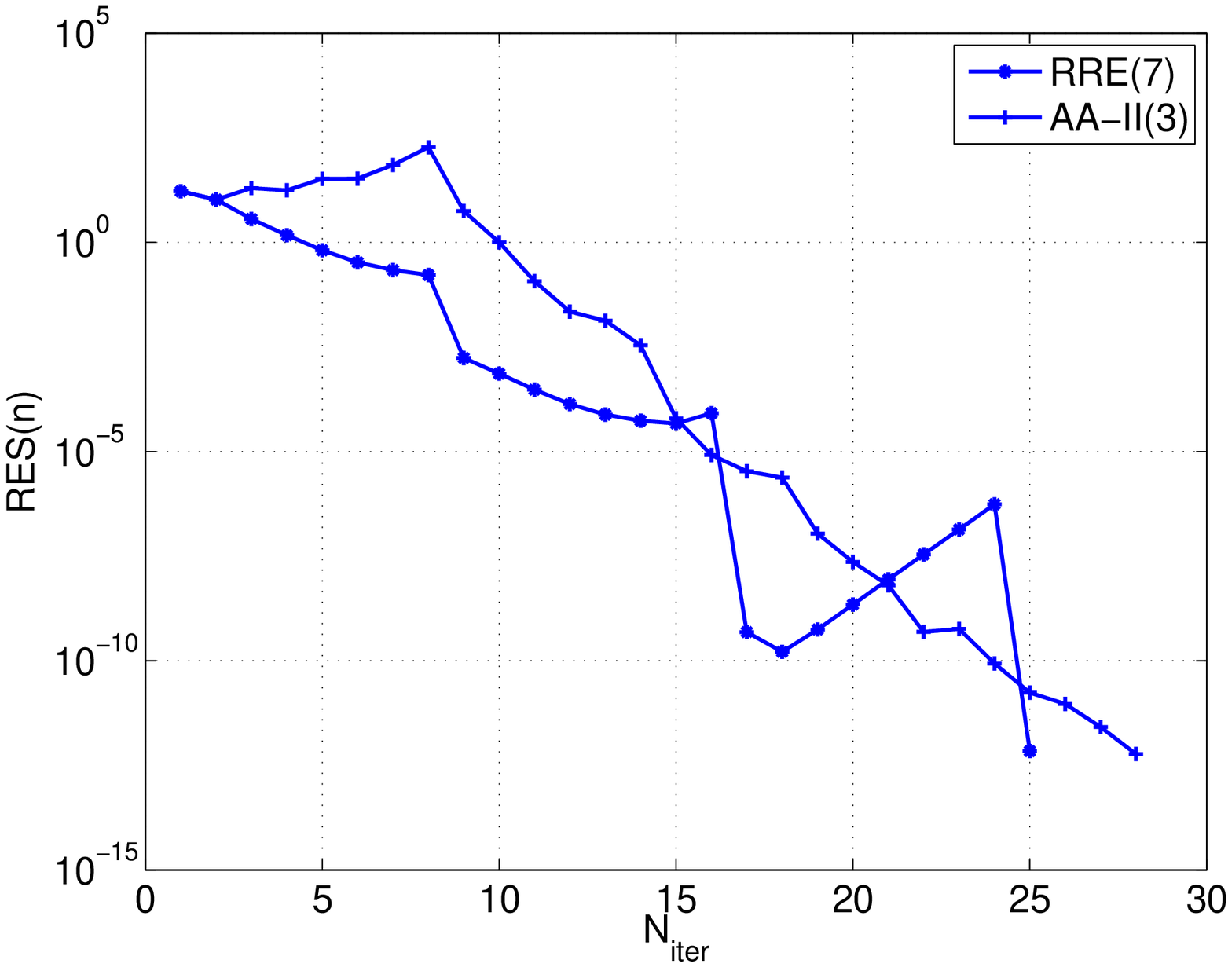} }
\subfigure[]{
\includegraphics[width=6.6cm]{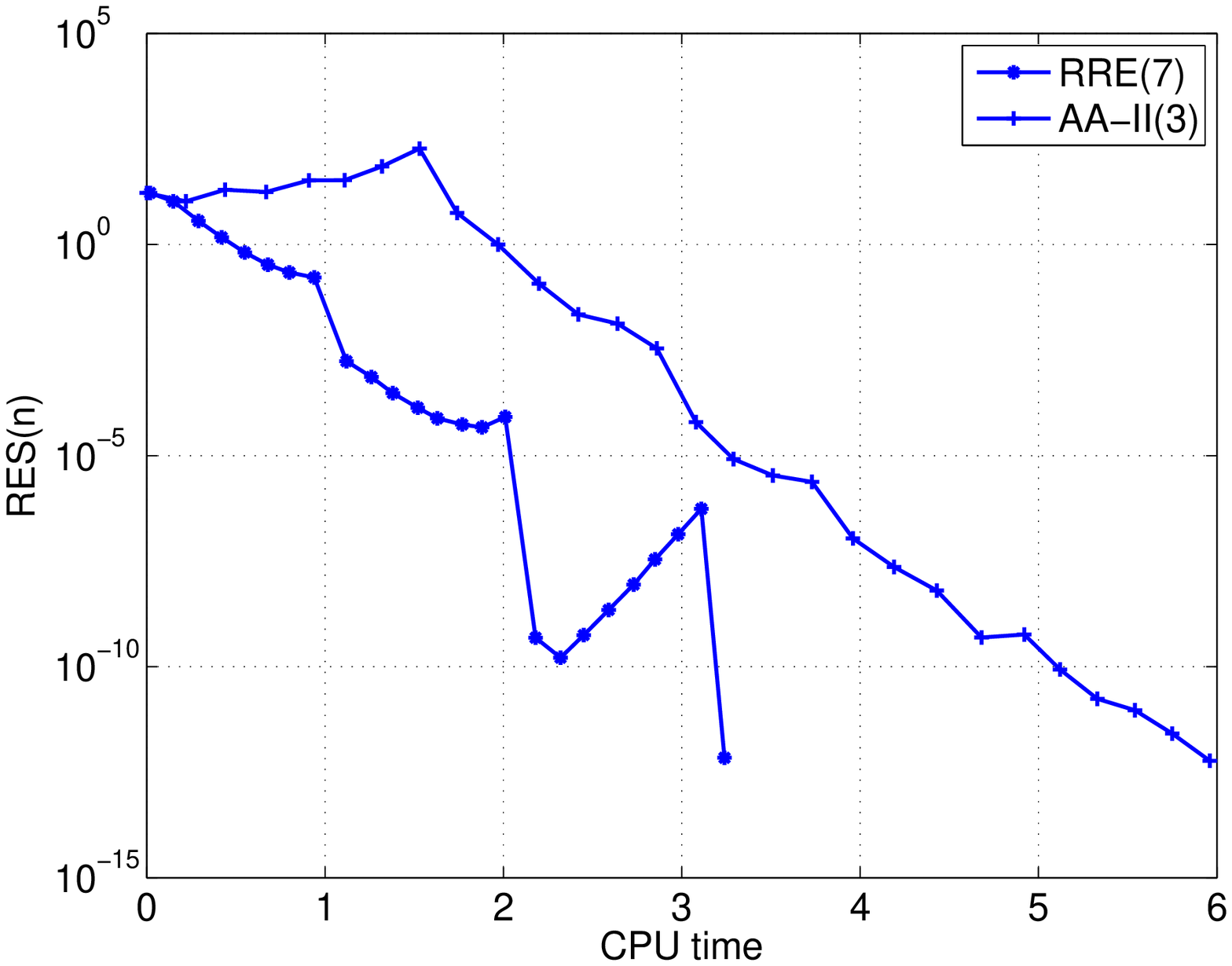} }
\caption{Numerical generation of asymmetric ground state of (\ref{doub_well12}) with $\mu=8.3$. Residual error (\ref{fsec32}) as function of the number of iterations (a) and CPU time in seconds (b) for the \PM method without acceleration (solid line) and accelerated with MPE(7) (circle symbols) and AA-II(3) (plus symbols).} \label{Fig329}
\end{figure}

It may be worth considering the case $\mu=6.3$ because of some relevant points. The first one is the generation of the asymmetric profile, Figure \ref{Fig321}(c), which in general is not possible with the \PM method without acceleration. The situation is similar to that of the previous case $\mu=3.3$ and it is shown in Table \ref{tav2} (first and second columns). In this case, the best results of the acceleration are given by MPE and AA-I (Figure \ref{Fig326}). The loss of performance of the $\epsilon$-algorithms and the improvement of AAM, observed in the previous experiments, are confirmed here and in the experiments for $\mu=8.3$ (Figures \ref{Fig328} and \ref{Fig329}). The comparison between MPE and AA-I, see Figures \ref{Fig327}(a), (b), reveals, ikn the authors' opinion,  a similar performance.

\begin{table}
\begin{center}
\begin{tabular}{|c|c||c|c|}
\hline\hline
$\mu=6.3$&$\mu=6.3$&$\mu=8.3$&$\mu=8.3$\\\hline\hline
eigs$(S)$& eigs$(F^{\prime}(u^{*}))$&eigs$(S)$& eigs$(F^{\prime}(u^{*}))$\\\hline
5.095370E+00&5.096207E+00&3.962824E+00&3.962824E+00\\
3.000000E+00&9.672929E-01&2.999999E+00&9.807797E-01\\
9.672929E-01&7.506018E-01&9.807797E-01&8.081404E-01\\
7.506018E-01&4.078905E-01&8.081404E-01&4.459845E-01\\
4.078905E-01&3.472429E-01&4.459845E-01&4.030040E-01\\
3.472429E-01&2.032986E-01&4.030040E-01&1.929797E-01\\
\hline\hline
\end{tabular}
\end{center}
\caption{Numerical generation of asymmetric profile of (\ref{doub_well12})
with $\mu=6.3$ and $\mu=8.3$. Six largest magnitude eigenvalues of the
approximated iteration matrix of the classical fixed-point method $S=L^{-1}N^{\prime}(U_{f})$ (left)  and of the \PM method (\ref{mm2}), (\ref{mm3c}), evaluated
at the last computed iterate $U_{f}$ obtained with MPE($10$).}\label{tav2}
\end{table}

\begin{figure}[htbp]
\centering 
\includegraphics[width=8cm]{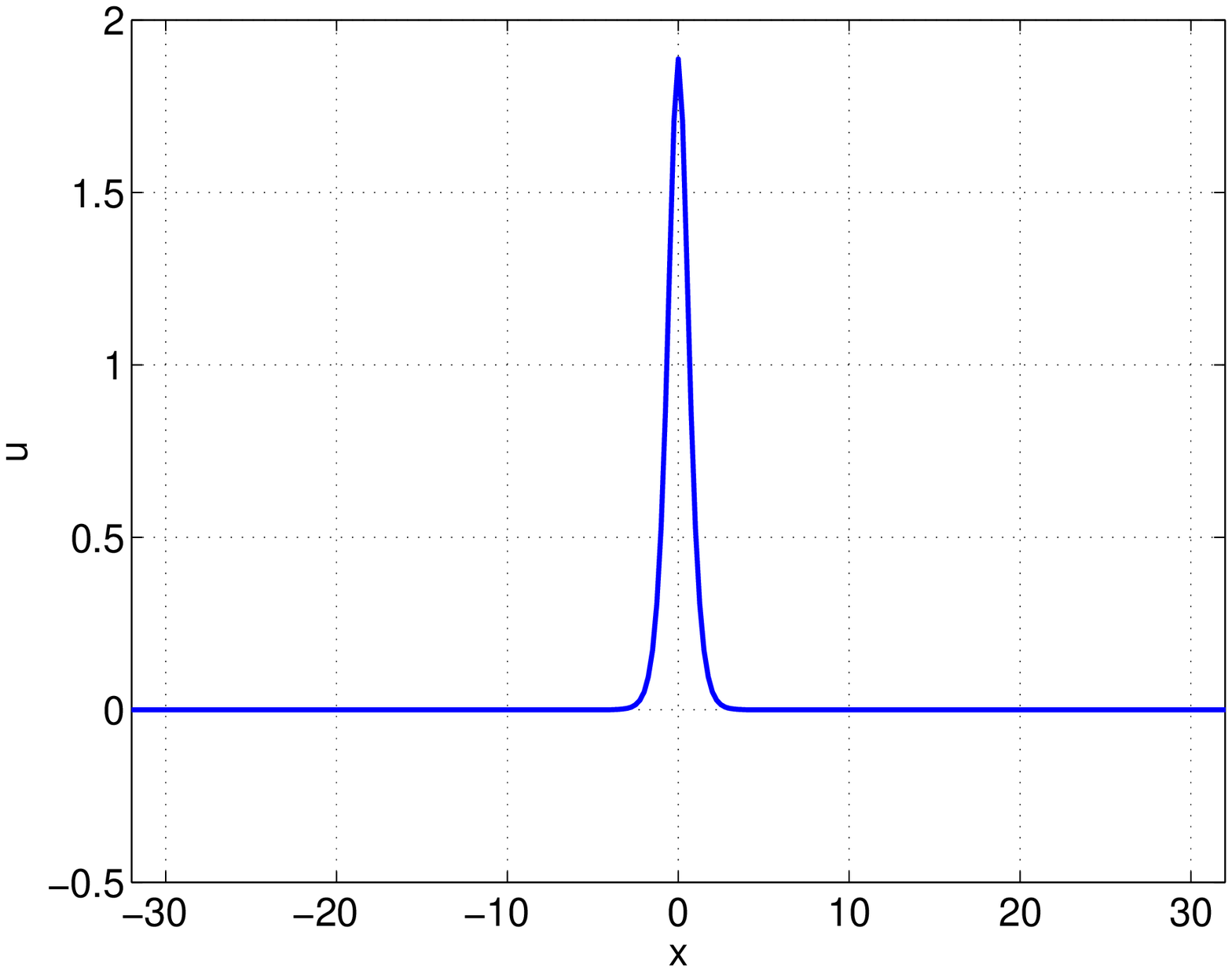}
\caption{Numerical generation of symmetric ground state of (\ref{doub_well12}) with $\mu=6.3$. Approximate profile with \PM method without acceleration.} \label{Fig3210}
\end{figure}

The second question with regard to the case $\mu=6.3$ concerns the behaviour of the \PM method without acceleration. In this case, the method is convergent, but to a symmetric localized wave, see Figure \ref{Fig3210}. 
\begin{table}
\begin{center}
\begin{tabular}{|c|c|}
\hline\hline
eigs$(S)$& eigs$(F^{\prime}(u^{*}))$\\\hline
3.000000E+00&6.098684E-01\\
6.098684E-01&2.696853E-01\\
2.696853E-01&1.518421E-01\\
1.518421E-01&1.039553E-01\\
1.039553E-01&6.046185E-02\\
6.046185E-02&5.492737E-02\\
\hline\hline
\end{tabular}
\end{center}
\caption{Numerical generation of symmetric profile of (\ref{doub_well12})
with $\mu=6.3$. Six largest magnitude eigenvalues of the
approximated iteration matrix of the classical fixed-point method $S=L^{-1}N^{\prime}(U_{f})$ (left)  and of the \PM method (\ref{mm2}), (\ref{mm3c}), evaluated
at the last computed iterate $U_{f}$ obtained with \PM method (\ref{mm2}), (\ref{mm3c}).}\label{tav3}
\end{table}
This can be explained by the first two columns of Table \ref{tav2}  and by Table \ref{tav3}. Note that, as mentioned before, for the asymmetric solution, the \PM method cannot be convergent. However, according to the information provided by Table \ref{tav3}, this is locally convergent to the symmetric solution. (In this case, the spectral radius of $F^{\prime}(u^{*})$ is below one.)
This profile can be indeed approximated by using acceleration techniques (and with the corresponding computational saving) but starting from a different initial iteration.

Finally, the case $\mu=8.3$ is also analyzed, see Figure \ref{Fig321}(d). The main reason we find to emphasize this case is to confirm the conclusions obtained from the experiments with the previous values of $\mu$:
\begin{itemize}
\item Among the VEM, the polynomial methods give a better performance, while the $\epsilon$-algorithms become less efficient as $\mu$ increases. As observed in Figure \ref{Fig321}, the larger $\mu$ the larger and narrower the asymmetric profile is. The computation becomes harder as is noticed by comparing the iterations required by the methods in Figures \ref{Fig322}, \ref{Fig324}, \ref{Fig326} and \ref{Fig328}. One can also note the increment of the magnitude of the eigenvalues of the corresponding iteration matrices of the \PM method in Tables \ref{tav1} and \ref{tav2}. Therefore, under more demanding conditions, the polynomial methods give a better answer than the $\epsilon$-algorithms.
\item Contrary to the $\epsilon$-algorithms, whose performance gets worse as $\mu$ increases, the AAM improve their behaviour up to being comparable with polynomial methods (cf. the periodic traveling wave generation in Section \ref{se3}). Furthermore, this is obtained with small values of the parameter $nw$, thus avoiding ill-conditioned problems.
\end{itemize}

\subsection{Example 3. Solitary wave solutions of the Benjamin equation}
An additional application of acceleration techniques concerns the oscillatory character of the wave to be numerically generated. This property has shown to make influence on the performance of the iteration with eventual loss of convergence in some cases, \cite{DougalisDM2012,DougalisDM2015}. Presented here is the use of acceleration as an alternative to overcome this difficulty and this will be illustrated with the numerical generation of solitary waves in one- and two-dimensional versions of the Benjamin equation.
\subsubsection{One-dimensional Benjamin equation}
A first example of the situation described above is given by the solitary wave solutions of the Benjamin equation, \cite{ben0}
\begin{equation}\label{E11}
u_t+\alpha u_x+\beta u u_x-\gamma \mathcal{H}u_{xx}-\delta u_{xxx}=0,
\end{equation}
where $u=u(x,t), x\in \mathbb{R}, t\geq 0, \alpha,\beta,\gamma,\delta$ are positive constants, and $\mathcal{H}$ denotes the Hilbert transform defined on the real line as
\begin{eqnarray}
\mathcal{H}f(x):=\frac{1}{\pi}p.v.\int_{-\infty}^{\infty}\frac{f(y)}{x-y}\,dy,\label{hilb}
\end{eqnarray} or through its Fourier transform as $$\widehat{\mathcal{H}f}(k)=-{\rm i}{\rm sign}(k)\widehat{f}(k), \quad k\in\mathbb{R}.$$
Equation (\ref{E11}) is a model for the propagation of internal waves along the interface of a two-layer fluid system and where gravity and surface tension effects are not negligible. It includes the limiting cases of negligible surface tension ($\delta=0$ or Benjamin-Ono equation) and a limit of a model with very thin upper fluid ($\gamma=0$ or KdV equation). Solitary-wave solutions of (\ref{E11}) with speed $c_{s}>0$ are determined by profiles  $u(x,t)=\varphi(x-c_{s}t), c_{s}>0$, such that $\varphi$ and its derivatives tend to zero  as $X=x-c_{s}t$ approaches $\pm\infty$ and satisfying
\begin{eqnarray}
(\alpha-c_{s})\varphi
+\frac{\beta}{2}\varphi^{2}-\gamma \mathcal{H}\varphi^{\prime}-\delta
\varphi^{\prime\prime}=0,\label{E15}
\end{eqnarray}
where ${}^{\prime}=d/dX$. Albert et al., \cite{AlbertBR1999} established a complete theory of existence and orbital stability of solitary waves of (\ref{E11}) for small $\gamma$, while Benjamin, \cite{ben1}, derived the oscillating behaviour of the waves, with the number of oscillations increasing as $\gamma$ approaches $\gamma^{\ast}=2\sqrt{\delta(\alpha-c_{s})}$, along with the asymptotic decay, as $|X|\rightarrow\infty$, like $1/X^{2}$.

Except in the limiting cases, solitary wave solutions are not analytically known. A standard way to generate solitary wave profiles numerically consists of considering (\ref{E15}) in the Fourier space
\begin{eqnarray}\label{Be3}
(-c_{s}+\alpha-\gamma |k|+\delta k^{2})\widehat{\varphi}
+\frac{\beta}{2}\widehat{\varphi^{2}}=0,\quad k\in\mathbb{R},
\end{eqnarray}
(where $\widehat{\varphi}(k)$ is the Fourier transform of $\varphi$), discretizing (\ref{Be3}) with periodic boundary conditions on a sufficiently long interval $(-l,l)$ and the use of discrete Fourier transform (DFT)
\begin{eqnarray}\label{E41}
(-c_{s}+\alpha-\gamma |k|+\delta k^{2})\widehat{\varphi^{N}}_{k}
+\frac{\beta}{2}\left(\widehat{\varphi^{N}\ast \varphi^{N}}\right)_{k}=0,
\end{eqnarray}
for $k=-\frac{N}{2},\ldots,\frac{N}{2}-1$, where $\varphi^{N}$ is a trigonometric polynomial of degree $N$ which approximates $\varphi$ and $\widehat{\varphi^{N}}_{k}$ denotes its $k^{\rm th}$ Fourier coefficient. Then (\ref{E41}) is numerically solved by incremental continuation with respect to $\gamma$ from $\gamma=0$ (which corresponds to KdV equation and for which solitary wave profiles are analytically known) and a nonlinear iteratively solver for each value of the homotopic path with respect to $\gamma$. For a more detailed description of the incremental continuation method and the performance of several nonlinear iterative solvers see \cite{AlbertBR1999,DougalisDM2012}. The experiments performed there reveal that the oscillatory behaviour of the wave increases the difficulty of its computation, even using numerical continuation. Our aim here is giving a computational alternative, based on the use of acceleration techniques.

\begin{figure}[htbp]
\centering \subfigure[]{
\includegraphics[width=6.6cm]{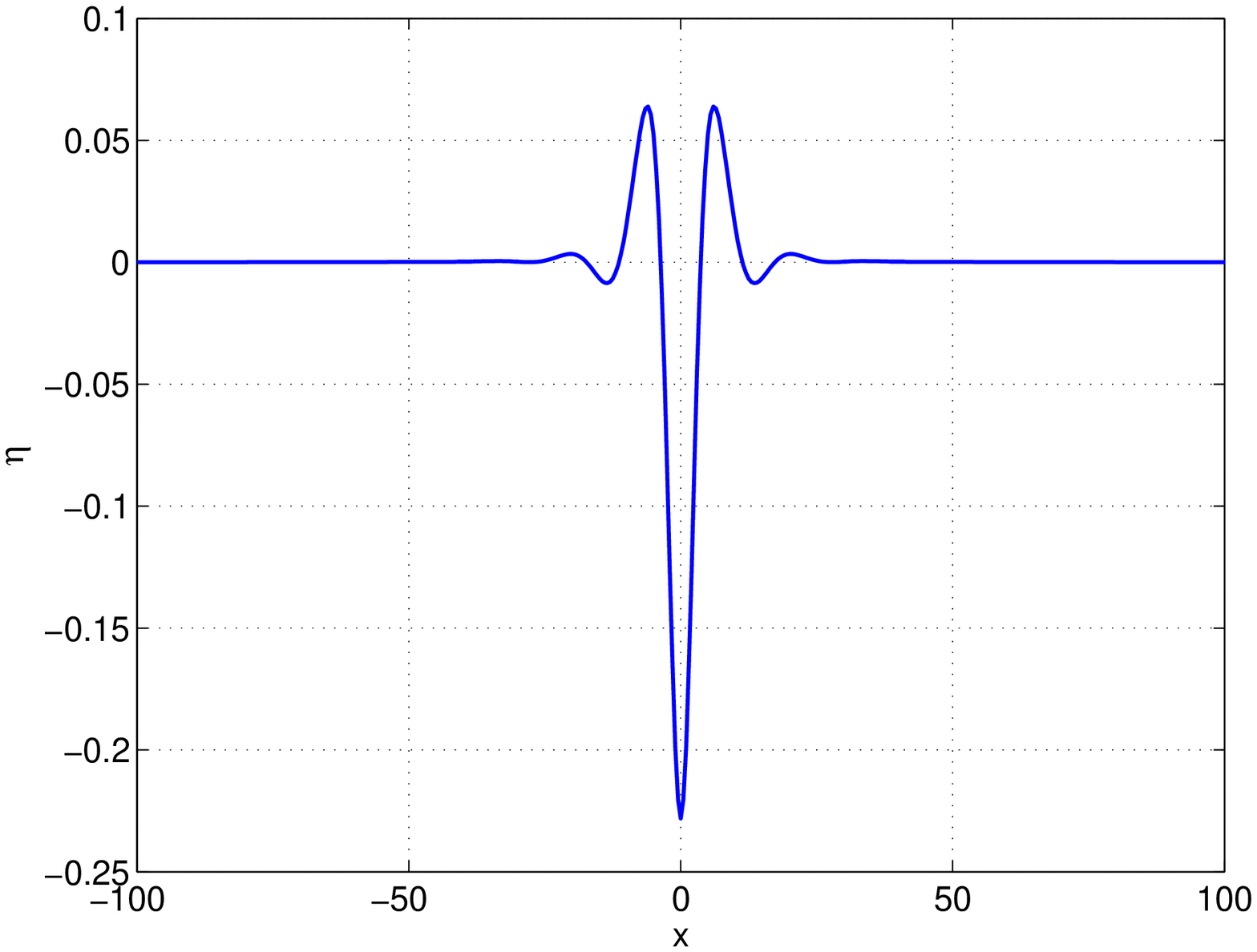} }
\subfigure[]{
\includegraphics[width=6.6cm]{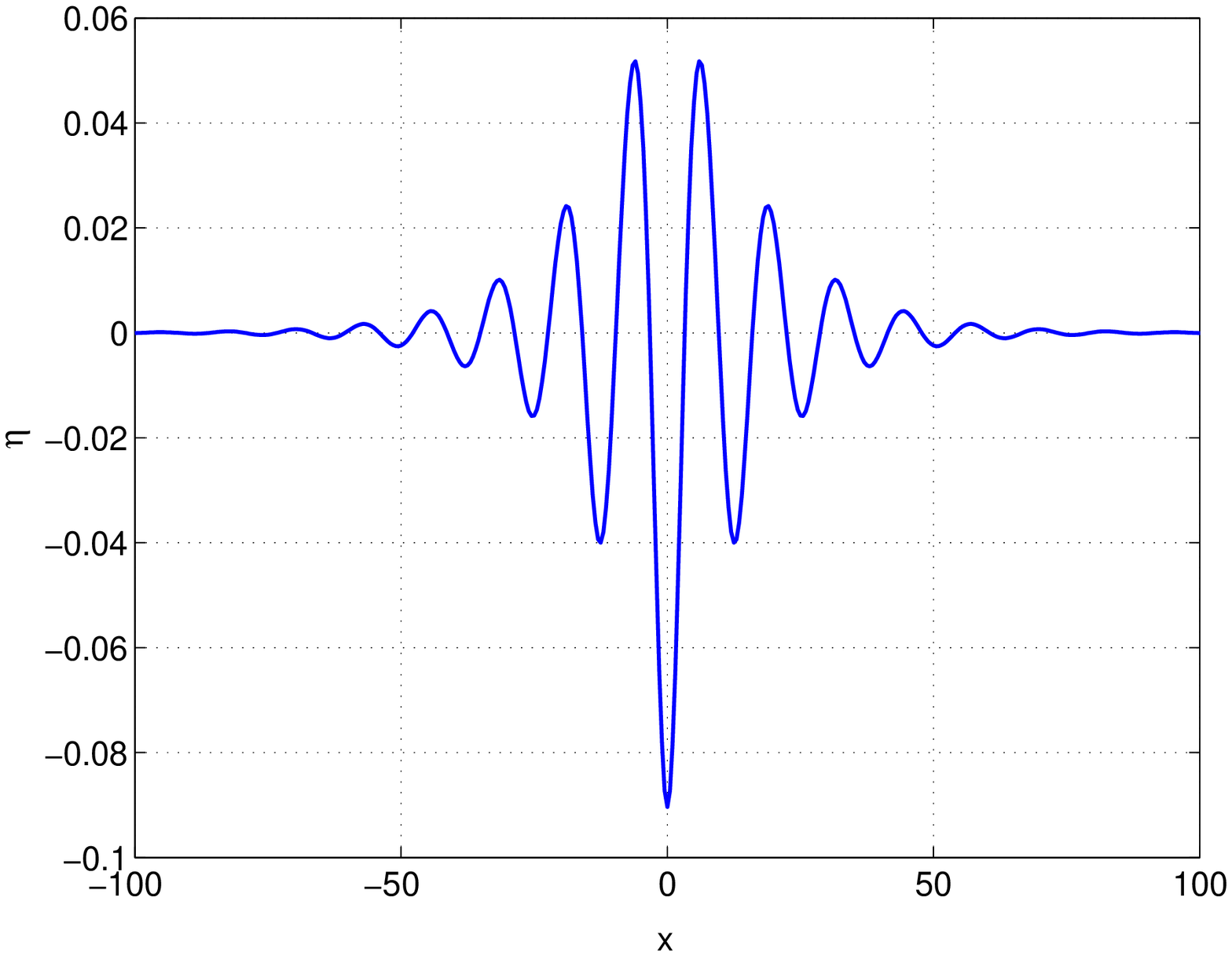} }
\subfigure[]{
\includegraphics[width=6.6cm]{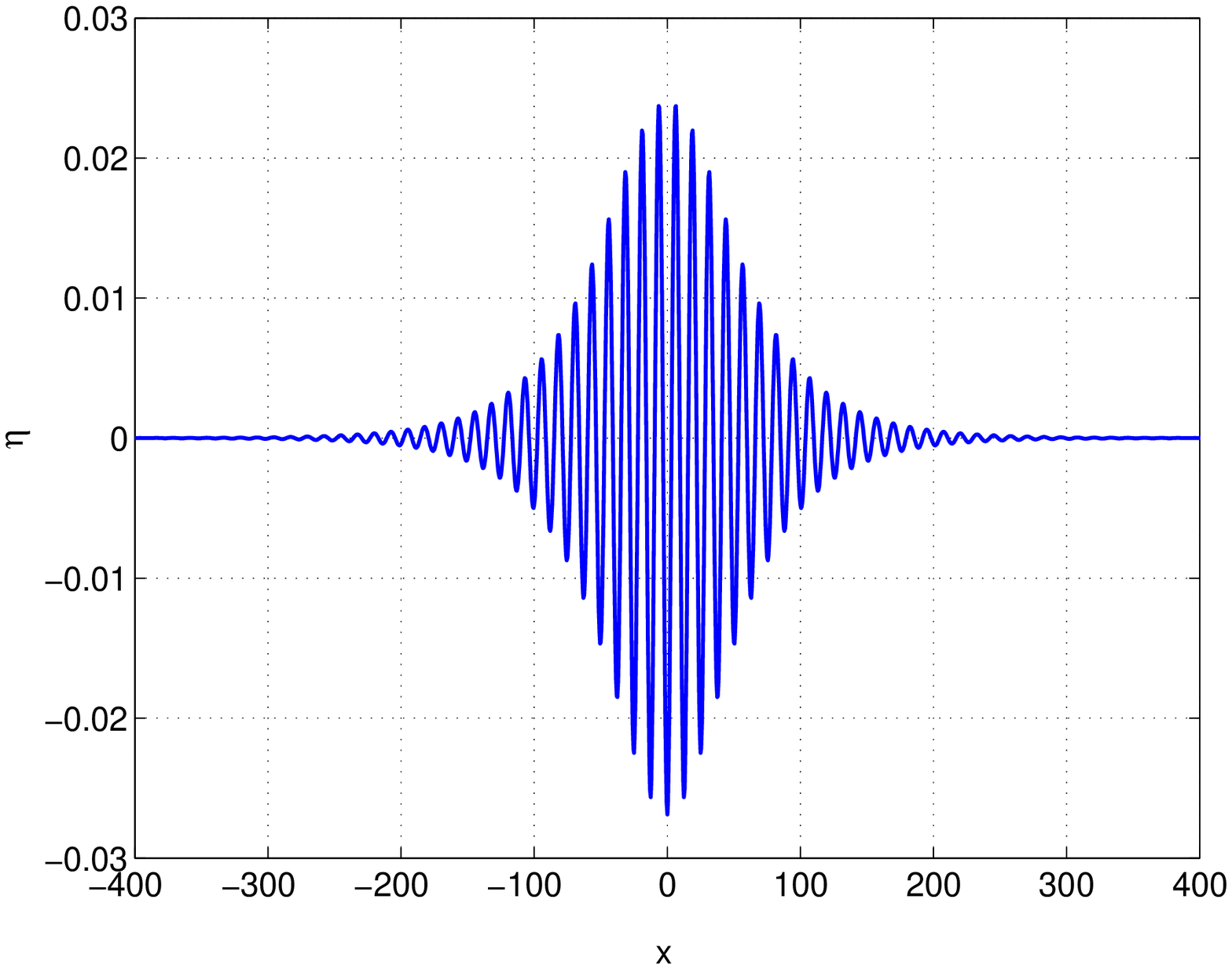} }
\subfigure[]{
\includegraphics[width=6.6cm]{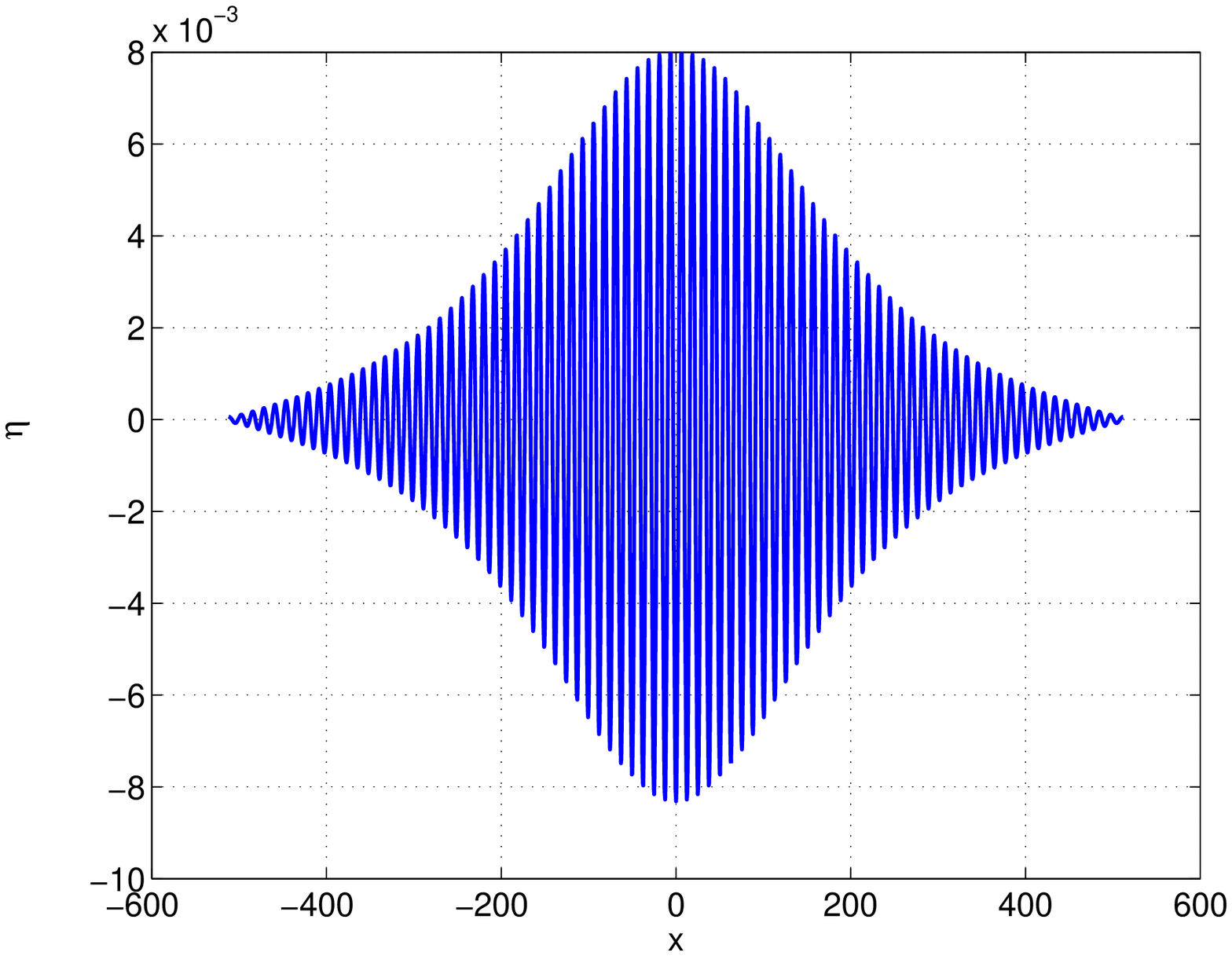} }
\caption{Approximate solitary-wave solutions of (\ref{E15}) given by MPE method with $\kappa=7$ and $c_{s}=0.75$, $\alpha=\beta=\delta=1$. (a)-(d) correspond, respectively, to $\gamma=0.9,0.99,0.999,0.9999$.} \label{Ben1}
\end{figure}

\begin{figure}[htbp]
\centering \subfigure[]{
\includegraphics[width=6.6cm]{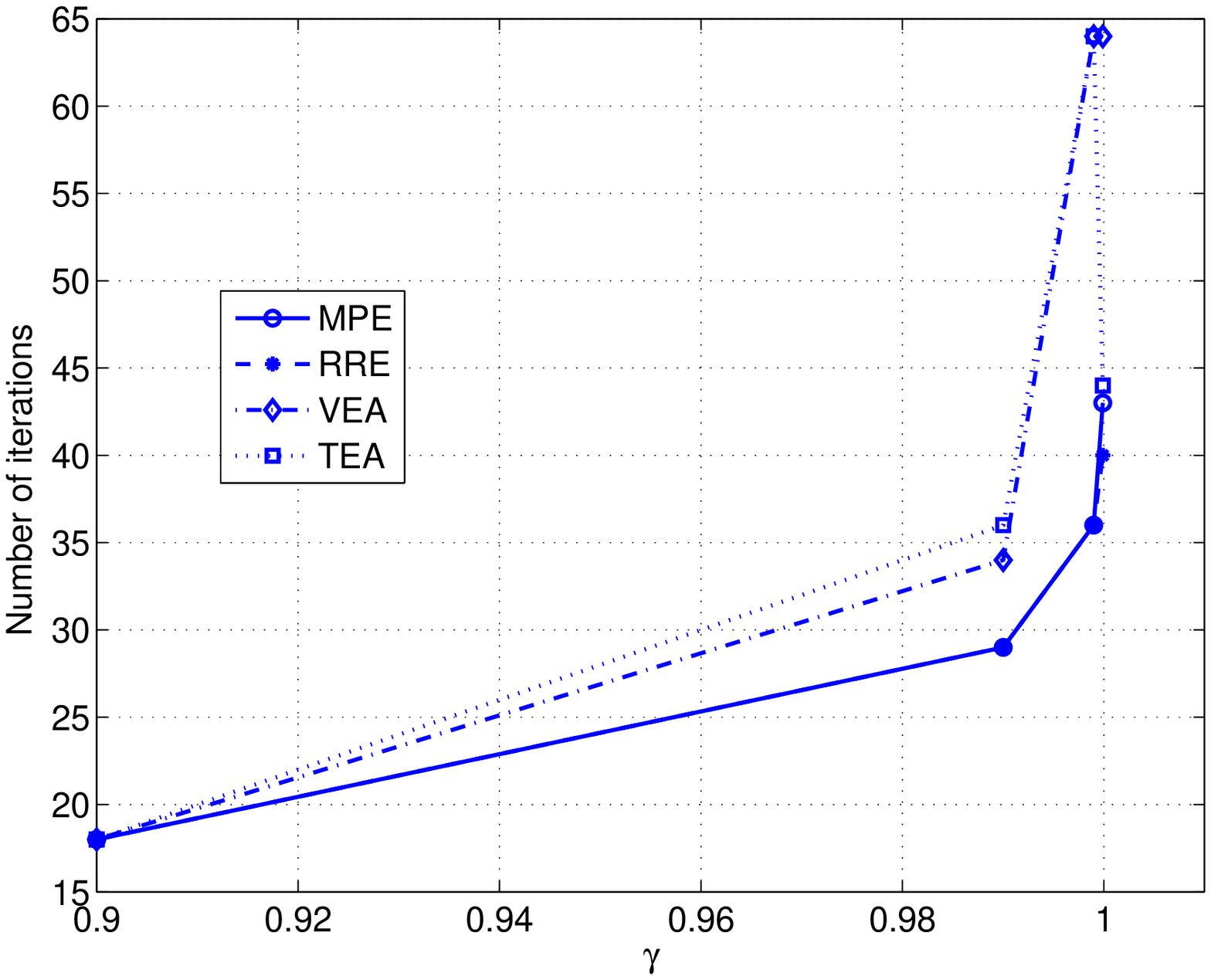} }
\subfigure[]{
\includegraphics[width=6.6cm]{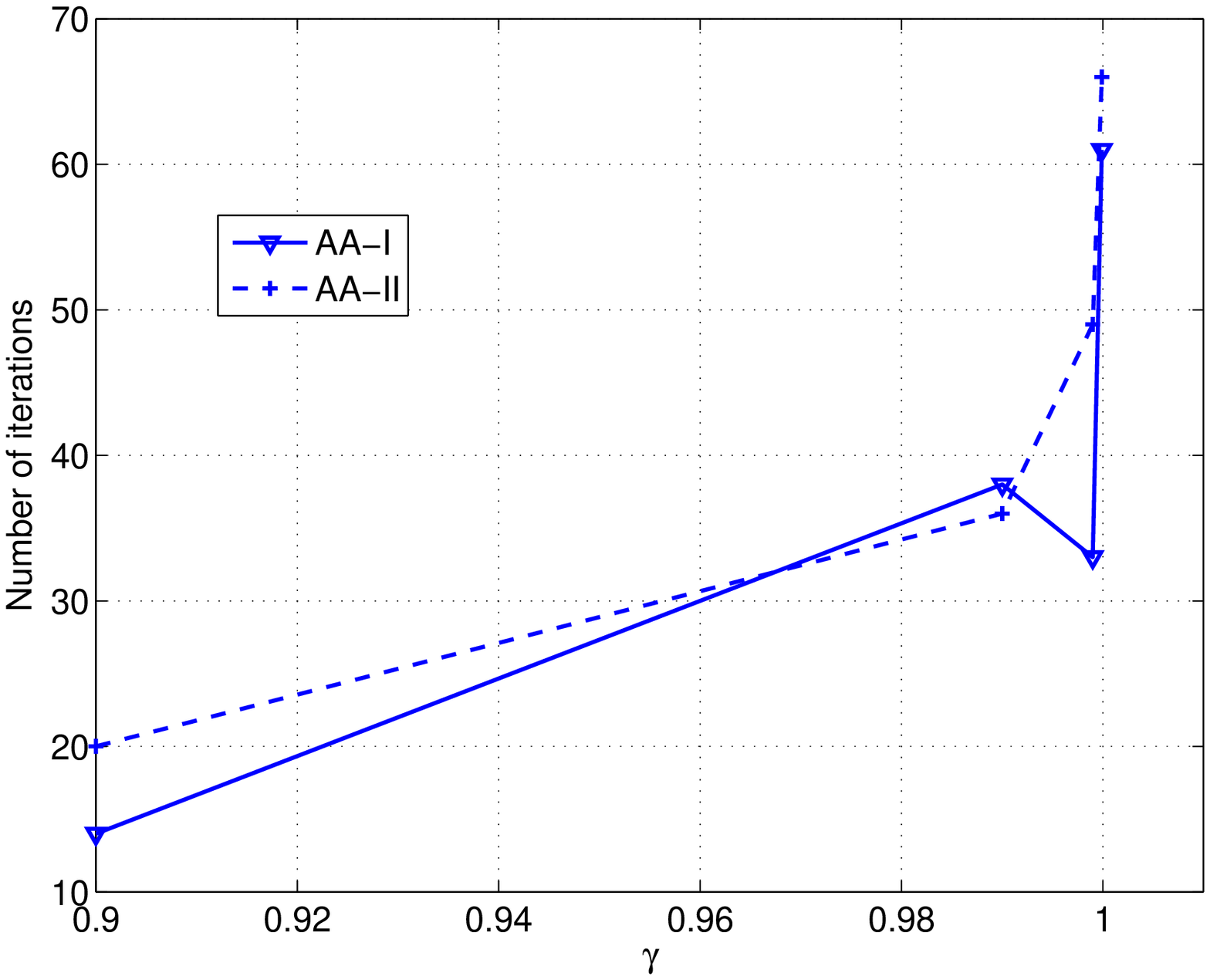} }
\caption{Solitary wave generation of (\ref{E15}), $c_{s}=0.75$, $\alpha=\beta=\delta=1$.  Number of iterations required  to achieve a residual error (in Euclidean norm) below $TOL=1E-13$ as function of $\gamma$. (a) MPE, RRE, VEA and TEA. (b) AA-I and AA-II.} \label{Ben2}
\end{figure}

To this end, we fix a speed $c_{s}=0.75$, the parameters $\alpha=\beta=\delta=1$ and generate numerically a solitary wave solution of (\ref{E15}) by combining the \PM method, standing for the family of iterative method (\ref{mm2}), along with the acceleration techniques considered in previous examples. We will take four values of $\gamma$, namely $0.9, 0.99, 0.999, 0.9999$ (which, for the considered values of the parameters, are close to $\gamma^{*}$, equals $1$ in our example), correspond to a more and more oscillating profile (with smaller and smaller amplitude, see Figures \ref{Ben1}(a)-(d); the computational window is $[-512,512]$ with $N=4096$ collocation points) and for which the \PM method with numerical continuation requires a long computation to converge or directly does not work. In all the experiments the initial iteration is the (analytically known) solitary wave profile corresponding to $\gamma=0$ (KdV equation).
As in the previous examples, we first estimate the performance of the acceleration techniques by comparing the number of iterations required by each of them to achieve a residual error (\ref{fsec32}) less than a tolerance $TOL=10^{-13}$.
\begin{table}
\begin{center}
\begin{tabular}{|c|c|c|c|c|}
\hline\hline  $\kappa$&MPE($\kappa$)&RRE($\kappa$)&VEA($\kappa$)&TEA($\kappa$)\\\hline
$2$&$35$&$35$&$33$&$36$\\
&($7.0418E-14$)&($5.3048E-14$)&($7.3270E-14$)&($7.5208E-14$)\\
$3$&$26$&$28$&$28$&$25$\\
&($8.5843E-14$)&($3.7474E-14$)&($4.3193E-14$)&($3.5514E-14$)\\
$4$&$25$&$27$&$23$&$21$\\
&($8.7162E-14$)&($4.0422E-14$)&($2.6305E-14$)&($2.7078E-14$)\\
$5$&$22$&$22$&$25$&$25$\\
&($2.3778E-14$)&($2.4540E-14$)&($2.4551E-14$)&($2.6612E-14$)\\
$6$&$19$&$19$&$24$&$24$\\
&($6.8569E-14$)&($6.3344E-14$)&($8.5839E-14$)&($8.8986E-14$)\\
$7$&$19$&$19$&$21$&$24$\\
&($4.8718E-15$)&($4.4709E-15$)&($5.1694E-14$)&($7.2978E-14$)\\
$8$&$21$&$21$&$19$&$19$\\
&($2.4289E-15$)&($2.3872E-15$)&($1.8080E-14$)&($6.5383E-14$)\\
$9$&$23$&$23$&$21$&$21$\\
&($2.2138E-15$)&($2.5529E-15$)&($3.5794E-15$)&($3.3103E-15$)\\
$10$&$25$&$25$&$23$&$23$\\
&($2.3084E-15$)&($2.2764E-15$)&($4.3858E-15$)&($3.3545E-15$)\\
\hline\hline
\end{tabular}
\end{center}
\caption{Solitary wave generation of (\ref{E15}) . Number of iterations required by MPE, RRE, VEA and TEA as function of $\kappa$ to achieve a residual error (\ref{fsec32}) below $TOL=10^{-13}$. The residual error at the last computed iterate is in parenthesis; $c_{s}=0.75$, $\alpha=\beta=\delta=1$, $\gamma=0.9$.}\label{tav7b}
\end{table}
For the case of the VEM and the four values of $\gamma$ considered, this information is given in Tables \ref{tav7b}-\ref{tav10b}. All the methods achieve convergence in the four cases. (The \PM method with continuation is not able to converge for the last two values of $\gamma$ and for the first two values the number of iterations required is prohibitive: for example, just going from $\gamma=0.98$ to $\gamma=0.99$ the method requires $266$ iterations to have a residual error of size $9.0634E-14$; the continuation process from the initial $\gamma=0$, where our computations start, requires a total number of iterations of about $4470$.) As expected the effort of VEM in number of iterations increases with $\gamma$, that is, with the oscillating character of the profile, see Figure \ref{Ben2}(a). 
\begin{table}
\begin{center}
\begin{tabular}{|c|c|c|c|c|}
\hline\hline  $\kappa$&MPE($\kappa$)&RRE($\kappa$)&VEA($\kappa$)&TEA($\kappa$)\\\hline
$2$&$77$&$67$&$147$&$93$\\
&($4.3301E-14$)&($5.8218E-14$)&($1.8199E-14$)&($8.9274E-14$)\\
$3$&$56$&$53$&$47$&$51$\\
&($2.9502E-14$)&($4.4310E-14$)&($2.1373E-14$)&($6.5537E-14$)\\
$4$&$43$&$49$&$48$&$41$\\
&($4.7982E-14$)&($1.5692E-14$)&($7.4855E-14$)&($2.4388E-14$)\\
$5$&$38$&$36$&$39$&$37$\\
&($3.8410E-14$)&($2.5110E-14$)&($4.0062E-14$)&($3.4320E-14$)\\
$6$&$33$&$33$&$43$&$43$\\
&($3.6079E-14$)&($2.3031E-14$)&($1.6073E-14$)&($1.6545E-14$)\\
$7$&$30$&$30$&$35$&$49$\\
&($6.3053E-14$)&($3.2262E-14$)&($3.3824E-14$)&($1.8129E-15$)\\
$8$&$38$&$31$&$40$&$37$\\
&($9.5443E-14$)&($6.3679E-14$)&($3.2121E-14$)&($1.4194E-14$)\\
$9$&$34$&$34$&$41$&$41$\\
&($1.8749E-15$)&($2.1169E-15$)&($4.1825E-15$)&($1.8746E-15$)\\
$10$&$37$&$37$&$45$&$45$\\
&($2.2503E-15$)&($1.7597E-15$)&($4.7245E-15$)&($1.8872E-15$)\\
\hline\hline
\end{tabular}
\end{center}
\caption{Solitary wave generation of (\ref{E15}) . Number of iterations required by MPE, RRE, VEA and TEA as function of $\kappa$ to achieve a residual error (\ref{fsec32}) below $TOL=10^{-13}$. The residual error at the last computed iterate is in parenthesis; $c_{s}=0.75$, $\alpha=\beta=\delta=1$, $\gamma=0.99$.}\label{tav8b}
\end{table}
Among them and except some particular cases (for example, when TEA is applied with $\gamma=0.9999$ and $\kappa=6$) the polynomial methods are more efficient than $\epsilon$-algorithms when $\gamma$ increases, although the difference is shorter than that was obtained in the examples of Section \ref{se31}. It is remarkable that in the case of a solitary wave profile with a small number of oscillations (for example, when $\gamma=0.9$) the performance of the methods is virtually the same: after one or two cycles, the improvement of the acceleration technique is good enough to not needing to complete the next cycle in order to achieve the tolerance for the residual error. This is particularly emphasized in the case of the $\epsilon$-algorithms, where the cycle is longer ($2\kappa$ against $\kappa+1$ for the case of the polynomial methods).
\begin{table}
\begin{center}
\begin{tabular}{|c|c|c|c|c|}
\hline\hline  $\kappa$&MPE($\kappa$)&RRE($\kappa$)&VEA($\kappa$)&TEA($\kappa$)\\\hline
$3$&$63$&$60$&$56$&$65$\\
&($2.2431E-14$)&($3.9572E-14$)&($9.8315E-14$)&($6.6026E-14$)\\
$4$&$79$&$49$&$77$&$70$\\
&($4.9749E-14$)&($4.1392E-14$)&($6.6180E-14$)&($9.3318E-14$)\\
$5$&$43$&$43$&$54$&$101$\\
&($5.8270E-14$)&($1.0725E-14$)&($5.9093E-14$)&($4.4323E-14$)\\
$6$&$45$&$45$&$71$&$71$\\
&($9.4355E-14$)&($3.5010E-14$)&($6.8677E-14$)&($1.9108E-14$)\\
$7$&$37$&$37$&$65$&$65$\\
&($8.3020E-14$)&($6.3139E-14$)&($1.1800E-14$)&($2.2238E-15$)\\
$8$&$51$&$41$&$73$&$73$\\
&($4.8702E-15$)&($1.3712E-14$)&($8.5215E-15$)&($3.7373E-15$)\\
$9$&$45$&$45$&$81$&$81$\\
&($2.8855E-15$)&($2.7603E-15$)&($3.4268E-15$)&($2.6722E-15$)\\
$10$&$51$&$51$&$111$&$69$\\
&($3.3110E-14$)&($3.2354E-14$)&($2.0694E-15$)&($4.0083E-14$)\\
\hline\hline
\end{tabular}
\end{center}
\caption{Solitary wave generation of (\ref{E15}) . Number of iterations required by MPE, RRE, VEA and TEA as function of $\kappa$ to achieve a residual error (\ref{fsec32}) below $TOL=10^{-13}$. The residual error at the last computed iterate is in parenthesis; $c_{s}=0.75$, $\alpha=\beta=\delta=1$, $\gamma=0.999$.}\label{tav9b}
\end{table}
\begin{table}
\begin{center}
\begin{tabular}{|c|c|c|c|c|}
\hline\hline  $\kappa$&MPE($\kappa$)&RRE($\kappa$)&VEA($\kappa$)&TEA($\kappa$)\\\hline
$3$&$60$&$56$&$65$&$49$\\
&($4.0201E-14$)&($2.3990E-14$)&($1.2911E-14$)&($1.4934E-14$)\\
$4$&$55$&$61$&$203$&$51$\\
&($2.4890E-14$)&($4.0300E-14$)&($6.8061E-14$)&($1.4970E-14$)\\
$5$&$45$&$42$&$52$&$73$\\
&($1.4175E-14$)&($6.6916E-14$)&($5.7187E-14$)&($1.2493E-15$)\\
$6$&$43$&$41$&$101$&$43$\\
&($3.6348E-14$)&($1.5143E-14$)&($7.7548E-14$)&($4.9096E-14$)\\
$7$&$46$&$46$&$65$&$65$\\
&($4.7457E-15$)&($2.3306E-15$)&($9.9131E-14$)&($8.3133E-16$)\\
$8$&$51$&$51$&$109$&$73$\\
&($1.9157E-15$)&($7.6397E-16$)&($3.0772E-15$)&($1.3583E-15$)\\
$9$&$56$&$56$&$101$&$81$\\
&($6.7580E-16$)&($7.1896E-16$)&($7.8767E-16$)&($2.4959E-15$)\\
$10$&$61$&$61$&$199$&$89$\\
&($7.8194E-16$)&($9.0591E-16$)&($7.7246E-16$)&($6.4547E-16$)\\
\hline\hline
\end{tabular}
\end{center}
\caption{Solitary wave generation of (\ref{E15}) . Number of iterations required by MPE, RRE, VEA and TEA as function of $\kappa$ to achieve a residual error(\ref{fsec32}) below $TOL=10^{-13}$. The residual error at the last computed iterate is in parenthesis; $c_{s}=0.75$, $\alpha=\beta=\delta=1$, $\gamma=0.9999$.}\label{tav10b}
\end{table}
The results from Tables \ref{tav7b}-\ref{tav10b} are in contrast with those from Tables \ref{tav11b}-\ref{tav14b}, that correspond to the AAM. The main conclusion here is that these methods are strongly affected by the oscillating character of the profiles and, compared to VEM, do not seem to be recommendable for this sort of computations, at least without a suitable choice of preconditioning. (Some of it was suggested by the previous experiments concerning the generalized solitary waves of some Boussinesq systems, see Section \ref{se3}.) Tables \ref{tav11b}-\ref{tav14b} show that as $\gamma$ increases, ill-conditioning of the corresponding  least-squares problem is observed from even moderate values of $nw$ in the case of AA-I (affecting the stability of the method, which is not able to converge or requires a great effort in number of iterations) while AA-II is not so affected.
\begin{table}
\begin{center}
\begin{tabular}{|c|c|c|}
\hline\hline  $nw$&AA-I($nw$)&AA-II($nw$)\\\hline
$1$&$20$($5.5496E-14$)&$25$($4.1169E-14$)\\
$2$&$17$($1.9530E-15$)&$21$($1.2546E-14$)\\
$3$&$15$($4.4608E-15$)&$20$($4.8161E-14$)\\
$4$&$15$($1.3401E-14$)&$20$($8.3292E-14$)\\
$5$&$14$($5.1185E-14$)&$21$($1.4961E-14$)\\
$6$&$14$($1.5211E-14$)&$21$($1.4958E-14$)\\
\hline\hline
\end{tabular}
\end{center}
\caption{Solitary wave generation of (\ref{E15}) . Number of iterations required by AA-I and AA-II as function of $nw$ to achieve a residual error (\ref{fsec32}) below $TOL=10^{-13}$. The residual error at the last computed iterate is in parenthesis. $c_{s}=0.75$, $\alpha=\beta=\delta=1$, $\gamma=0.9$.}\label{tav11b}
\end{table}
However, when AAM work, they exhibit a competitive performance, as shown in Figure \ref{Ben2}(b) when comparing with Figure \ref{Ben2}(a). (Our implementation follows that of described in \cite{walkern}, which uses the unconstrained form of the least-squares problem and that also was suggested by some other authors, \cite{fangs}. For the numerical resolution we have used some other alternatives, with $QR$ decomposition with pivoting, \cite{walkern}, and SVD, \cite{fangs}. The results of Tables \ref{tav11b}-\ref{tav14b} correspond to the first implementation, while the second one overcomes ill-conditioning in some more cases of AA-I, but at the cost of an important increase of the iterations.)
\begin{table}
\begin{center}
\begin{tabular}{|c|c|c|}
\hline\hline  $nw$&AA-I($nw$)&AA-II($nw$)\\\hline
$1$&$58$($8.4209E-14$)&$38$($6.6459E-14$)\\
$2$&$44$($3.1097E-14$)&$44$($5.6598E-14$)\\
$3$&$38$($1.0963E-14$)&$36$($9.8103E-14$)\\
$4$&$53$($2.9810E-14$)&$47$($3.1570E-14$)\\
$5$&&$57$($4.4383E-14$)\\
$6$&&$47$($5.3103E-14$)\\
$7$&&$53$($7.2250E-14$)\\
$8$&&$64$($3.6054E-14$)\\
\hline\hline
\end{tabular}
\end{center}
\caption{Solitary wave generation of (\ref{E15}) . Number of iterations required by AA-I and AA-II as function of $nw$ to achieve a residual error (\ref{fsec32}) below $TOL=10^{-13}$. The residual error at the last computed iterate is in parenthesis. $c_{s}=0.75$, $\alpha=\beta=\delta=1$, $\gamma=0.99$.}\label{tav12b}
\end{table}
\begin{table}
\begin{center}
\begin{tabular}{|c|c|c|}
\hline\hline  $nw$&AA-I($nw$)&AA-II($nw$)\\\hline
$1$&$63$($8.2392E-14$)&$77$($9.2229E-14$)\\
$2$&$57$($9.5654E-14$)&$83$($8.6540E-14$)\\
$3$&$225$($4.7416E-14$)&$49$($7.4595E-14$)\\
$4$&$33$($2.0366E-14$)&$91$($6.5926E-14$)\\
$5$&&$64$($1.1616E-14$)\\
$6$&&$70$($5.9306E-14$)\\
\hline\hline
\end{tabular}
\end{center}
\caption{Solitary wave generation of (\ref{E15}) . Number of iterations required by AA-I and AA-II as function of $nw$ to achieve a residual error (\ref{fsec32}) below $TOL=10^{-13}$. The residual error at the last computed iterate is in parenthesis. $c_{s}=0.75$, $\alpha=\beta=\delta=1$, $\gamma=0.999$.}\label{tav13b}
\end{table}
\begin{table}
\begin{center}
\begin{tabular}{|c|c|c|}
\hline\hline  $nw$&AA-I($nw$)&AA-II($nw$)\\\hline
$1$&$800$($2.7029E-14$)&$84$($9.3995E-15$)\\
$2$&$61$($3.1858E-14$)&$193$($3.5965E-14$)\\
$3$&$99$($2.2611E-14$)&$98$($1.3640E-14$)\\
$4$&$144$($8.0198E-14$)&$82$($4.5365E-14$)\\
$5$&&$94$($1.8114E-14$)\\
$6$&&$66$($8.1966E-14$)\\
\hline\hline
\end{tabular}
\end{center}
\caption{Solitary wave generation of (\ref{E15}) . Number of iterations required by AA-I and AA-II as function of $nw$ to achieve a residual error (\ref{fsec32}) below $TOL=10^{-13}$. The residual error at the last computed iterate is in parenthesis. $c_{s}=0.75$, $\alpha=\beta=\delta=1$, $\gamma=0.9999$.}\label{tav14b}
\end{table}

\subsubsection{Lump solitary waves of 2D Benjamin equation}
In order to finish off this example we study the performance of the acceleration techniques when generating numerically lump solitary wave solutions of the 2D Benjamin equation \cite{kim,kima1,kima2}
\begin{eqnarray}
\label{ben2d0}
\left(\eta_{t}+ \alpha\eta\eta_{x}-\beta \mathcal{H}(\eta_{xx})+\delta\eta_{xxx} \right)_{x}-\eta_{zz}=0,
\end{eqnarray}
where $\alpha, \beta, \delta\geq 0$ and $\mathcal{H}$ is the Hilbert transform (\ref{hilb}) with respect to $x$. In (\ref{ben2d0}), as in the one-dimensional case, $\eta=\eta(x,z,t)$ stands for the interfacial deviation wave between two ideal fluids with a bounded upper
layer and the heavier one with infinite depth, and under the
presence of interfacial tension. The two-dimensional version
incorporates weak transverse variations. For the experiments below we will consider a normalized version of (\ref{ben2d0}), \cite{kim}
\begin{eqnarray}
\label{ben2d1}
\left(\eta_{t}+ (\eta^{2})_{x}-2\Gamma \mathcal{H}(\eta_{xx})+\eta_{xxx} \right)_{x}-\eta_{zz}=0,
\end{eqnarray}
where $\Gamma\geq 0$. (The case $\Gamma=0$ corresponds to the KP-I equation, \cite{kadomtsevp}.) For localized solutions, the zero total mass
\begin{eqnarray}
\int_{-\infty}^{\infty} \eta(x,z,t)dx=0,\label{ben2d2}
\end{eqnarray}
is also assumed. Lump solitary wave solutions of (\ref{ben2d1}), (\ref{ben2d2}) are solutions of the form $\eta(x,z,t)=\eta(X,Z), X=x-c_{s}t,  Z=z$ for some $c_{s}>0$. Substitution into (\ref{ben2d1}) leads to
\begin{eqnarray}
\label{lumpsw} \left(-c_{s}\eta+\eta^{2}-2\Gamma
\mathcal{H}(\eta_{X})+\eta_{XX}\right)_{XX}-\eta_{ZZ}=0,
\end{eqnarray}
As shown in \cite{kima2}, the value $\Gamma=1$ marks a bifurcation point as for the type of lump solutions of (\ref{ben2d1}) between lumps of KP-I type and of wavepacket type. This implies in particular that as $\Gamma<1$ approaches one the lump wave increases the oscillations.

The numerical procedure used in \cite{kim,kima2} to generate approximate lump waves combines numerical continuation in $\Gamma$, pseudospectral approximation to (\ref{lumpsw}) (where constraint (\ref{ben2d2}) is imposed) and Newton's method for the resolution of the corresponding system of equations in each step of the $\Gamma$-homotopic path. The use of the \PM methods (instead of Newton's) was suggested in \cite{alvarezd}. (For the use of the \PM method in the generation of two-dimensional solitary waves see e.~g. \cite{AbramyanS1985,VoronovichSS1998}.) As in the one-dimensional case, the computation of approximate lump profiles comes up two main difficulties: the use of numerical continuation and the oscillating behaviour of the lump.  These problems can be overcome with the use of acceleration techniques, especially VEM.

In order to illustrate this we will take $c_{s}=1$ and generate approximate lump solitary waves for $\Gamma=0.99, 0.999, 0.9999$. As described in \cite{alvarezd}, the periodic problem on a square $[-L_{x},L_{x}]\times [-L_{z},L_{z}]$ of (\ref{ben2d1}) is discretized by using a Fourier collocation method, generating approximations $(\eta_{h})_{i,j}$ to the lump profile $\eta(x_{i},z_{j})$ at the collocation points $x_{i}=-L_{x}+ih_{x}, z_{j}=-L_{z}+jh_{z}, h_{x}=2L_{x}/N_{x}, h_{z}=2L_{z}/N_{z}, i=1,\ldots, N_{x}, j=1,\ldots, N_{z}$. The system for the discrete Fourier coefficients of the approximation is of the form
\begin{eqnarray}
(k_{x}^{2}(c_{s}+2\Gamma |k_{x}|+k_{z}^{2}))\widehat{\eta}_{h}(k_{x},k_{z})=k_{x}^{2}\left(\widehat{\eta_{h}^{2}}\right)(k_{x},k_{z}),\label{ben2d3}
\end{eqnarray}
for $k_{x}=-N_{x}/2,\ldots,N_{x}/2, k_{z}=-N_{z}/2,\ldots,N_{z}/2$ and where $\widehat{\eta}_{h}(k_{x},k_{z})$ stands for discrete $(k_{x},k_{z})$-Fourier component of $\eta_{h}$. The zero total mass condition (\ref{ben2d2}) is imposed as
\begin{eqnarray}
\widehat{\eta}_{h}(0,0)=0.\label{ben2d4}
\end{eqnarray}
When (\ref{ben2d4}) is included into (\ref{ben2d3}), the resulting system for the rest of Fourier components is nonsingular and it is iteratively solved, for fixed $\Gamma$, by using:
\begin{itemize}
\item[(i)] The \PM method with numerical continuation from the initial iteration given by the exact profile for $\Gamma=0$ 
\begin{eqnarray}
\eta_{0}(x,z)=12c_{s}\frac{3+c_{s}^{2}z^{2}-c_{s}x^{2}}{(3+c_{s}x^{2}+c_{s}^{2}z^{2})^{2}}.\label{IL}
\end{eqnarray}
\item[(ii)] The \PM method without numerical continuation but accelerated with the six techniques MPE, RRE, TEA, VEA, AA-I and AA-II and the same initial iteration (\ref{IL}).
\end{itemize} 

The experiments below follow a similar design to that of the one dimensional case. We have taken $N_{x}=N_{z}=1024$ with $L_{x}=L_{z}=256$ and a tolerance of $TOL=10^{-8}$ for the control of the iteration. As before, the number of iterations shown in the numerical results correspond to the total account, including the iterations of each cycle. From this value, one can obtain the iterations exclusively due to the corresponding acceleration. We think that this way of counting the iterations makes the comparison with the results without acceleration more realistic.
\begin{figure}[htbp]
\centering \subfigure{\subfigure{
\includegraphics[width=5.5cm]{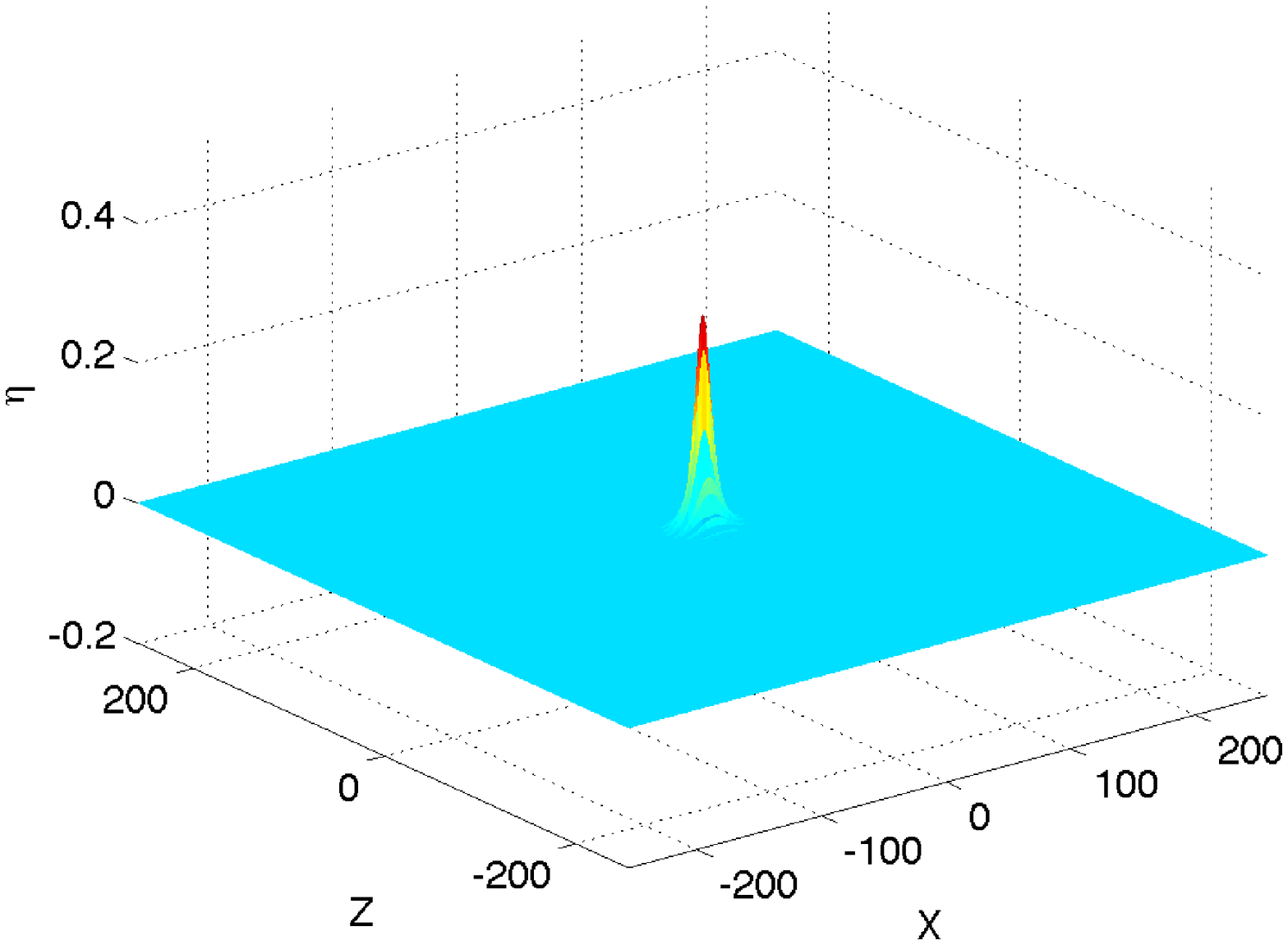}}
\subfigure{
\includegraphics[width=5.5cm]{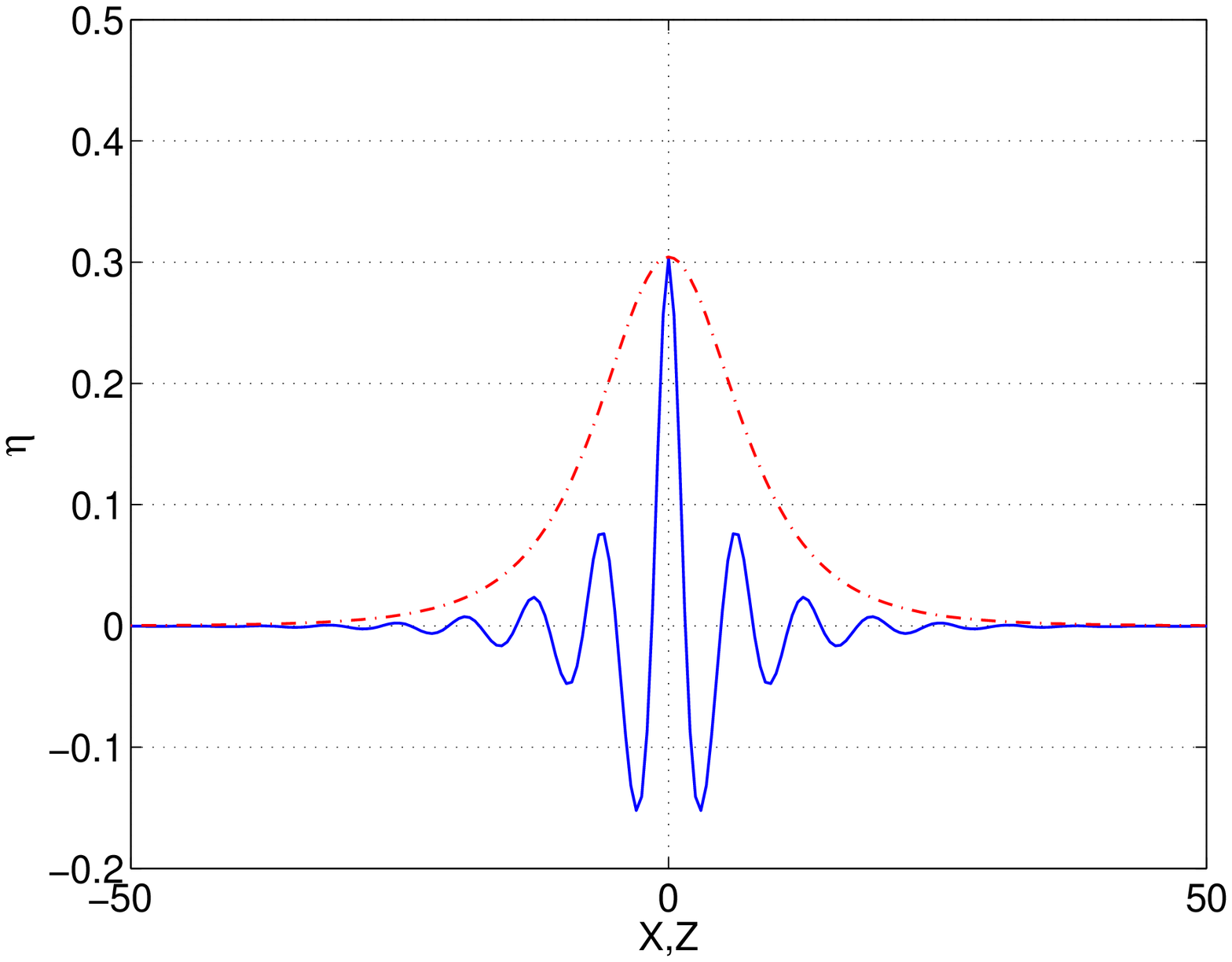}}
}
\end{figure}
\begin{figure}[htbp]
\centering \subfigure{\subfigure{
\includegraphics[width=5.5cm]{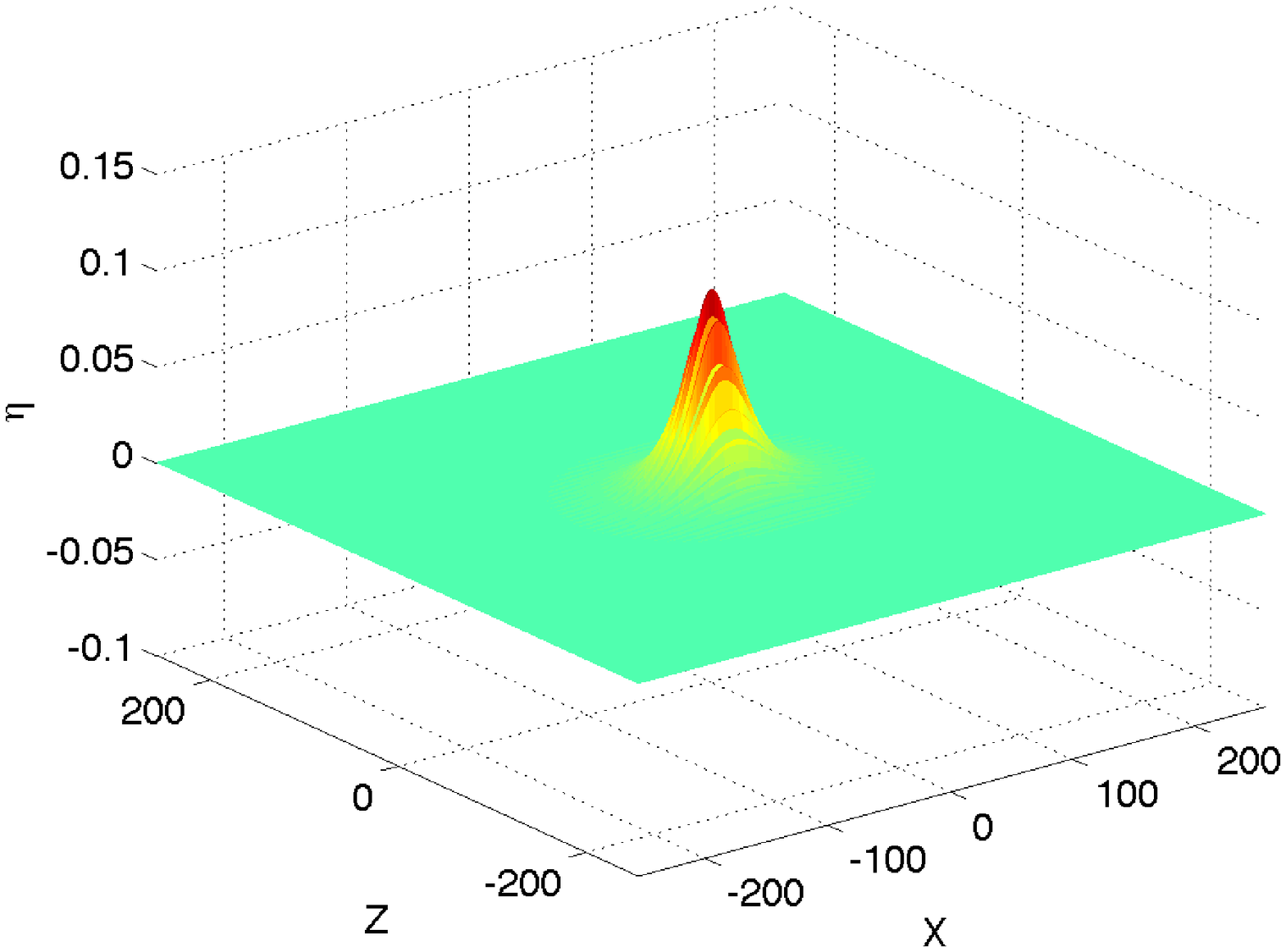}}
\subfigure{
\includegraphics[width=5.5cm]{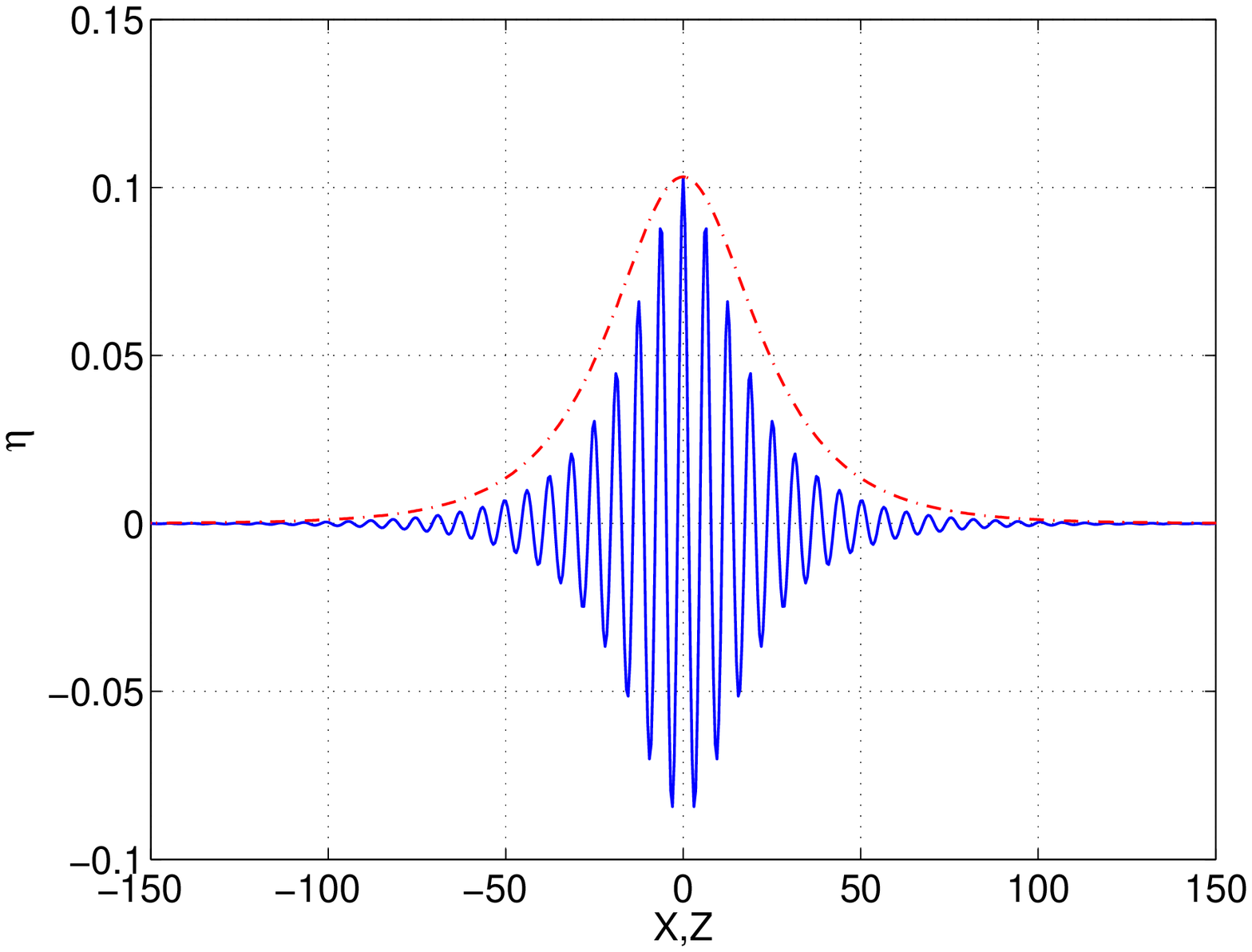}} }
\end{figure}
\begin{figure}[htbp]
\centering \subfigure{\subfigure{
\includegraphics[width=5.5cm]{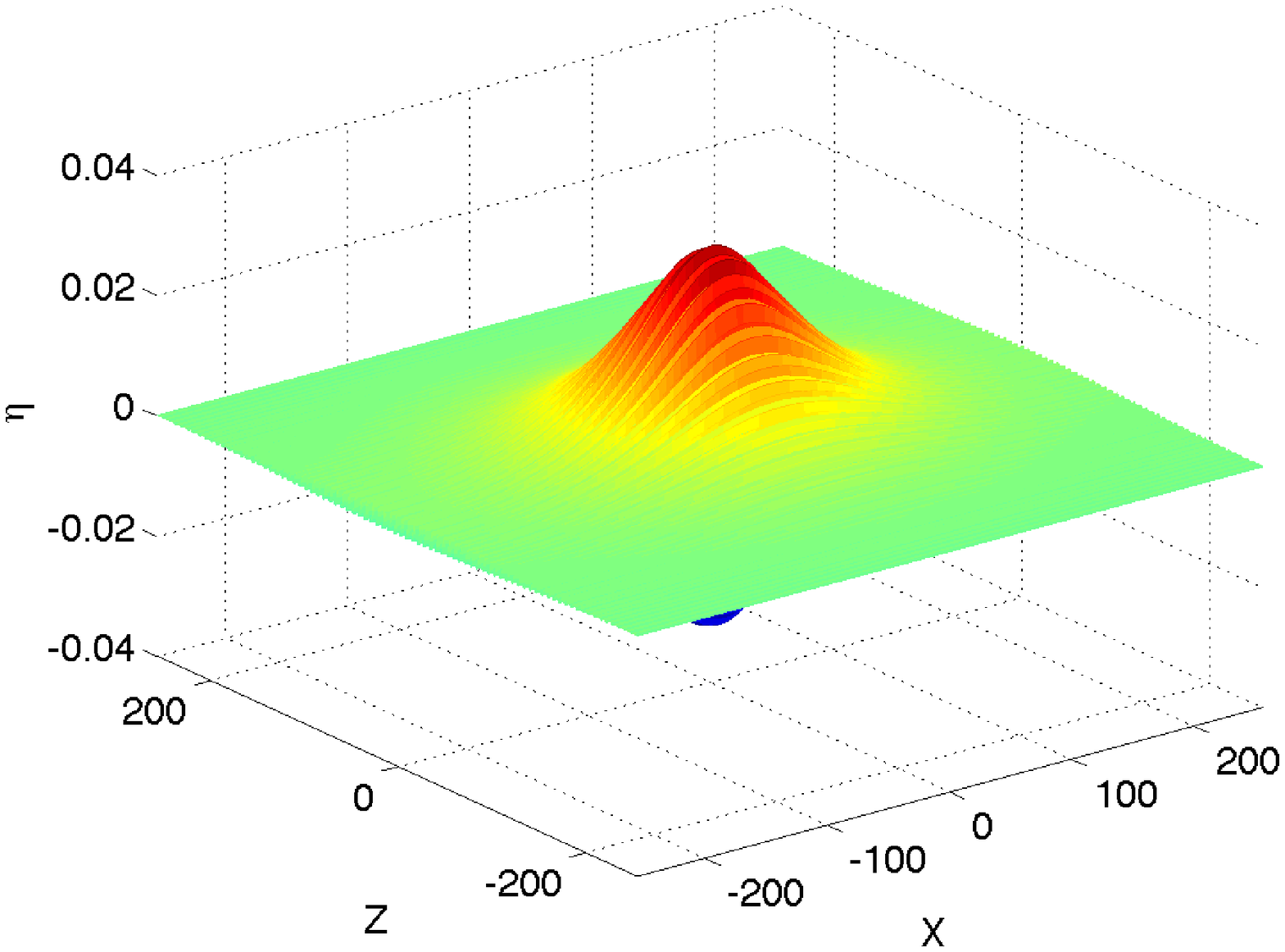}}
\subfigure{
\includegraphics[width=5.5cm]{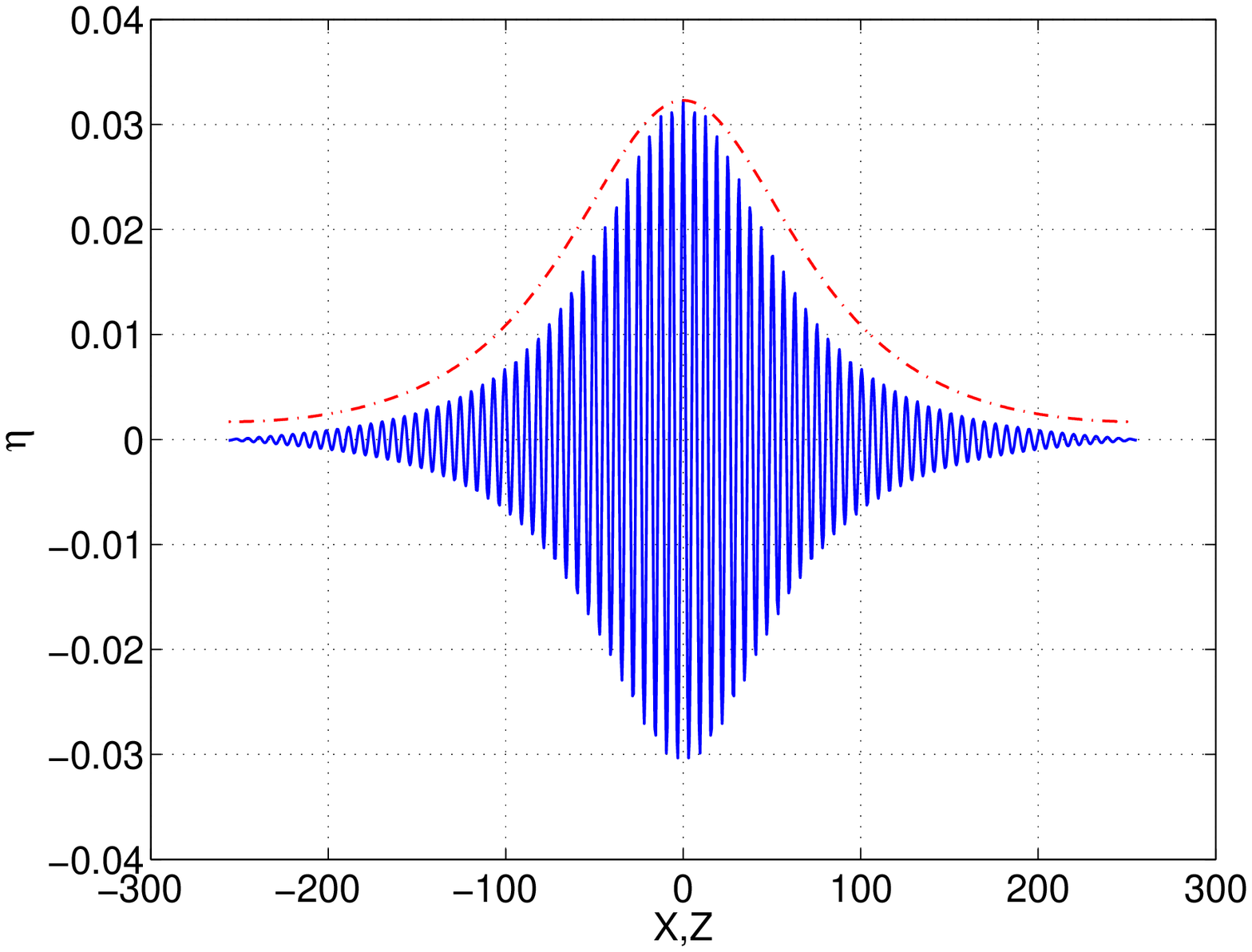}} }
\caption{Solitary wave generation of (\ref{ben2d1}) with \PM
method, accelerated with MPE. The approximated profiles correspond to $\Gamma=0.99,0.999,0.9999$ (left).
On the right, the corresponding $X$ and $Z$ cross sections {are shown} (solid
and dashed-dotted lines, respectively) .} \label{fexample41}
\end{figure}
We also remark that the use of the same initial iteration (\ref{IL})
is against the alternative technique with acceleration methods since, according to the form of the resulting waves, the initial profile is not close. This should be observed in the behaviour of the residual error with respect to the number of iterations: the main effort is at the beginning; once the error is small enough, all the techniques accelerate the convergence in a more important way. Here three values of $\Gamma=0.99, 0.999, 0.9999$ are considered. The corresponding approximate waves (confirming the highly oscillatory behaviour) are computed with the acceleration procedure of (ii) and can be observed in Figure \ref{fexample41}.
The first procedure in (i), based on continuation with respect to $\Gamma$, is totally inefficient for the the first value and does not work for the other two. The performance of the acceleration techniques is compared in Figure \ref{Ben2d} and Tables \ref{tav15b}-\ref{tav18b}.

\begin{table}
\begin{center}
\begin{tabular}{|c|c|c|c|c|}
\hline\hline  $\kappa$&MPE($\kappa$)&RRE($\kappa$)&VEA($\kappa$)&TEA($\kappa$)\\\hline
$2$&$38$&$41$&$63$&$49$\\
&($7.5320E-09$)&($7.7910E-09$)&($5.8233E-09$)&($8.1460E-09$)\\
$3$&$32$&$75$&$49$&$56$\\
&($6.6515E-09$)&($3.3504E-09$)&($9.0167E-09$)&($9.4419E-09$)\\
$4$&$32$&$37$&$33$&$38$\\
&($5.0660E-09$)&($3.8941E-09$)&($9.4694E-09$)&($8.9305E-09$)\\
$5$&$26$&$30$&$35$&$34$\\
&($4.9957E-09$)&($8.4036E-09$)&($1.3179E-09$)&($1.1643E-09$)\\
$6$&$29$&$30$&$30$&$27$\\
&($4.5762E-09$)&($1.4663E-09$)&($9.7151E-09$)&($2.8493E-09$)\\
$7$&$25$&$33$&$31$&$31$\\
&($2.8813E-09$)&($5.5247E-09$)&($1.8370E-09$)&($6.0375E-09$)\\
$8$&$21$&$29$&$35$&$35$\\
&($6.8132E-09$)&($3.7396E-09$)&($6.7544E-09$)&($2.3207E-09$)\\
$9$&$22$&$25$&$35$&$40$\\
&($7.0811E-09$)&($6.1487E-09$)&($7.6658E-09$)&($5.3006E-09$)\\
$10$&$23$&$25$&$39$&$43$\\
&($9.6176E-10$)&($7.8664E-09$)&($9.0031E-09$)&($2.7713E-09$)\\
\hline\hline
\end{tabular}
\end{center}
\caption{Solitary wave generation of (\ref{ben2d1}) . Number of iterations required by MPE, RRE, VEA and TEA as function of $\kappa$ to achieve a residual error below $TOL=10^{-8}$. The residual error (\ref{fsec32}) at the last computed iterate is in parenthesis; $c_{s}=1$, $\Gamma=0.99$.}\label{tav15b}
\end{table}

\begin{table}
\begin{center}
\begin{tabular}{|c|c|c|c|c|}
\hline\hline  $\kappa$&MPE($\kappa$)&RRE($\kappa$)&VEA($\kappa$)&TEA($\kappa$)\\\hline
$3$&$47$&$54$&$59$&$66$\\
&($8.2832E-09$)&($5.6681E-09$)&($9.5469E-09$)&($7.6937E-09$)\\
$4$&$46$&$53$&$87$&$65$\\
&($9.8170E-09$)&($8.4899E-09$)&($7.9785E-09$)&($1.6182E-09$)\\
$5$&$43$&$44$&$56$&$56$\\
&($2.7604E-09$)&($6.3987E-09$)&($7.2821E-09$)&($6.2222E-09$)\\
$6$&$36$&$38$&$105$&$67$\\
&($5.3518E-09$)&($5.7851E-09$)&($1.2954E-10$)&($3.0156E-10$)\\
$7$&$41$&$43$&$76$&$77$\\
&($3.7360E-09$)&($3.1947E-09$)&($1.1064E-10$)&($4.9070E-10$)\\
$8$&$38$&$56$&$86$&$87$\\
&($1.9864E-09$)&($1.9460E-09$)&($4.8067E-11$)&($4.6454E-10$)\\
$9$&$51$&$51$&$96$&$59$\\
&($4.6989E-12$)&($2.2425E-11$)&($2.8033E-11$)&($6.7102E-09$)\\
$10$&$56$&$56$&$85$&$65$\\
&($5.5395E-10$)&($1.0988E-11$)&($2.7617E-09$)&($1.3631E-09$)\\
\hline\hline
\end{tabular}
\end{center}
\caption{Solitary wave generation of (\ref{ben2d1}) . Number of iterations required by MPE, RRE, VEA and TEA as function of $\kappa$ to achieve a residual error below $TOL=10^{-8}$. The residual error (\ref{fsec32}) at the last computed iterate is in parenthesis; $c_{s}=1$, $\Gamma=0.999$.}\label{tav16b}
\end{table}

\begin{table}
\begin{center}
\begin{tabular}{|c|c|c|c|c|}
\hline\hline  $\kappa$&MPE($\kappa$)&RRE($\kappa$)&VEA($\kappa$)&TEA($\kappa$)\\\hline
$3$&$53$&$54$&$57$&$72$\\
&($5.7252E-09$)&($6.4717E-09$)&($4.9757E-09$)&($8.5993E-09$)\\
$4$&$51$&$47$&$82$&$65$\\
&($6.5144E-09$)&($6.5364E-09$)&($1.2292E-09$)&($7.8352E-10$)\\
$5$&$43$&$49$&$67$&$68$\\
&($1.4420E-09$)&($5.7277E-09$)&($9.3104E-09$)&($3.3648E-09$)\\
$6$&$50$&$50$&$105$&$105$\\
&($2.8215E-12$)&($1.3107E-10$)&($4.7779E-12$)&($3.2453E-09$)\\
$7$&$49$&$50$&$91$&$122$\\
&($8.1696E-09$)&($5.7674E-09$)&($1.4031E-11$)&($6.5282E-10$)\\
$8$&$55$&$55$&$103$&$155$\\
&($2.5994E-10$)&($9.0192E-09$)&($1.4896E-10$)&($2.9075E-10$)\\
\hline\hline
\end{tabular}
\end{center}
\caption{Solitary wave generation of (\ref{ben2d1}) . Number of iterations required by MPE, RRE, VEA and TEA as function of $\kappa$ to achieve a residual error below $TOL=10^{-8}$. The residual error (\ref{fsec32}) at the last computed iterate is in parenthesis; $c_{s}=1$, $\Gamma=0.9999$.}\label{tav17b}
\end{table}

\begin{table}
\begin{center}
\begin{tabular}{|c|c|c|}
\hline\hline  $nw$&AA-I($nw$)&AA-II($nw$)\\\hline
$2$&$23$($2.5840E-09$)&$25$($5.5612E-09$)\\
$3$&&$20$($2.5548E-09$)\\
$4$&$36$($1.5755E-09$)&$18$($8.3616E-09$)\\
$5$&$18$($6.8363E-09$)&$17$($6.2358E-09$)\\
$6$&$19$($4.5057E-09$)&$17$($6.4962E-09$)\\
$7$&$19$($3.3613E-09$)&$17$($3.1330E-09$)\\
\hline\hline
\end{tabular}
\end{center}
\caption{Solitary wave generation of (\ref{ben2d1}) . Number of iterations required by AA-I and AA-II as function of $nw$ to achieve a residual error below $TOL=10^{-8}$. The residual error (\ref{fsec32}) at the last computed iterate is in parenthesis. $c_{s}=1$,  $\Gamma=0.99$.}\label{tav18b}
\end{table}

\begin{figure}[htbp]
\centering \subfigure[]{
\includegraphics[width=6.6cm]{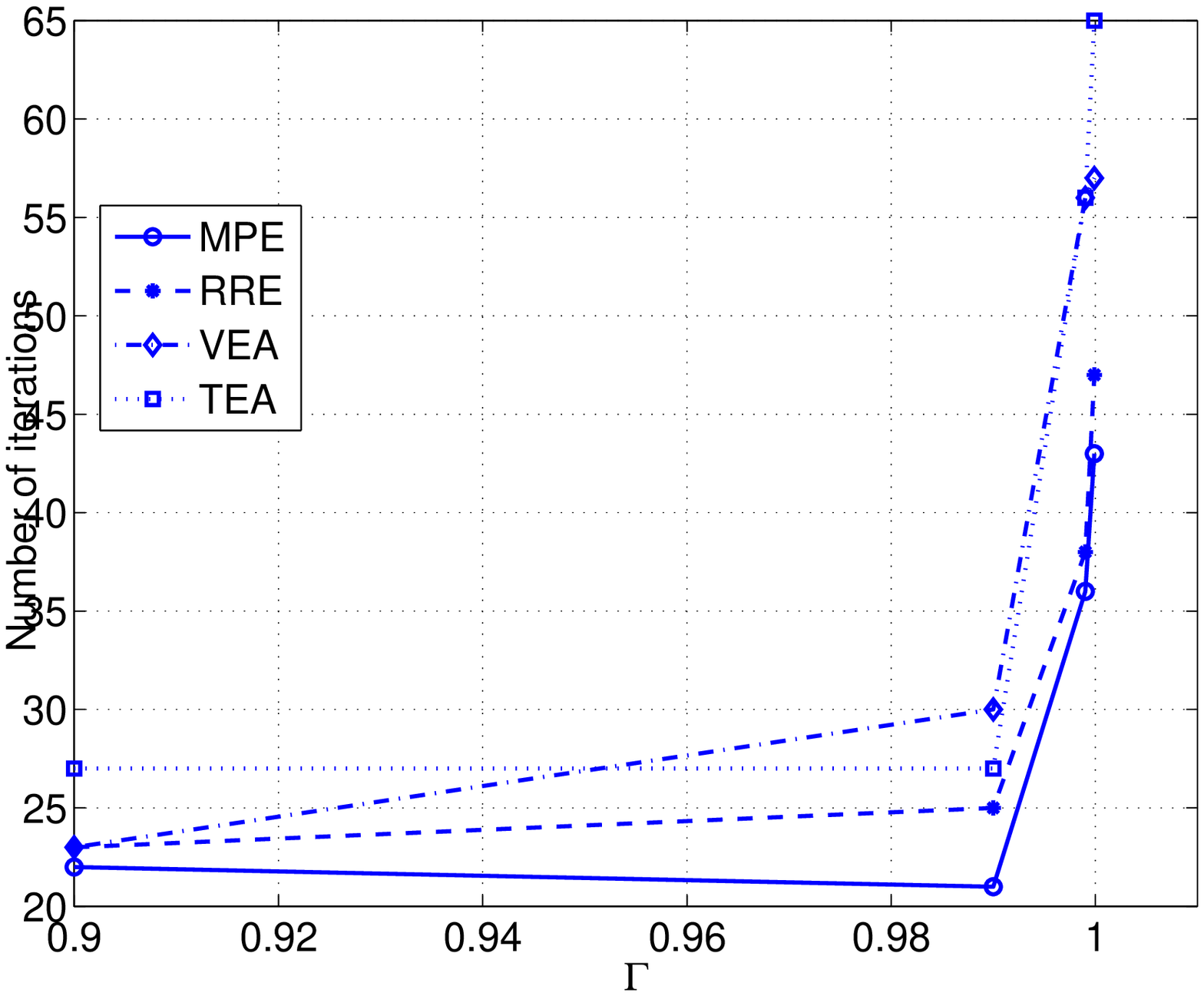} }
\subfigure[]{
\includegraphics[width=6.6cm]{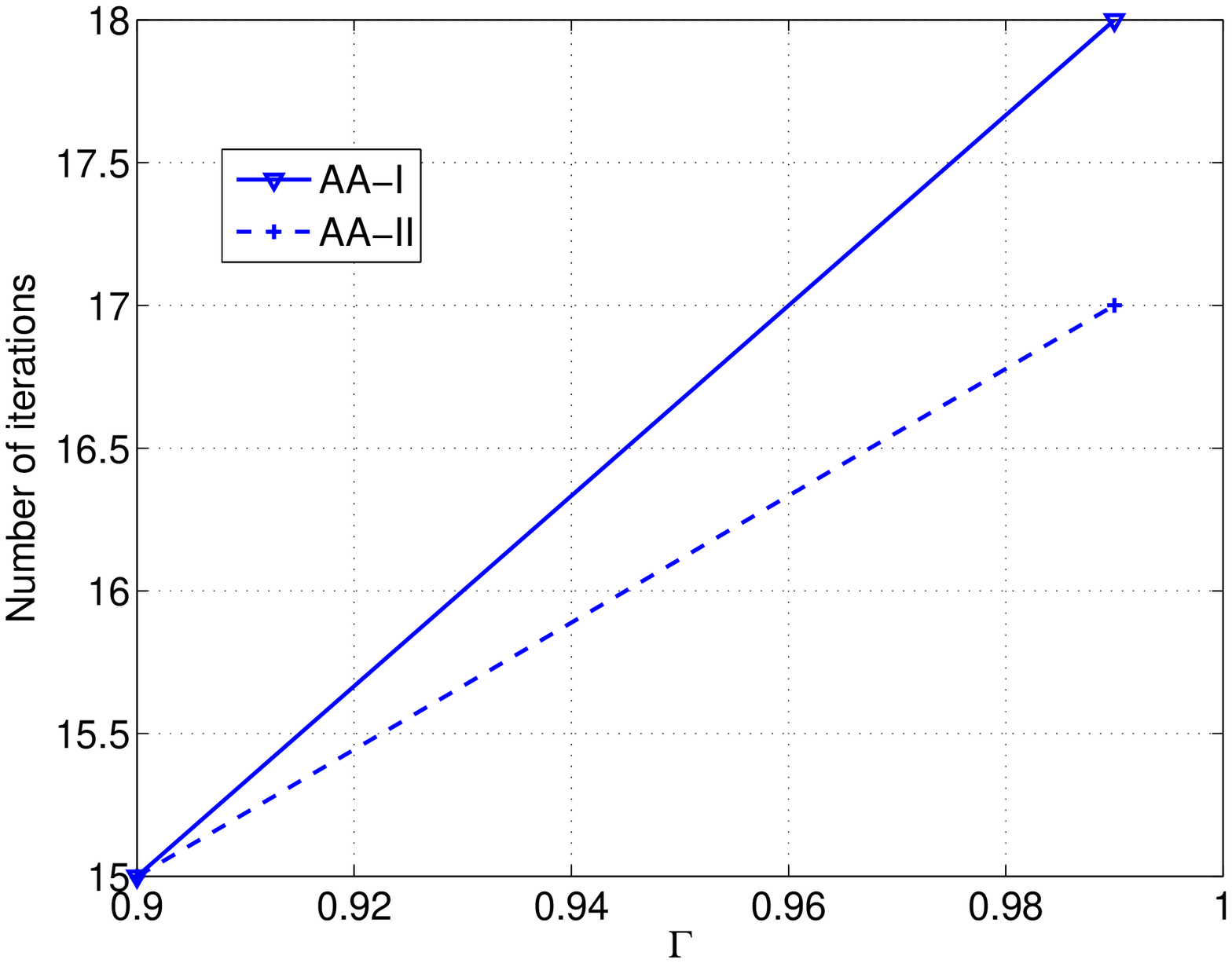} }
\caption{Solitary wave generation of (\ref{ben2d1}), $c_{s}=1$.  Number of iterations required  to achieve a residual error (in Euclidean norm) below $TOL=10^{-8}$ as function of $\Gamma$. (a) MPE, RRE, VEA and TEA. (b) AA-I and AA-II.} \label{Ben2d}
\end{figure}
The comparison of the techniques in this case confirms the conclusions obtained in the one-dimensional version, namely:
\begin{itemize}
\item The best performance is given by the polynomial methods (MPE in this case).
\item The $\epsilon$-algorithms, although less efficient, are also competitive (contrary to what was observed in some previous examples).
\item AAM only work correctly up to a moderate value of $\Gamma<1$. When $\Gamma$ approaches one they cannot get the performance of VEM or directly fail.
\end{itemize}
\section{Acceleration techniques with extended \PM type methods}
\label{se4}
One of the drawbacks of the \PM type methods (\ref{mm2}) in traveling wave generation is their limitation to some specific problems, namely those with homogeneous nonlinearities. When the nonlinear term is not homogeneous but a combination of homogeneous functions of different degree, these methods can be extended by adapting the stabilizing function $s$ to each homogeneous part. This leads to the so-called e-\PM type methods, derived in \cite{alvarezd2014b}. In this section and in order to improve the traveling wave generation for problems with this type of nonlinearities, we will apply the acceleration techniques to the e-\PM method as initial iterative procedure. This will be illustrated with the numerical generation of
localized ground state solutions of
the following generalized nonlinear Schr\"{o}dinger equation
\begin{eqnarray}
iu_{t}+u_{xx}-V(x)u+|u|^{2}u-0.2|u|^{4}u+\nu
|u|^{6}u=0,\label{gnls2}
\end{eqnarray}
with $V(x)=-3.5{\rm sech}^{2}(x+1.5)-3{\rm sech}^{2}(x-1.5)$ and $\nu$ a real constant. Equation (\ref{gnls2}) was studied in \cite{Yang2012} (see also references therein) where a bifurcation of solitary waves for $\nu= \nu_{c}\approx 0.01247946$ was analyzed. The bifurcation is of transcritical type with two tangentially connected branches of smooth solutions. This can be characterized by using the behaviour of the power
\begin{eqnarray}
P(\mu)=\int_{-\infty}^{\infty} U^{2}(x,\mu)dx,\label{power}
\end{eqnarray}
as function of $\mu$ for any localized ground state solution $u(x,t)=U(x,\mu)e^{i\mu t}, \mu\in\mathbb{R}$. The two branches are connected at some $(\mu_{0},P(\mu_{0}))\approx (3.28,14.35)$.
The numerical generation of localized ground state profiles of (\ref{gnls2}) with e-\PM type methods was treated in \cite{alvarezd2014b} where the equation for the profiles $U(x,\mu)$
\begin{eqnarray*}
-\mu U+u^{\prime\prime}-V(x)U+|U|^{2}U-0.2|U|^{4}U+\nu
|U|^{6}U=0,\label{gnls2b}
\end{eqnarray*}
was discretized by Fourier collocation techniques, leading to the system $LU_{h}=N(U_{h})$ for the vector approximation $U_{h}$ at the grid points $x_{j}$ and where
\begin{eqnarray}
L&=&\mu I-D_{h}^{2}+diag(V(x_{0}),\ldots,V(x_{m-1})),\nonumber\\
N(U_{h})&=&N_{1}(U_{h})+N_{2}(U_{h})+N_{3}(U_{h})\nonumber\\
&=&\left(|U_{h}|.^{2}\right).U_{h}-0.2\left(|U_{h}|.^{4}\right).U_{h}+\kappa
\left(|U_{h}|.^{6}\right).U_{h}.\label{lab41}
\end{eqnarray}
The nonlinearity in (\ref{lab41}) contains three homogeneous terms with degrees $p_{1}=3, p_{2}=5, p_{3}=7$ and the e-\PM method
\begin{eqnarray}
  LU_{h}^{n+1}&=&\sum_{j=1}^{3}s_{j}(U_{h}^{n})N_{j}(U_{h}^{n}), n=0,1,\ldots,\label{lab22e}\\
s_{j}(u)&=&\left(\frac{\langle Lu,u\rangle}{\langle N(u),u\rangle}\right)^{\gamma_{j}},\quad \gamma_{j}=\frac{p_{j}}{p_{j}-1},\quad j=1,2,3,\label{lab25e}
\end{eqnarray}
is applied. The iteration (\ref{lab22e}), (\ref{lab25e}) will be considered as the method to be complemented with acceleration techniques. Finally, the quantity (\ref{power}) has been approximated by
\begin{eqnarray}
P_{h}(U_{h})=h\sum_{j} U_{h,j}^{2}.\label{power2}
\end{eqnarray}
\begin{table}
\begin{center}
\begin{tabular}{|c|c|}\hline\hline
Classical fixed point&e-\PM { method} (\ref{lab22e}), (\ref{lab25e})\\\hline
1.687048E+00&9.829607E-01\\
9.834930E-01&4.740069E-01\\
4.793766E-01&3.616157E-01\\
3.747266E-01&2.606251E-01+1.734293E-01i\\
1.979766E-01&2.606251E-01-1.734293E-01i\\
1.426054E-01&1.764488E-01\\\hline
\end{tabular}
\end{center}
\caption{Six largest magnitude eigenvalues of the iteration matrices
of classical fixed point algorithm and of e-\PM method (\ref{lab22e}), (\ref{lab25e}) for $\mu=3.281$ at the last computed iterate. The dominant egienvalue in the column on the right justifies the slow performance of the method.}\label{tav_epet1}
\end{table}
\begin{table}
\begin{center}
\begin{tabular}{|c|c|c|c|c|}
\hline\hline  $\kappa$&MPE($\kappa$)&RRE($\kappa$)&VEA($\kappa$)&TEA($\kappa$)\\\hline
$2$&$185$&$260$&$1250$&$187$\\
&($5.3669E-11$)&($9.9728E-11$)&($9.5278E-11$)&($5.5118E-11$)\\
$3$&$135$&$135$&$209$&$155$\\
&($5.1486E-11$)&($7.4460E-11$)&($8.6758E-11$)&($6.5522E-11$)\\
$4$&$118$&$64$&$167$&$109$\\
&($6.3944E-11$)&($7.6714E-11$)&($8.7013E-11$)&($5.6136E-11$)\\
$5$&$64$&$69$&$75$&$78$\\
&($7.0759E-11$)&($8.7100E-11$)&($7.1288E-11$)&($8.7351E-11$)\\
$6$&$55$&$65$&$85$&$795$\\
&($8.4602E-11$)&($4.5311E-11$)&($7.6708E-11$)&($6.9120E-11$)\\
$7$&$53$&$58$&$80$&$91$\\
&($6.6155E-11$)&($1.1617E-11$)&($9.8698E-11$)&($3.4911E-11$)\\
$8$&$49$&$58$&$89$&$70$\\
&($5.6418E-11$)&($9.0426E-11$)&($9.0902E-11$)&($6.9110E-11$)\\
$9$&$52$&$62$&$82$&$78$\\
&($1.3763E-11$)&($4.1423E-11$)&($7.4598E-11$)&($7.1601E-11$)\\
$10$&$47$&$67$&$88$&$106$\\
&($6.1296E-11$)&($4.6194E-11$)&($6.9450E-11$)&($4.9290E-11$)\\
\hline\hline
\end{tabular}
\end{center}
\caption{Ground state generation of (\ref{gnls2}). Number of iterations required by MPE, RRE, VEA and TEA as function of $\kappa$ to achieve a residual error below $TOL=10^{-10}$. The residual error (\ref{fsec32}) at the last computed iterate is in parenthesis; $\mu=3.281$. For the e-\PM method (\ref{lab22e}), (\ref{lab25e}) without acceleration $n=1023$ iterations are required for a residual error of $9.9939E-11$.}\label{tav_epet2}
\end{table}

\begin{table}
\begin{center}
\begin{tabular}{|c|c|c|c|c|}
\hline\hline  $\kappa$&MPE($\kappa$)&RRE($\kappa$)&VEA($\kappa$)&TEA($\kappa$)\\\hline
$2$&$14.446162E+00$&$14.446162E+00$&$14.446162E+00$&$14.458882E+00$\\
$3$&$14.446162E+00$&$14.446162E+00$&$14.446162E+00$&$14.458882E+00$\\
$4$&$14.446162E+00$&$14.446162E+00$&$14.446162E+00$&$14.458882E+00$\\
$5$&$14.446162E+00$&$14.446162E+00$&$14.446162E+00$&$14.458882E+00$\\
$6$&$14.446162E+00$&$14.446162E+00$&$14.446162E+00$&$14.458882E+00$\\
$7$&$14.446162E+00$&$14.446162E+00$&$14.446162E+00$&$14.446162E+00$\\
$8$&$14.458882E+00$&$14.446162E+00$&$14.446162E+00$&$14.446162E+00$\\
$9$&$14.446162E+00$&$14.446162E+00$&$14.446162E+00$&$14.446162E+00$\\
$10$&$14.446162E+00$&$14.446162E+00$&$14.446162E+00$&$14.446162E+00$\\
\hline\hline
\end{tabular}
\end{center}
\caption{Ground state generation of (\ref{gnls2}). Values of (\ref{power2}) for each iteration from Table \ref{tav_epet2}; $\mu=3.281$.}\label{tav_epet2b}
\end{table}
\begin{figure}[htbp]
\centering 
\subfigure[]{
\includegraphics[width=6.6cm]{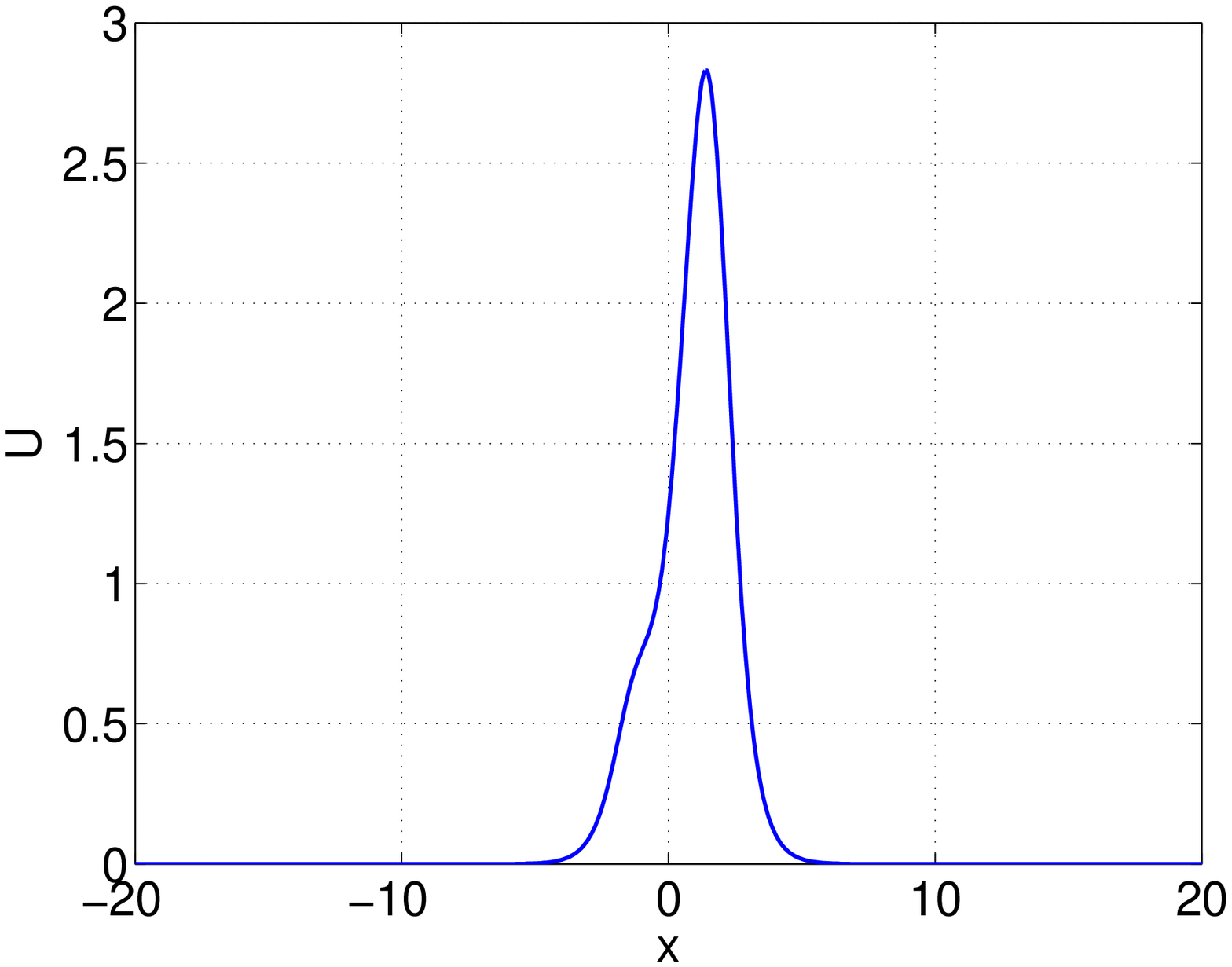} }
\subfigure[]{
\includegraphics[width=6.6cm]{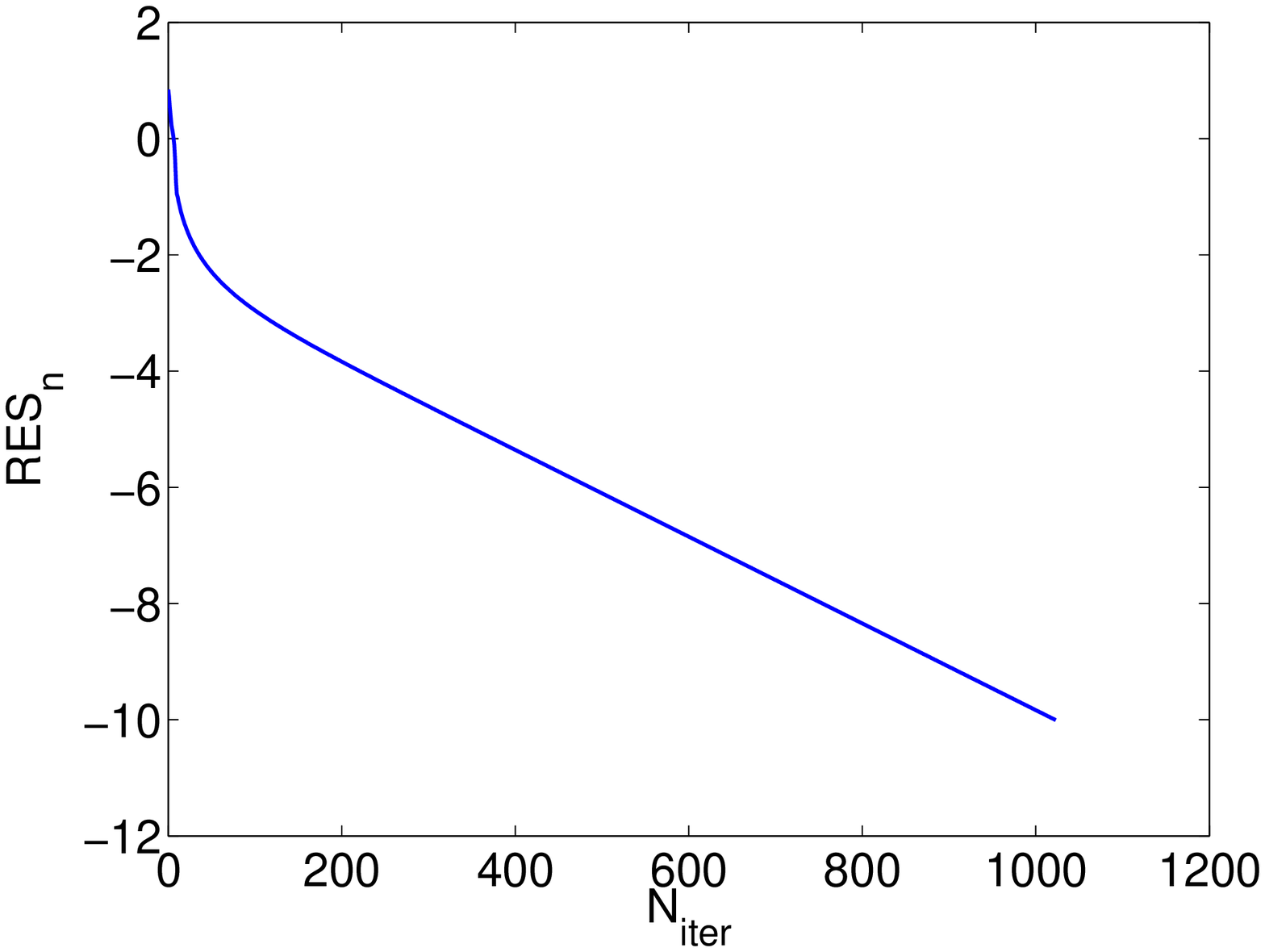} }
\caption{Numerical ground state generation of  (\ref{gnls2}) with $\mu=3.281$. (a) Approximate asymmetric profile; (b) Residual error as function of the number of iterations for the e-\PM method (\ref{lab22e}), (\ref{lab25e})  without acceleration. } \label{Figpet1}
\end{figure}
\begin{figure}[htbp]
\centering 
\subfigure[]{
\includegraphics[width=8.6cm]{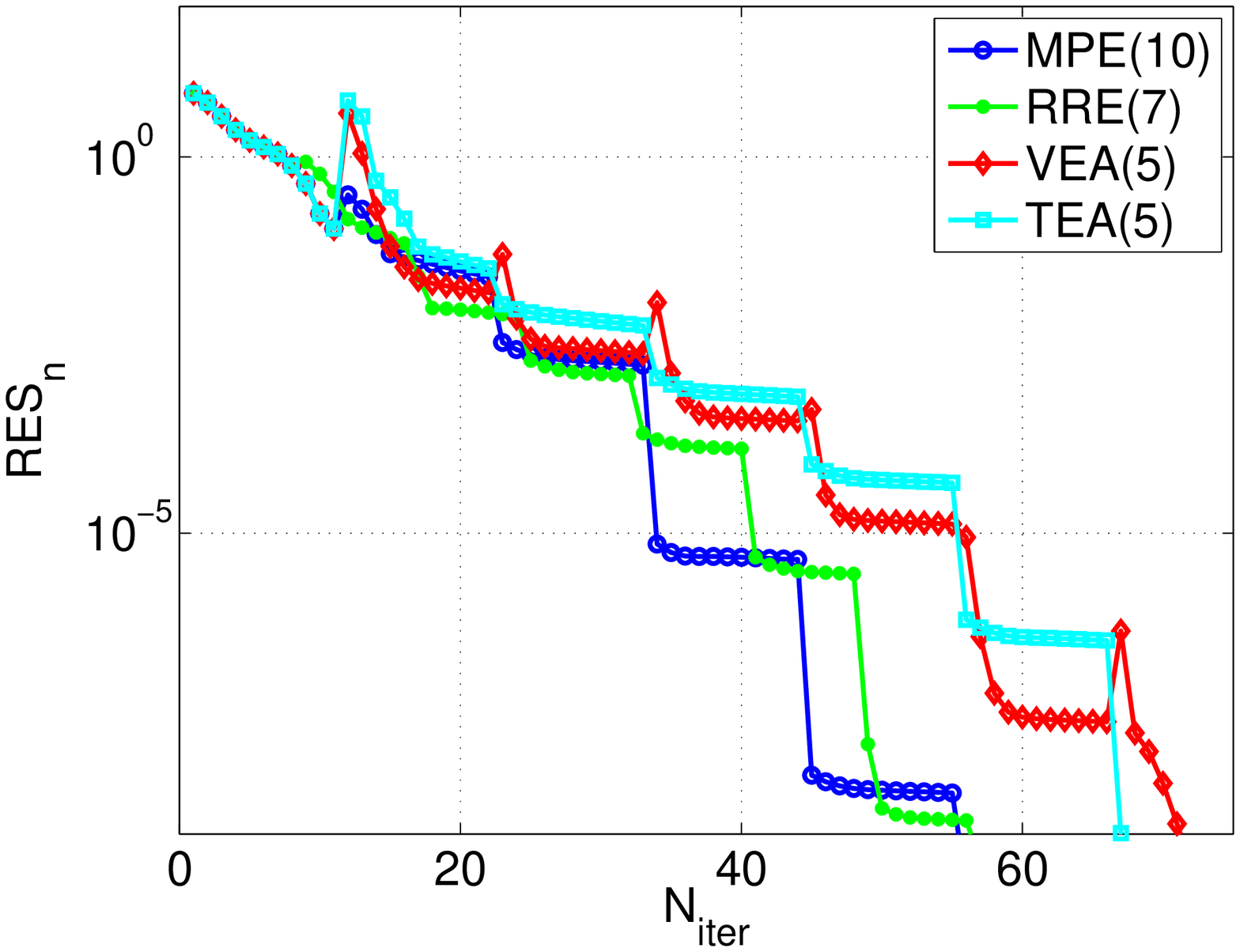} }
\subfigure[]{
\includegraphics[width=8.6cm]{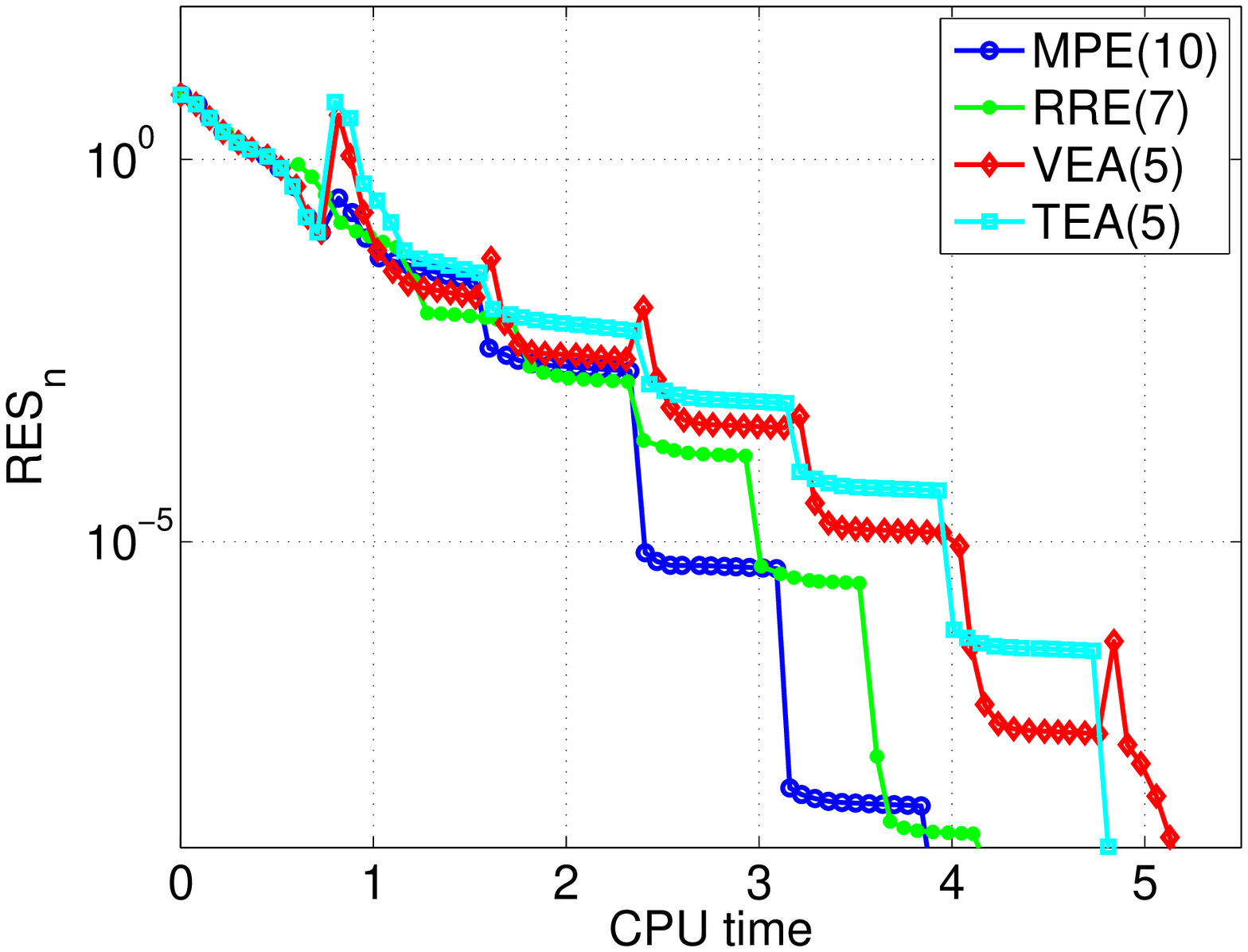} }
\caption{Numerical ground state generation of  (\ref{gnls2}) with $\mu=3.281$. Residual error (\ref{fsec32}) as function of the number of iterations (a) and CPU time in seconds (b) for the e-\PM method (\ref{lab22e}), (\ref{lab25e}) accelerated with MPE(10), RRE(7), VEA(5) and TEA(5). } \label{Figpet2}
\end{figure}

The numerical illustration of this case takes $\mu=3.281$, $TOL=10^{-10}$ and a superposition of squared hyperbolic secant functions as initial iteration. The numerical profile $U_{h}$ generated by (\ref{lab22e}), (\ref{lab25e}) is shown in Figure \ref{Figpet1}(a). The corresponding value for (\ref{power2}) is $P_{h}(U_{h})=14.446162$ and the poor performance of the method is made clear in Figure \ref{Figpet1}(b) which displays the behaviour of the residual error (\ref{fsec32}) as function of the number of iterations and shows that the method requires $n=1023$ iterations to achieve a residual error below $TOL$. (See Table \ref{tav_epet1}, first eigenvalue of the second column, to explain this slow behaviour.)

\begin{figure}[htbp]
\centering 
\subfigure[]{
\includegraphics[width=6.6cm]{epet1.eps} }
\subfigure[]{
\includegraphics[width=6.6cm]{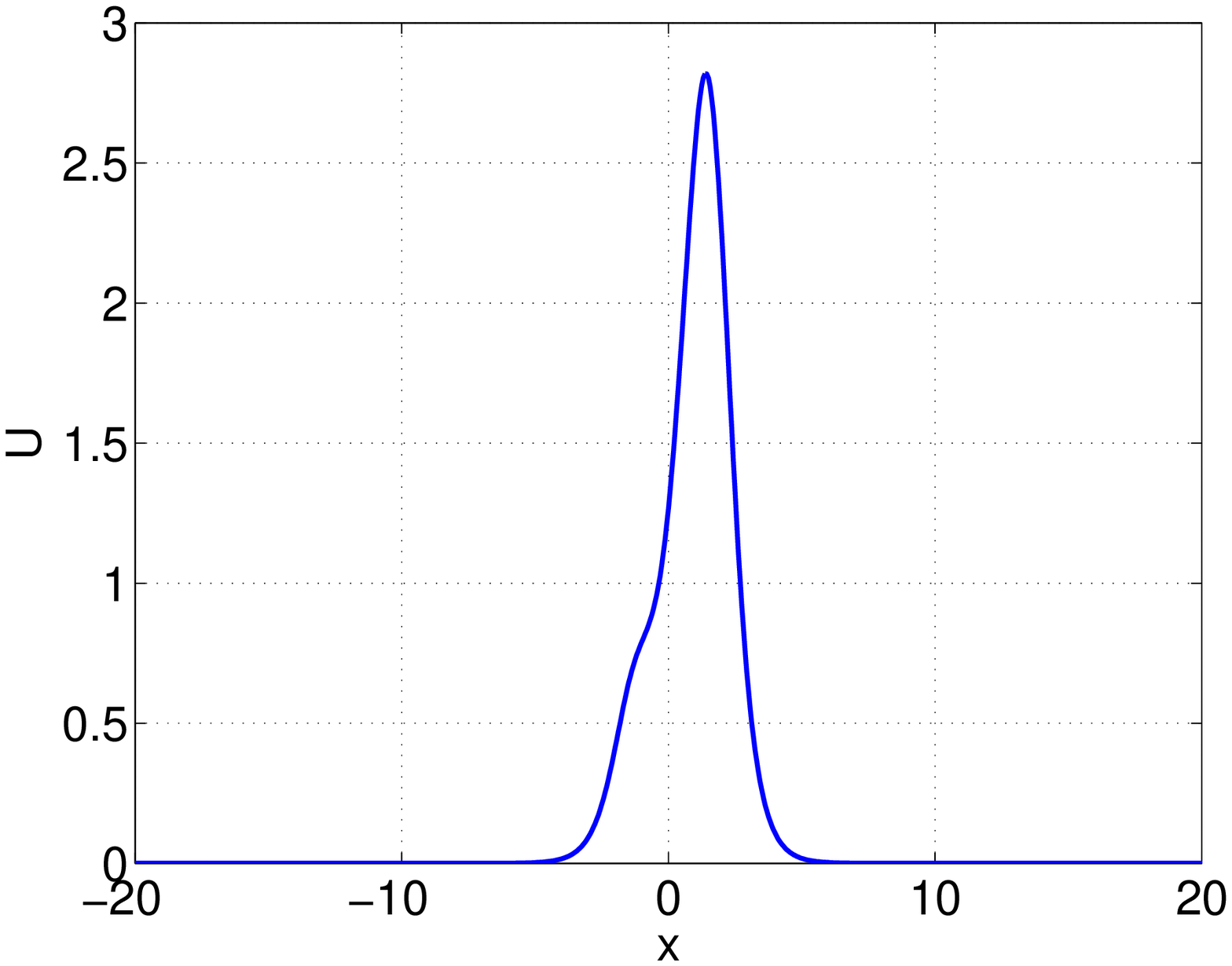} }
\subfigure[]{
\includegraphics[width=6.6cm]{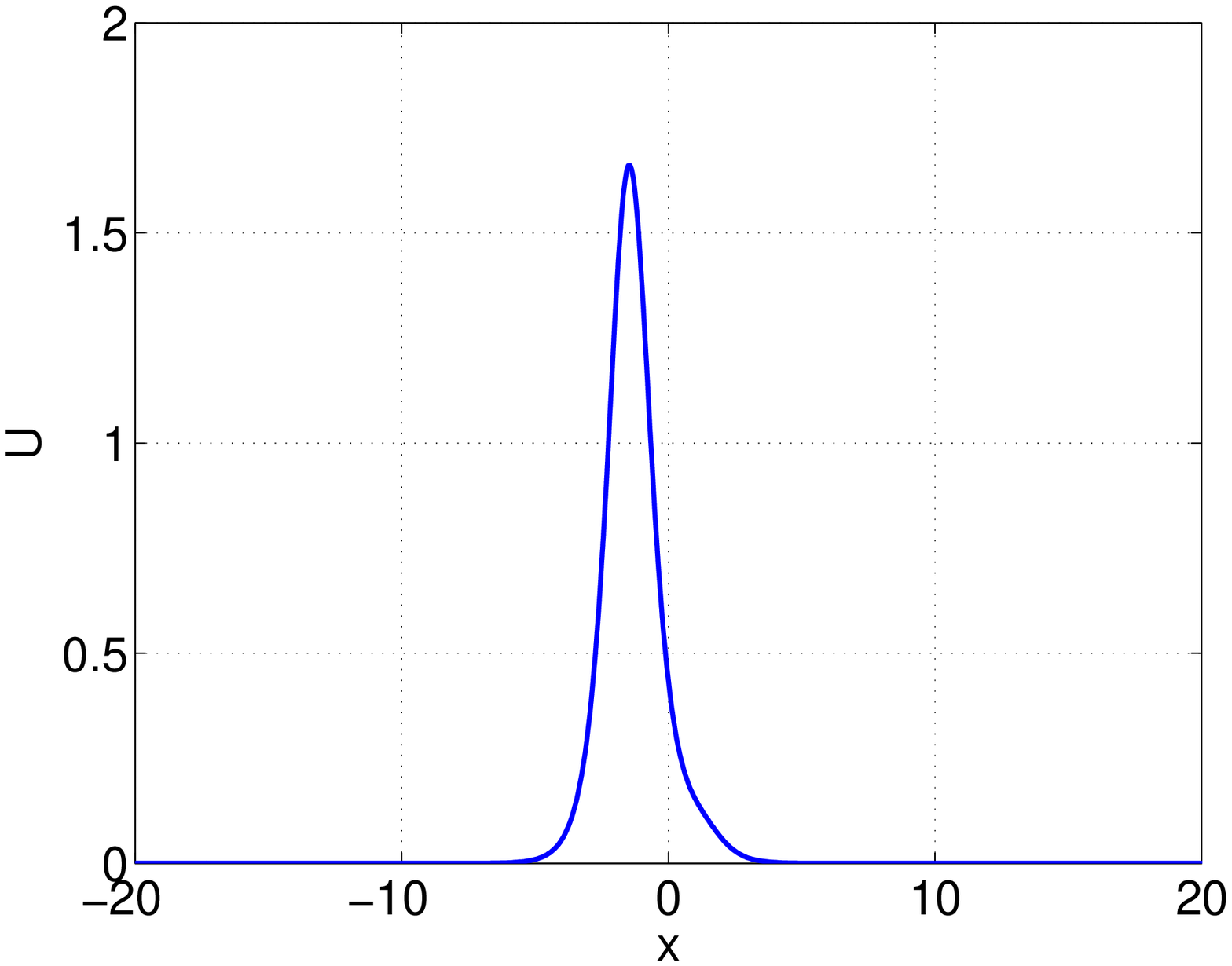} }
\subfigure[]{
\includegraphics[width=6.6cm]{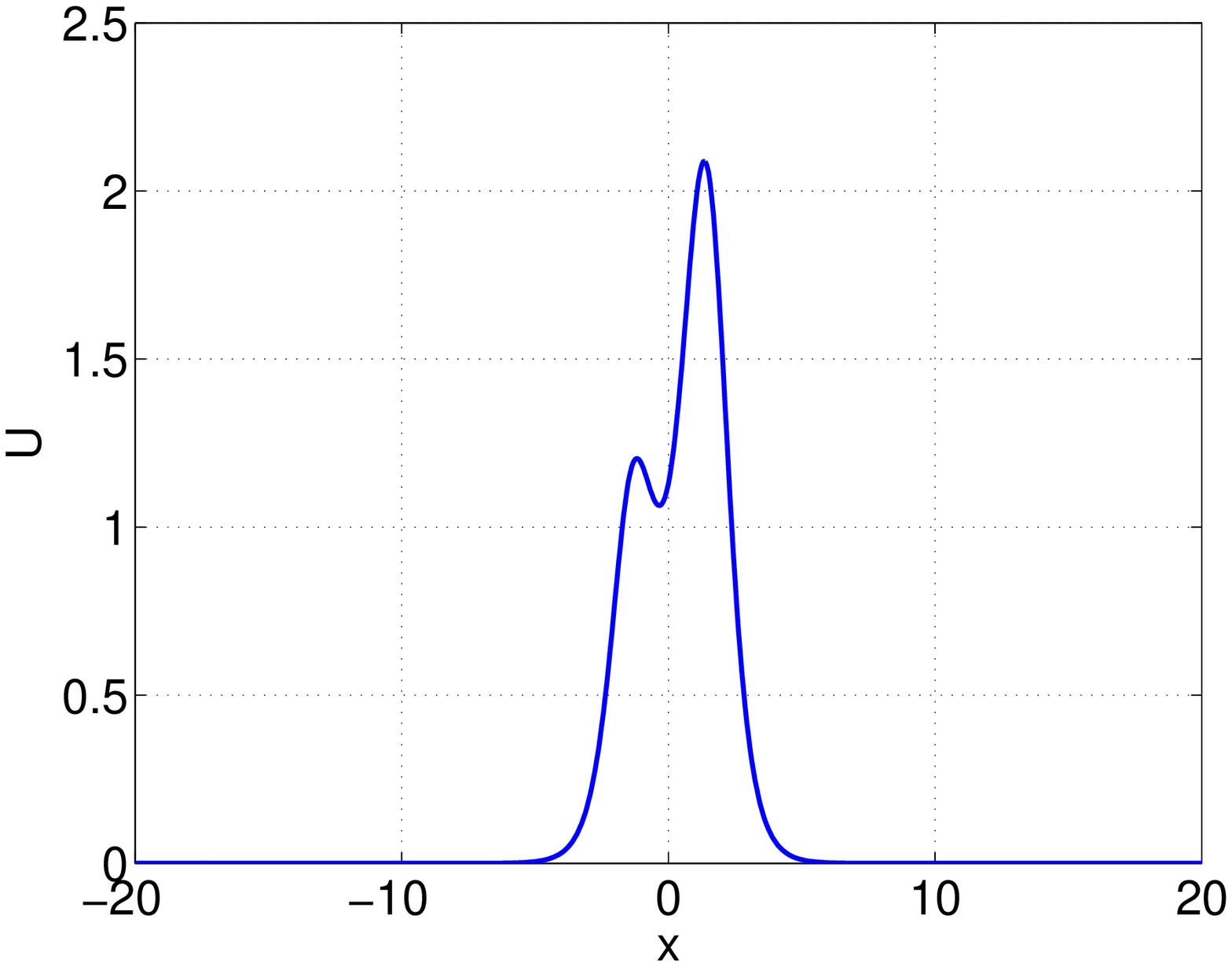} }
\caption{Numerical ground state generation of  (\ref{gnls2}) with $\mu=3.281$. (a) Approximate asymmetric profile; (a) $P=14.4462$; (b) $P=3.9918$; (c) $P=14.4559$; (d) $P=9.7217$.} \label{Figpet3}
\end{figure}

The application of acceleration with VEM to this example is displayed in Table \ref{tav_epet2} and Figure \ref{Figpet2}. The following points are emphasized:
\begin{itemize}
\item The acceleration leads to a great improvement with respect to the e-\PM method (\ref{lab22e}), (\ref{lab25e}). In order to have a residual error below $TOL$, the reduction in the number of iterations is above $90\%$.
\item As in the previous examples, polynomial methods work better than $\epsilon$-algorithms. By comparing the two polynomial methods, MPE is more efficient. Its best performance requires a large number of $\kappa$ (which means a long cycle, above eight). On the other hand, the best results for the $\epsilon$-algorithms are obtained with moderate values of $\kappa$, around five. (This also happens in general in the previous examples.)
\item The value of $\mu$ considered for the experiments is close to the one corresponding to the bifurcation point, that is it is close to the tangential point of the two branches of solitary wave solutions. The computation of the quantity (\ref{power2}) for each acceleration, shown in Table \ref{tav_epet2b}, attempts to study the behaviour of the iterations close to the bifurcation. In most of the cases, the computed value coincides to that of the profile generated by (\ref{lab22e}), (\ref{lab25e}) without acceleration. In the case of MPE(8) and TEA(2)-TEA(6), the value changes to $P_{h}(U_{h})=14.458882$. This suggests that for these cases the accelerated iteration converges to the profile of the upper branch  while in most of the cases (including the one without acceleration) the limit profile belongs to the lower branch,\cite{Yang2012}. (Close to the bifurcation indeed, the form of the profiles is very similar, see Figures \ref{Figpet3}(a) and (b). Note however from Table \ref{tav_epet1c} that the dominant eigenvalue of the iteration matrix of (\ref{lab22e}), (\ref{lab25e}) is above one. This and Table \ref{tav_epet1} may explain the convergence of this method to the profile with $P_{h}(U_{h})=14.4446162$. 
\item The comparison with the best choices of the VEM is illustrated in Figure \ref{Figpet2}, which compares the behaviour of the residual error (\ref{fsec32}) as function of the number of iterations and of CPU time in seconds. The results reveal again the better performance of the polynomial methods when the residual error starts to be below $10^{-5}$.
\end{itemize}
\begin{table}
\begin{center}
\begin{tabular}{|c|c|c|c|c|}
\hline\hline  $nw$&AA-I($nw$)&$P$&AA-II($nw$)&$P$\\\hline
$1$&$78$($5.2802E-11$)&$14.4589$&$59$($5.4557E-11$)&$14.4589$\\
$2$&$50$($7.8887E-11$)&$14.4462$&$51$($5.5612E-09$)&$14.4462$\\
$3$&$74$($4.4814E-11$)&$14.4589$&$28$($1.7879E-11$)&$3.9918$\\
$4$&$93$($9.0311E-11$)&$14.4462$&$32$($2.6548E-11$)&$3.9918$\\
$5$&$49$($2.0421E-11$)&$3.9918$&$64$($1.2763E-11$)&$3.9918$\\
$6$&$64$($3.3540E-11$)&$3.9918$&$80$($3.0285E-11$)&$3.9918$\\
$7$&$61$($4.5136E-11$)&$14.4462$&$29$($9.1724E-11$)&$3.9918$\\
$8$&$57$($4.7025E-11$)&$14.4662$&$54$($3.7450E-12$)&$9.7217$\\
$9$&&&$55$($5.6993E-12$)&$3.9918$\\
$10$&&&$108$($3.4559E-11$)&$3.9918$\\
\hline\hline
\end{tabular}
\end{center}
\caption{Ground state generation of  (\ref{gnls2}) with $\mu=3.281$.  Number of iterations required by AA-I and AA-II as function of $nw$ to achieve a residual error below $TOL=10^{-10}$. The residual error (\ref{fsec32}) at the last computed iterate is in parenthesis.}\label{tav_epet3}
\end{table}
When the iteration (\ref{lab22e}), (\ref{lab25e}) is accelerated with the AAM, we obtain the results displayed in Table \ref{tav_epet3}:
\begin{itemize}
\item The behaviour of the methods in this case looks similar to that of some previous examples as far as the general performance is concerned: they are competitive for moderate values of $nw$, with a better performance of AA-II, less affected by ill-conditioning.
\item In some cases the AMM approximate profiles (see Figures \ref{Figpet3}(c) and (d)) which correspond to values of (\ref{power2}) out of the branches. This uncertain behaviour provides the main drawback of the methods. The spectral information for the two additional approximate profiles is given in Tables \ref{tav_epet1b} and \ref{tav_epet1d}. The results suggest the lack of preservation of (\ref{power2}) through the iterative process.
\end{itemize}
\begin{table}
\begin{center}
\begin{tabular}{|c|c|}\hline\hline
Classical fixed point&e-\PM { method} (\ref{lab22e}), (\ref{lab25e})\\\hline
1.994420E+00&4.824994E-01\\
4.828929E-01&2.921950E-01-4.208258E-02i\\
2.155227E-01&2.921950E-01+4.208258E-02i\\
1.266822E-01&1.270267E-01\\
8.385958E-02&8.457380E-02\\
6.009441E-02&6.014597E-02\\\hline
\end{tabular}
\end{center}
\caption{Six largest magnitude eigenvalues of the iteration matrices
of classical fixed point algorithm and of e-\PM method (\ref{lab22e}), (\ref{lab25e}) for $\mu=3.281$ at the last computed iterate for $P=3.9918$.}\label{tav_epet1b}
\end{table}
\begin{table}
\begin{center}
\begin{tabular}{|c|c|}\hline\hline
Classical fixed point&e-\PM { method} (\ref{lab22e}), (\ref{lab25e})\\\hline
1.643665E+00&1.016836E+00\\
1.015912E+00&4.764502E-01\\
4.862159E-01&3.707178E-01\\
3.756950E-01&2.883210E-01-1.334802E-01i\\
2.022740E-01&2.883210E-01+1.334802E-01i\\
1.417784E-01&1.748132E-01\\\hline
\end{tabular}
\end{center}
\caption{Six largest magnitude eigenvalues of the iteration matrices
of classical fixed point algorithm and of e-\PM method (\ref{lab22e}), (\ref{lab25e}) for $\mu=3.281$ at the last computed iterate for $P=14.4559$.}\label{tav_epet1c}
\end{table}
\begin{table}
\begin{center}
\begin{tabular}{|c|c|}\hline\hline
Classical fixed point&e-\PM { method} (\ref{lab22e}), (\ref{lab25e})\\\hline
1.919934E+00&1.418879E+00\\
1.201956E+00&6.091959E-01-4.750746E-02i\\
6.134143E-01&6.091959E-01+4.750746E-02i\\
3.716346E-01&4.517846E-01\\
2.394216E-01&2.884719E-01\\
1.733253E-01&1.799915E-01\\\hline
\end{tabular}
\end{center}
\caption{Six largest magnitude eigenvalues of the iteration matrices
of classical fixed point algorithm and of e-\PM method (\ref{lab22e}), (\ref{lab25e}) for $\mu=3.281$ at the last computed iterate for $P=9.7217$.}\label{tav_epet1d}
\end{table}
\section{Concluding remarks and future work}
\label{se5}
In this paper we have studied numerically the use of acceleration techniques applied to fixed point algorithms of \PM type to generate numerically traveling waves in nonlinear dispersive wave equations. The comparison has been established between vector extrapolation methods and Anderson acceleration methods for different types of traveling waves. From the plethora of numerical experiments, our main conclusions are:
\begin{itemize}
\item The use of acceleration techniques improves the performance of the \PM type methods in all the cases. This improvement is observed in two main points: first, when the \PM type method is convergent, the acceleration reduces the number of iterations in a relevant way. (In some cases, this is really important: in some one-dimensional problems, the reduction is at least of $50\%$ and attains up to $75\%$.) On the other hand, the mechanism of acceleration, especially in the case of VEM, allows to transform initially divergent sequences into convergent processes. This is particularly relevant in traveling waves with high oscillations. Furthermore, acceleration has been shown to be more efficient than other alternatives for some cases, like numerical continuation.
\item In general, VEM provide better results and among them, polynomial methods such as MPE and RRE are more efficient than $\epsilon$-algorithms like VEA and TEA, although in some convergent cases the acceleration in terms of the number of iterations is very similar among all the methods while in computational time the $\epsilon$-algorithms work better.
\item The AAM are competitive in some cases but they are mostly affected by ill-conditioning and a more computational effort due to their longer implementation. The best results of these methods are obtained when generating numerically periodic traveling waves in some nonlinear dispersive systems and ground state profiles in NLS type equations, while their performance is poor when computing highly oscillatory traveling waves. Their application to the e-\PM type methods suggests an uncertain behaviour with respect to relevant quantities of the problem through the iteration.
\end{itemize}
The main question to that this comparative study has not been able to answer is in the authors' opinion finding a deeper understanding of the way how acceleration techniques (especially VEM) work on these problems. In particular, we miss some conclusions about the width of extrapolation (which is related to the extrapolation step for convergence) to be used a priori (if it is possible to do that). We have observed that this looks to be strongly dependent on the problem under study. This might be though a good starting point for a future research.

\section*{Acknowledgements}
This research has been supported by  project
MTM2014-54710-P.

\end{document}